\documentclass[10pt,leqno]{article}

\usepackage[francais]{babel}
\usepackage{amsmath,amssymb,amscd}
\usepackage{a4wide}

\usepackage[latin1]{inputenc}

\usepackage[all]{xy} 

\makeatletter
\@addtoreset{equation}{subsection}         
\makeatother

\newenvironment{paragr}[1][]{\refstepcounter{subsection} \noindent \textbf{\thesubsection . \ #1}}{\medskip}

\newenvironment{theoreme}{ \medskip\refstepcounter{theo}  \noindent\textbf{Th\'eor\`eme \thetheo}. ---\em}{\em \medskip}
\newenvironment{proposition}{\medskip\refstepcounter{theo}   \noindent\textbf{Proposition \thetheo}. ---\em}{\em\medskip}
\newenvironment{corollaire}{\medskip\refstepcounter{theo}  \noindent\textbf{Corollaire \thetheo}. ---\em}{\em\medskip}

\newenvironment{lemme}{\medskip\refstepcounter{theo}   \noindent\textbf{Lemme \thetheo}. ---\em}{\em\medskip}

\newenvironment{definition}{\medskip\refstepcounter{theo}  \noindent\textbf{D\'efinition \thetheo}. ---}{\medskip}

\newenvironment{preuve}[1][]{\noindent \textbf{Démonstration.} #1 --- }{\hfill
  \ensuremath{\square} \medskip}

\newenvironment{remarque}{\medskip\refstepcounter{theo}  \noindent\textbf{Remarque \thetheo}. ---}{\medskip}
\newenvironment{remarques}{\medskip\refstepcounter{theo}  \noindent\textbf{Remarques \thetheo}. ---}{\medskip}

\DeclareMathOperator{\reg}{reg}

\DeclareMathOperator{\vol}{vol}
\DeclareMathOperator{\rang}{rang}
\DeclareMathOperator{\Ad}{Ad}

\DeclareMathOperator{\End}{End}

\DeclareMathOperator{\Norm}{Norm}

\DeclareMathOperator{\Aut}{Aut}

\DeclareMathOperator{\Hom}{Hom}

\DeclareMathOperator{\inv}{inv}

\DeclareMathOperator{\Id}{Id}
\DeclareMathOperator{\stab}{stab}

\DeclareMathOperator{\Ker}{Ker}
\DeclareMathOperator{\Ima}{Im}

\newcommand{\ZZ}{\mathbb{Z}}

\newcommand{\SG}{\mathfrak{S}}
\newcommand{\NN}{\mathbb{N}}
\newcommand{\RR}{\mathbb{R}}
\newcommand{\AAA}{\mathbb{A}}
\newcommand{\CC}{\mathbb{C}}

\newcommand{\QQ}{\mathbb{Q}}

\newcommand{\oc}{\mathcal{O}}
\newcommand{\rc}{\mathcal{R}}
\newcommand{\Sc}{\mathcal{S}}
\newcommand{\tc}{\mathcal{T}}

\newcommand{\vc}{\mathcal{V}}

\newcommand{\ec}{\mathcal{E}}
\newcommand{\gc}{\mathcal{G}}

\newcommand{\yc}{\mathcal{Y}}
\newcommand{\lc}{\mathcal{L}}

\newcommand{\fc}{\mathcal{F}}
\newcommand{\pc}{\mathcal{P}}

\newcommand{\bc}{\mathcal{B}}

\newcommand{\ggo}{\mathfrak{g}}

\newcommand{\of}{\mathfrak{o}}

\newcommand{\mgo}{\mathfrak{m}}
\newcommand{\ngo}{\mathfrak{n}}
\newcommand{\ago}{\mathfrak{a}}

\newcommand{\pgo}{\mathfrak{p}}
\newcommand{\qgo}{\mathfrak{q}}

\newcommand{\rgo}{\mathfrak{r}}
\newcommand{\hgo}{\mathfrak{h}}

\newcommand{\lgo}{\mathfrak{l}}

\newcommand{\al}{\alpha}
\newcommand{\be}{\beta}

\newcommand{\la}{\lambda}

\newcommand{\back}{\backslash}

\newcommand{\Cc}{C_c^\infty}

\newcommand{\bg}{\langle}
\newcommand{\bd}{\rangle}

\newcommand{\eps}{\varepsilon}

\renewcommand{\leq}{\leqslant}
\renewcommand{\geq}{\geqslant}

\title{Sur la contribution unipotente dans la formule des traces d'Arthur pour les groupes généraux linéaires}
\author{Pierre-Henri Chaudouard}
\date{}
\begin{document}

\maketitle

\begin{abstract}

 Le thème de l'article est  l'étude de la partie unipotente de la formule des traces d'Arthur pour un groupe général linéaire. La contribution de l'orbite unipotente régulière ou de ses variantes par blocs  a été essentiellement traitée dans un précédent article (cf. \cite{scfhn}). Ici on s'intéresse à la contribution des orbites unipotentes simples c'est-à-dire aux orbites de Richardson induites à partir d'un sous-groupe de Levi dont les blocs sont de tailles deux-à-deux distinctes. De manière remarquable, la contribution s'exprime à l'aide d'une intégrale orbitale pondérée globale unipotente. Comme corollaire, on obtient des formules intégrales pour certains coefficients globaux d'Arthur. Cet article comprend également une construction nouvelle des intégrales orbitales pondérées unipotentes locales d'Arthur.
  
\end{abstract}

\renewcommand{\abstractname}{Abstract}
\begin{abstract}
The theme of the article is  the study of the unipotent part of Arthur's trace formula for general linear groups. The case  of regular (or ''regular by blocks'') unipotent orbits has been essentially done in a previous paper (cf. \cite{scfhn}). Here we are interested by the contribution of Richardson orbits that are induced by Levi subgroups with two-by-two distinct blocks. In this case, the contribution is remarkably given by a global  unipotent weighted orbital integral. As a corollary, we get integral formulas for some of Arthur's global coefficients. We also present a new construction of Arthur's local unipotent weighted orbital integral.
\end{abstract}

\tableofcontents

\section{Introduction}

\begin{paragr}
  La formule des traces d'Arthur-Selberg est un outil central en théorie des formes automorphes et en fonctorialité de Langands. Certaines de ses plus spectaculaires applications reposent  sur une comparaison de plusieurs formules des traces. Ce sont celles qui ont trait à la partie dite endoscopique du programme de Langlands. D'autres reposent sur l'utilisation d'une formule des traces isolée, par exemple le calcul de dimension de certains espaces de formes automorphes (cf. \cite{Lan-dim}), les lois de Weyl avec reste (par exemple cf. \cite{Muller} ou \cite{Lapid-Muller}). Pour ces dernières, il est intéressant de savoir si la formule des traces vaut encore pour des fonctions test qui ne sont pas à support compact (cf. les résultats de \cite{Hoffmann} et dans une autre direction \cite{FL-GL2}, \cite{FLss} et  \cite{FLM}). Par ailleurs, Langlands a suggéré d'introduire dans la formule des traces des fonctions qui ne sont pas à support compact mais qui permettraient de discriminer le spectre automorphe suivant les pôles de ses fonctions $L$ associées (cf. \cite{Lan-beyond}). Cela soulève alors des problèmes de convergence  du côté géométrique dans la formule des traces. Ces problèmes sont d'autant plus difficiles que la partie géométrique non-semi-simple de la formule des traces d'Arthur n'est pas sous une forme aussi aboutie que ce que l'on pourrait désirer (même si elle convient tout-à-fait au programme endoscopique). En effet, dans ses travaux (cf. \cite{dist_orb}), Arthur montre que la partie géométrique s'exprime comme une certaine combinaison linéaire d'intégrales orbitales pondérées semi-locales. Les coefficients  qui apparaissent sont de nature plus globale et seule leur existence est affirmée : ils ne sont jamais explicités sauf pour des orbites semi-simples.  Il y a en fait une procédure de descente au centralisateur de la partie semi-simple ; le point crucial est donc de comprendre ces coefficients dans le cas unipotent. Mentionnons les  articles récents \cite{Matz-bounds} et \cite{HW} où l'étude de ces coefficients est abordée. 
\end{paragr}

\begin{paragr}
Dans cet article, on considère les groupes généraux linéaires et on explicite la contribution de certaines orbites unipotentes en termes de nouveaux objets : des intégrales orbitales pondérées unipotentes globales. Comme corollaire, on obtient des formules intégrales pour les coefficients globaux associés à certaines orbites unipotentes. Au passage, nous donnons  une construction très simple des intégrales orbitales pondérées unipotentes locales d'Arthur. 
\end{paragr}

\begin{paragr}[Intégrales orbitales pondérées unipotentes globales.] --- Présentons maintenant plus en détail nos résultats. Soit $G=GL(n)$ sur un corps de nombres $F$ dont on note $\AAA$ l'anneau des adèles. Dans le reste de l'article, on travaille en fait sur la partie nilpotente de la variante de la formule des traces sur l'algèbre de Lie $\ggo$ de $G$  (cf. \cite{PH1}). On passe aisément pour les groupes linéaires généraux de cette partie nilpotente à la partie unipotente de la formule des traces puisque pour toute fonction $f\in \Cc(G(\AAA))$, il existe $\phi\in \Cc(\ggo(\AAA))$ telle que pour $g\in G(\AAA)$ et tout $Y\in \ggo(\AAA)$ nilpotent, on ait 
$$f(g^{-1}(I+Y)g)=\phi(g^{-1}Yg).$$

Pour les besoins de l'introduction, on donne tous les énoncés dans leur version sur $G$.  Soit $\of$ une orbite unipotente. Pour les groupes généraux linéaires, toutes les orbites sont de Richardson : il existe donc $P$ un sous-groupe parabolique de $G$ tel que $\of\cap N_P$ soit ouvert et dense dans $N_P$ où $N_P$ est le radical unipotent de $P$. Soit $M$ un facteur de Levi de $P$. Un tel $P$ n'est pas unique mais la classe de conjugaison de $M$ l'est. Quitte à conjuguer, on peut et on va  supposer que $M$ est le sous-groupe de Levi standard $GL(n_1)\times\ldots\times GL(n_r)$ avec $n=n_1+\ldots+n_r$ et $n_1\geq n_2 \geq \ldots \geq n_r$.

Soit 
$$\gamma=  \begin{pmatrix} I_{n_1} & \underset{0}{I_{n_2}} & & &   \\   &I_{n_2}  &  \underset{0}{I_{n_3}} & & \\ & & \ddots & \ddots& \\ & & & I_{n_{r-1}} & \underset{0}{I_{n_r}} \\ & & & & I_{n_r}\end{pmatrix}  \in GL(n,F) $$
où $I_k$ désigne la matrice identité carrée de taille $k$. Alors $\gamma \in\of$. Soit $W\subset G$ le sous-groupe fini des matrices de permutation. Soit $\pc(M)$ l'ensemble fini des sous-groupes paraboliques de $G$ de facteur de Levi $M$. Soit $P\in \pc(M)$.
 On dispose pour tout $P\in \pc(M)$ de l'application de Harish-Chandra (relative au sous-compact maximal standard $K$ de $G(\AAA)$)
$$H_P: G(\AAA) \to \Hom(X^*(M),\RR)$$
où $X^*(M)$ est le groupe des caractères de $M$. On pose alors pour tout $g\in G(\AAA)$
$$R_P(g)=H_P(wg)$$
où $w\in W$ est un élément qui vérifie $\gamma \in w^{-1}N_P w$. Un tel élément est bien défini à une translation à droite près par $W^M=W\cap M$ de sorte que $R_P(g)$ est bien défini. La famille 
$$(\exp(-\bg \la, R_P(g)\bd))_{P\in \pc(M)}$$
de fonctions de la variable $\la \in X^*(M)\otimes_{\ZZ}\CC$ est en fait une $(G,M)$-famille au sens d'Arthur. Par un procédé standard d'Arthur, on sait lui associer un poids (pas forcément positif) noté
$$v_{L,\gamma}(g)$$
pour tout sous-groupe de Levi $L$ de $G$ contenant $M$. On montre que le poids $v_{L,\gamma}^G$ vu comme fonction sur $G(\AAA)$ est invariant par $G_\gamma(\AAA)$ où $G_\gamma$ est le centralisateur de $\gamma$ dans $G$.
Le premier théorème que l'on démontre est le suivant.

\begin{theoreme} (cf. théorème \ref{thm:cv-pondere})
Supposons $n_1>n_2>\ldots >n_r$. 
  Pour toute fonction $f\in \Cc(G(\AAA))$ l'intégrale orbitale pondérée globale
$$J_{M,\gamma}(f)=\int_{G_\gamma(\AAA)\back G(\AAA)} f(g^{-1}\gamma g)\, v_{M,\gamma}(g) \, dg$$
est absolument convergente.
\end{theoreme}
 
S'il existe $i<r$ tel que $n_i=n_{i+1}$, l'intégrale ne converge pas en général. 

\end{paragr}

\begin{paragr}[Contribution de l'orbite $\of$ dans la formule des traces.] --- Soit $f\in \Cc(G(\AAA))$. Pour définir la contribution $J_\of(f)$ de l'orbite $\of$ dans la formule des traces pour $G$, Arthur commence par introduire l'intégrale orbitale tronquée (donc convergente)
$$\int_{G_\gamma(F)\back G(\AAA)^1} F(g,T) f(g^{-1}\gamma g) \, dg$$
où $F(\cdot,T)$ est la fonction caractéristique d'un compact de $G(F)\back G(\AAA)^1$ construit à l'aide de la théorie de la réduction et dépendant d'un paramètre $T$. Pour tout groupe $H$, on note 
$$H(\AAA)^1=\bigcap_{\chi\in X^*(H)} \ker(|\chi|_\AAA).$$
Arthur montre que cette intégrale orbitale tronquée est asymptotique à un polynôme en $T$ quand le paramètre \og tend vers l'infini\fg{} dans une certaine chambre. Par définition, $J_\of(f)$ est le terme constant de ce polynôme.  Notre second résultat est le suivant.

\begin{theoreme}\label{thm:Jo=JM}(cf. théorème \ref{thm:asymp} et proposition \ref{thm:nilp-globale})
Supposons $n_1>n_2>\ldots >n_r$. 
  Pour toute fonction $f\in \Cc(G(\AAA))$, on a 
$$ J_\of(f) =\vol(G_\gamma(F)\back G_\gamma(\AAA)^1)\cdot J_{M,\gamma}(f)$$
\end{theoreme}

Le théorème ne vaut pas s'il existe $1\leq i <r-1$ tel que $n_i=n_{i+1}$ ne serait-ce que parce que dans ce cas l'intégrale $J_{M,\gamma}(f)$ n'est pas convergente. Dans le cas extrême $n_1=n_2=\ldots=n_r$, on peut obtenir une expression pour $J_\of(f)$ qui  fait intervenir des intégrales orbitales régularisées (le cas régulier $r=n$ et $n_1=\ldots=n_r=1$ a été obtenu par Finis et Lapid dans un travail non publié ; le cas d'un corps de fonctions $F$ et d'une fonction test $f$ très simple a été traité  dans \cite{scfhn}). 
\end{paragr}

\begin{paragr}[Développement d'Arthur.] --- On déduit du théorème le corollaire suivant. Soit $S$ un ensemble fini de places contenant les places archimédiennes. On décompose $\AAA=\AAA_S\times \AAA^S$, l'exposant $S$ désignant la partie \og hors $S$\fg{}. Le sous-groupe compact maximal standard $K$ de $G(\AAA)$ se décompose aussi en $K_S\times K^S$. Soit $ \mathbf{1}^S$ la fonction caractéristique de $K^S$.

  \begin{corollaire}\label{cor:dvptA} (cf. théorème \ref{thm:dvpt})
    On suppose $n_1>\ldots>n_r$. Pour toute fonction $f_S\in \Cc(G(\AAA_S))$, on a
$$
  J_{\of}(f_S\otimes \mathbf{1}^S)=\sum_{(L,\of_L)} \frac{|W^L|}{|W|}\,  a^L(S,\of_L) \, J_L(\of_L, f_S),
$$
où
\begin{itemize}
\item la somme porte sur les couples $(L,\of_L)$ formés d'un sous-groupe de Levi semi-standard $L$ de $G$ et d'une orbite unipotente $\of_L$ dont l'induite de Lusztig-Spaltenstein est $\of$ ;
\item $J_L(\of_L, f_S)$ est une intégrale orbitale unipotente semi-locale d'Arthur définie dans \cite{ar_unipvar} ;
\item le coefficient est défini par
$$ a^L(S,\of_L)=\vol(L_1(F)\back L_1(\AAA)^1) \cdot J_{L_1}^{L}(1,\mathbf{1}^S_{L})$$
où $L_1$ est un sous-groupe de Levi inclus dans $L$ dont l'orbite de Richardson est $\of_L$ ; la distribution  $J_{L_1}^{L}(1,\mathbf{1}^S_{L})$ vaut $1$ si $L_1=L$ et si $L_1\subsetneq L$ c'est un nouvel objet : une intégrale orbitale pondérée globale hors $S$ sur $L$ relative  à l'orbite triviale $1$ de $L$ et la fonction unité $\mathbf{1}^S_{L}$ sur $L(\AAA)$.
\end{itemize}
  \end{corollaire}

Par exemple, la contribution associée à $(G,\of)$ est $a^G(S,\of) J_G(\of, f_S)$ où l'on a pour tout  $P\in \pc(M)$
$$ J_G(\of, f_S)=\int_{K_S} \int_{N_P(\AAA_S)} f_S(k^{-1}uk)\, du$$
et
$$a^G(S,\of)= \vol(M(F)\back M(\AAA)^1) \cdot   \int_{N_P(\AAA^S)} \mathbf{1}^S(u ) \,  w_M(u) \, du$$ 
où la dernière intégrale est convergente relativement à un certain poids $w_M$. 

Le corollaire \ref{cor:dvptA} se déduit du théorème \ref{thm:Jo=JM} au moyen de la combinatoire des $(G,M)$-familles et de l'énoncé suivant. Pour toute place $v$ de $F$, soit $F_v$ le complété de $F$ en $v$. 

\begin{theoreme}\label{thm:IOversusAr}(cf. théorème \ref{thm:egalite-IOP})
  Il existe une constante $c_{\gamma,v}$ qui dépend des choix de mesures telles que pour tout sous-groupe de Levi $L$ contenant $M$, on a 
$$\int_{G_\gamma(F_v)\back G(F_v)} f_v(g^{-1}\gamma g) v_{L,\gamma}(g)\, dg = c_{\gamma,v} \cdot J_L( \of^L,f_v)$$
où
\begin{itemize}
\item $\of^L$ est l'induite de Richardson de l'orbite unipotente triviale de $M$ ;
\item l'intégrale de gauche est absolument convergente ;
\item dans le membre de droite, $J_L( \of^L,f_v)$ est une intégrale orbitale pondérée unipotente d'Arthur.
\end{itemize}
\end{theoreme}

Pour des choix \og naturels\fg{} de mesure, la constante ne vaut pas $1$ mais admet une expression intéressante en termes de facteurs locaux de la fonction zêta du corps de base $F$. Remarquons que le théorème donne une façon rapide de définir les intégrales orbitales unipotentes pondérées d'Arthur. Par ailleurs, on déduit du théorème \ref{thm:IOversusAr} des majorations pour les coefficients $a^L(S,\of_L)$ (cf. corollaire \ref{cor:calcul}). Comme autre application du corollaire  \ref{cor:dvptA}, dans le cas local, on obtient (cf. théorème \ref{thm:calcul-IOPlocale}) des formules pour certaines intégrales orbitales unipotentes pondérées  d'une fonction test très simple (unité de l'algèbre de Hecke sphérique dans le cas local non-archimédien).
\end{paragr}

\begin{paragr}[Plan de l'article.] --- Les sections \ref{sec:Richard} à \ref{sec:fam-ortho} sont des préliminaires qui  sont consacrés à la combinatoire de certains sous-groupes paraboliques (dits de Richardson) associés à un élément nilpotent de $G$. On y donne aussi la définition des fonctions $R_P$ et leurs principales propriétés. On déduit de ces fonctions $R_P$ la construction de $(G,M)$-familles à la section \ref{sec:GMfam} qu'on compare à une autre construction due à Arthur. Dans les sections \ref{sec:IOloc} à \ref{sec:uncalculIOPloc} qui sont purement locales, on introduit une nouvelle construction des poids et des intégrales orbitales pondérées locales. On compare aussi nos intégrales à celles d'Arthur. Le reste de l'article est dans un cadre global et se limite aux orbites nilpotentes pour lesquelles on a la condition $n_1>\ldots>n_r$ (cf. plus haut) qu'on appelle simples. Dans la section \ref{sec:IOPglob}, on prouve la convergence des intégrales orbitales globales qui leur sont associées et qui dépendent du paramètre de troncature $T$. Dans la section \ref{sec:asymp}, on montre que ces intégrales sont asymptotiques en $T$ aux intégrales orbitales tronquées d'Arthur. On peut alors identifier la contribution d'une telle orbite dans la formule des traces d'Arthur à une intégrale orbitale pondérée. Dans la section finale \ref{sec:dvpt-Arthur}, on énonce quelques corollaires comme une version précisée du développement d'Arthur.
\end{paragr}

\begin{paragr}[Remerciements.] --- Lors de l'élaboration de cet article, j'ai reçu le soutien de l'Institut Universitaire de France et des  projets Ferplay ANR-13-BS01-0012 et  Vargen ANR-13-BS01-0001-01 de l'ANR. 
J'ai terminé cet article lors d'un séjour très agréable à l'Institute for Advanced Study à Princeton, durant lequel j'ai bénéficié du soutien de la National Science Foundation (\emph{agreement} No. DMS-1128155). Bien sûr les résultats du présent article n'engagent que leur auteur et ne reflètent pas nécessairement les opinions de la  National Science Foundation.
  
C'est donc un plaisir pour moi de remercier toutes ces institutions. Je remercie enfin Gérard Laumon pour les nombreuses conversations que nous avons pu avoir autour des thèmes abordés dans cet article ou dans un précédent article en commun (\cite{scfhn}).
\end{paragr}

\section{Sous-groupes de Richardson d'un élément nilpotent}\label{sec:Richard}

\begin{paragr}Soit $F$ un corps qui dans tout le texte sera de caractéristique $0$. Pour tout groupe algébrique défini sur $F$ noté $G$, on note par la même lettre en minuscule gothique, ici $\ggo$, son algèbre de Lie. Soit $N_G$ son radical unipotent et $A_G$ son centre. Soit $X^*(G)$ le groupe des caractères rationnels de $G$ défini sur $F$. L'action adjointe de $G$ sur $\ggo$ est noté $\Ad$. Sauf mention contraire, un sous-groupe de $G$ signifie un sous-groupe algébrique de $G$ défini sur $F$. Pour tous sous-groupes $H$ et $M$ de  $G$, soit $\Norm_H(M)$ le normalisateur de $M$ dans $H$. 
 \end{paragr}

\begin{paragr}\label{S:fM}
  Soit $E$ un $F$-espace vectoriel. Soit $G$ le groupe algébrique défini sur $F$ des automorphismes de $E$. Son algèbre de Lie s'identifie à l'algèbre des endomorphismes de $E$.
\end{paragr}

\begin{paragr}[Notations d'Arthur.] ---  Tout au long de l'article, on s'efforcera de suivre les notations qui se sont imposées dans les travaux d'Arthur. Pour la commodité du lecteur, on les rappelle brièvement ici. Soit $M\subset H$ des sous-groupes de $G$. Soit $\fc^H(M)$ l'ensemble des sous-groupes paraboliques de $G$ inclus dans $H$ et contenant $M$. Supposons de plus que  $M$ est un sous-groupe de Levi de $G$, ce par quoi on entend un facteur de Levi d'un sous-groupe parabolique de $G$. Soit  $\pc^H(M)$ l'ensemble des sous-groupes paraboliques de $G$ inclus dans $H$, dont $M$ est un facteur de Levi.  Tout sous-groupe parabolique $Q\in \fc^H(M)$, possède un unique facteur de Levi noté $M_Q$ qui contient $M$. En particulier, $Q=M_Q N_Q$ est une décomposition de Levi de $Q$. L'image de l'application $Q\in\fc^H(M) \mapsto M_Q$ est notée $\lc^H(M)$. Lorqu'on a  $H=G$, on omet l'exposant $G$ dans la notation.

 Soit $M$ un sous-groupe de Levi de $G$. Pour tout sous-groupe $H$ de $G$ stable par conjugaison par $A_M$, on note $\Sigma(H,A_M)$ l'ensemble des racines de $A_M$ sur $\hgo$. Pour $\al \in  \Sigma(H,A_M)$, on note encore $\al$ la forme linéaire sur l'algèbre de Lie $\ago_M$ de $A_M$ obtenue par différentiation. Soit  $P\in \pc(M)$ et $Q$ un sous-groupe parabolique contenant $P$.  Soit $\Delta_P^Q$ l'ensemble des racines simples dans $\Sigma(M_Q\cap N_P,A_M)$. Lorsque $Q=G$, on note simplement $\Delta_P=\Delta_P^G$.   On a $\Delta_P^Q\subset \Delta_P$.

Pour tout $P\in \pc(M)$, soit $a_P=\Hom(X^*(P),\RR)$ et son dual $a_P^*=X^*(P)\otimes_\ZZ\RR$. Soit $\bg \cdot, \cdot \bd$ l'accouplement canonique entre $a_P^*$ et $a_P$. Le morphisme de  restriction $X^*(P)\to X^*(M)$ est bijectif ce qui permet d'identifier $a_P$ à $a_M=\Hom(X^*(M),\RR)$, de même pour les espaces duaux. On note $a_{P,\CC}$, $a_{P,\CC}^*$ etc. les espaces complexes obtenus par extension des scalaires.

Soit $M\subset L$ des sous-groupes de Levi et $P\in \pc(M)$ et $Q\in \pc(L)$ tels que $P\subset Q$. En utilisant le morphisme de  restriction $X^*(M)\to X^*(A_M)$, on identifie $a_M^*$ à $X^*(A_M)\otimes \RR$. De la sorte on a $\Sigma(P,A_M)\subset a_M^*$. Soit $a_P^{Q,*}$ le sous-espace de $a_P^*$ engendré par $\Delta_P^Q$ ; il ne dépend que de $M$ et $L$, on le note encore $a_M^{L,*}$. On note $\rho_P^Q$ la demi-somme des éléments de  $\Sigma(M_Q\cap N_P,A_M)$. À chaque $\al\in \Delta_P$ est associée une coracine $\al^\vee\in a_P$ (cf. \cite{ar1}). On note $\Delta_P^{Q,\vee}$ l'ensemble des coracines des éléments de $\Delta_P^Q$. Dualement au morphisme de  restriction $X^*(L)\to X^*(M)$, on a une projection  $a_{M} \to a_{L}$ dont on note $a_M^{L}$ le noyau. L'ensemble $\Delta_P^{Q,\vee}$ forme une base de $a_M^{L}$. Ce dernier est en dualité parfaite avec $a_M^{L,*}$.  Soit $\hat{\Delta}_P^Q$ la base des poids de  $a_M^{L,*}$, c'est-à-dire la base duale de $\Delta_P^{Q,\vee}$. Soit  $\hat{\Delta}_P^{Q,\vee}$ la base des copoids, c'est-à-dire la base de  $a_M^{L}$ duale de  $\Delta_P^{Q}$. Pour tout $\Delta\subset a_M$, on note $\ZZ(\Delta)$ le sous-groupe engendré par $\Delta$.
\end{paragr}

\begin{paragr}
  Soit $X\in \ggo$ un élément nilpotent.  Soit $G_X$ le centralisateur de $X$ dans $G$. C'est un sous-groupe connexe. 
\end{paragr}

\begin{paragr}[Induite de Lusztig-Spaltenstein.] --- \label{S:ind-LS}Soit $P\subset G$ un sous-groupe parabolique tel que son algèbre de Lie $\pgo$, contienne $X$. Soit $L=P/N_P$ le plus grand quotient réductif de $P$. Le groupe $L$ agit par adjonction sur son algèbre de Lie $\lgo=\pgo/\ngo_P$. Soit 
$$\pi: \pgo\to \lgo$$
 la projection canonique. 
Soit $\oc^L_X$ l'orbite de $\pi(X)$ sous $L$. 

Soit $I_P^G(X)$ (qu'on notera encore  $I_P^G(\oc^L_X)$) l'unique orbite nilpotente dans $\ggo$ pour l'action adjointe de $G$ telle que l'intersection
\begin{equation}
  \label{eq:intersection}
  I_P^G(X)\cap \pi^{-1}(\oc^L_X)
\end{equation}
soit un ouvert de Zariski  dense dans $\oc^L_X+\ngo_P$. C'est \emph{l'induite de Lusztig-Spaltenstein de $X$} (cf.  \cite{lus-spal}). L'intersection \eqref{eq:intersection}, lorsqu'elle n'est pas vide,  est exactement la $P$-orbite de $X$ (\emph{ibid.}).

Pour tout sous-groupe parabolique $P\subset Q\subset G$, on peut définir de la même façon une orbite induite $I_P^Q(X)$ qui est une orbite nilpotente dans $\qgo/\ngo_Q$ (pour l'action adjointe de $Q/N_Q$). L'induction est transitive au sens où l'on a la formule suivante (\emph{ibid.}, §1.7)
\begin{equation}
  \label{eq:transitivite}
  I_Q^G(I_P^Q(X))=I_P^G(X).
\end{equation}
 
Lorsque l'orbite $\oc_X^L$ est l'orbite nulle $(0)$, l'induite $I_P^G(0)$ est \emph{l'orbite de Richardson} de $P$.
\end{paragr}

\begin{paragr}[Invariants d'un endomorphisme nilpotent $X$.] ---  Soit $r$ l'indice de nilpotence de $X$. Pour $1\leq j \leq r$ soit $d_j$ la multiplicité du bloc de taille $j$ dans la décomposition de Jordan de $X$. Soit
  \begin{equation}
    \label{eq:invX}
    \inv(X)= |\{1\leq j \leq r  \mid d_j\not=0   \}  |.
  \end{equation}
  C'est un entier compris entre $1$ et $j$.

\begin{lemme}(\cite{Spal} II. proposition 5.14)\label{lem:Richard2}
  Soit $P$ un sous-groupe parabolique de $G$ et $X\in I_P^G(0)$.  Alors $P$ est le stabilisateur d'un drapeau de sous-$F$-espaces
$$(0)\subsetneq E_1 \subsetneq \ldots  \subsetneq E_r=E$$
dont la suite ordonnée des dimensions des quotients $E_{i}/E_{i-1}$ pour $1\leq i\leq r$ est 
\begin{equation}
  \label{eq:ordonnee}
  d_r\leq d_r+d_{r-1}\leq \ldots \leq d_r+\ldots +d_1.
\end{equation}
En particulier deux sous-groupes paraboliques $P$ et $Q$ satisfont $I_P^G(0)=I_Q^G(0)$ si et seulement s'ils sont associés.
\end{lemme}

\begin{lemme}\label{lem:proj-induite}
  Soit $P$ un sous-groupe parabolique. Soit $\oc=I_P^G(0)$ et $X\in \oc\cap \ngo_P$. Pour tout sous-groupe parabolique $Q$ contenant $P$, la projection de $X$ sur $\qgo/\ngo_Q$ appartient à $I_P^Q(0)$.
\end{lemme}

\begin{preuve}
L'intersection $\oc\cap \ngo_P$ est ouverte et dense dans $\ngo_P$. De plus, c'est une $P$-orbite. Sa projection sur $\ngo_P/\ngo_Q$ est un ouvert dense de $\ngo_P/\ngo_Q$ qui contient la projection $Y$ de $X$ et elle est incluse dans la $Q/N_Q$-orbite de $Y$. Donc $Y\in I_P^Q(0)$.
\end{preuve}
\end{paragr}

\begin{paragr}[Ensembles $\rc(X)$ et $\lc\Sc(X)$.] --- Soit
  \begin{equation}
    \label{eq:sgp-LS}
    \lc\Sc(X)=\{P\subset G \textrm{ sous-groupe parabolique} \mid X\in \pgo \textrm{ et }X\in I_P^G(X)  \}
  \end{equation}

  \begin{equation}
    \label{eq:sgp-Ric}
     \rc(X)=\{P\subset G \textrm{ sous-groupe parabolique} \mid X\in \ngo_P \textrm{ et }X\in I_P^G(0)  \}.
  \end{equation}

Les éléments de  $\lc\Sc(X)$ et  $\rc(X)$ sont appelés respectivement les sous-groupes de Lusztig-Spaltenstein et de Richardson de $X$. Ce sont des ensembles non vides, en particulier $\rc(X)$ contient toujours le drapeau des images itérées et celui des noyaux itérés (comme il résulte du  lemme \ref{lem:Richard2}). Il résulte de la transitivité de l'induction qu'un sous-groupe parabolique $Q$ appartient à $\lc\Sc(X)$ si et seulement s'il contient un élément de $\rc(X)$. En particulier, on a $\rc(X)\subset \lc\Sc(X)$. Le but de cette section est de fournir une première description de ces ensembles. 

\begin{lemme}(cf. \cite{lus-spal}, théorème 1.3) \label{rq:centralisateur}
 Pour tout $P\in \lc\Sc(X)$, sous-groupe de  Lusztig-Spaltenstein de $X$, le centralisateur  $G_X$ de $X$ est inclus dans $P$ .
\end{lemme}
\end{paragr}

\begin{paragr}[Endomorphisme nilpotent simple.] --- On introduit la définition suivante.

\begin{definition}\label{def:simple}
  On dit que $X$ est \emph{simple} s'il vérifie l'une des deux conditions équivalente du lemme  \ref{lem:norm-jord} ci-dessous.
\end{definition}

\begin{lemme}\label{lem:norm-jord}
Soit $P\in \rc(X)$ et $M$ un facteur de Levi de $P$. Les assertions suivantes sont équivalentes.
 \begin{enumerate}
  \item $\Norm_G(M)=M$ ;
  \item L'invariant $\inv(X)$ de \eqref{eq:invX} est maximal c'est-à-dire égal à l'indice $r$ de nilpotence de $X$.
  \end{enumerate}
\end{lemme}

\begin{preuve}
D'après le lemme  \ref{lem:Richard2}, le groupe $M$ est le stabilisateur de $r$ sous-espaces de $E$ en somme directe et dont la somme vaut $E$. De plus, la suite ordonnée de ces dimensions est exactement la suite \eqref{eq:ordonnee}. On a donc $\Norm_G(M)=M$  si et seulement cette suite est strictement croissante.  Le lemme s'ensuit.
\end{preuve}
  
\end{paragr}

\begin{paragr}[Paramétrage de $\rc(X)$.] --- C'est l'objet de la proposition ci-dessous.

  \begin{proposition}\label{prop:P-R}  Soit $M$ un facteur de Levi d'un élément de $\rc(X)$. Alors

  \begin{enumerate}
  \item Pour tout $P\in \pc(M)$, il existe un unique élément $\tilde{P}\in \rc(X)$ tel que $P$ et $\tilde{P}$ soient conjugués sous $G$.
  \item L'application
    \begin{equation}
      \label{eq:lappli}
      \begin{array}{ccc}
  \pc(M) &\to& \rc(X)\\ P &\mapsto &\tilde{P}
\end{array}
\end{equation}

est surjective.
\item Les fibres de l'application \eqref{eq:lappli} sont naturellement des torseurs sous le groupe 
$$\Norm_G(M)/M.$$
\item L'application \eqref{eq:lappli} est bijective si et seulement si $X$ est simple au sens de la définition \ref{def:simple}. 
 
\end{enumerate}
\end{proposition}

 \begin{remarque}
    Comme $\pc(M)$ est fini, on obtient en particulier l'ensemble $\rc(X)$, et donc aussi $\lc\Sc(X)$, est fini.
  \end{remarque} 

  \begin{preuve} Soit $\oc$ la $G$-orbite de $X$.

Commençons par montrer l'unicité dans l'assertion 1 : autrement dit deux éléments de $\rc(X)$ qui sont conjugués sont égaux. Soit $P_1$ et $P_2$ dans $\rc(X)$. Pour alléger les notations, on pose $\ngo_i=\ngo_{P_i}$.  Par les propriétés de l'induction de Lusztig-Spaltenstein qui ont été  rappelées au §\ref{S:ind-LS}, on sait que 
$$\oc \cap \ngo_i$$
est exactement la $P_i$-orbite de $X$.  Soit $g\in G$ tel que $P_1=gP_2 g^{-1}$. Comme on a $X\in \ngo_2$, on a aussi $gXg^{-1}\in \ngo_1\cap \oc$.  Il existe donc $p\in P_1$ tel que 
$$pXp^{-1}=gXg^{-1}$$
soit encore $p^{-1}g\in G_X$. Or on a $G_X\subset P_1$ (cf. lemme \ref{rq:centralisateur}). On en déduit  $g\in P_1$ d'où 
$$P_2= g^{-1}P_1 g=P_1.$$

Montrons ensuite l'existence dans l'assertion 1 c'est-à-dire pour tout $P\in \pc(M)$, il existe un élément de $\rc(X)$ qui est conjugué à $P$. Soit $P\in \pc(M)$. D'après le lemme \ref{lem:Richard2}, on  a $I_P^G(0)=\oc$. En particulier, il existe donc $g\in  G$ tel que $gXg^{-1}\in \ngo_P$. Pour $Q=g^{-1}Pg$, on a donc $X\in \ngo_Q$ et $I_Q^G(0)=\oc$ (toujours par  le lemme \ref{lem:Richard2}). On a donc $Q\in \rc(X)$.

Continuons par la surjectivité de \eqref{eq:lappli}. Soit $P$ et $Q$ dans $\rc(X)$. D'après le lemme \ref{lem:Richard2}, les sous-groupes paraboliques $P$ et $Q$ sont associés. Si l'on suppose de plus que $P\in \pc(M)$, on en déduit l'existence de $g\in G$ tel que $Q\in \pc(gMg^{-1})$ c'est-à-dir $g^{-1}Qg\in \pc(M)$ ce qui montre la surjectivité.

Calculons une fibre de \eqref{eq:lappli}. Soit $P\in \pc(M)$.  La fibre de l'image de $P$ par  \eqref{eq:lappli} est formée des éléments  $\pc(M)$ qui sont $G$-conjugués à $P$. Soit $g\in G/P$ tel que $Q=gPg^{-1}\in \pc(M)$. Identifions $g$ à un représentant dans $G$. Les groupes $g^{-1}Mg$ et $M$ sont deux facteurs de Levi de $P$. Ils sont donc conjugués par un élément de $P$. Quitte à changer le représentant $g$, on peut supposer que  $g^{-1}Mg=M$ c'est-à-dire $g\in \Norm_G(M)$. Comme  $\Norm_P(M)=M$, une fibre de \eqref{eq:lappli} est bien un torseur sous $\Norm_G(M)/M$.

Enfin 4 est simplement la combinaison de 3 et 4. 
  \end{preuve} 
\end{paragr}

\begin{paragr}[Paramétrage de $\lc\Sc(X)$.] --- C'est la proposition suivante.

 \begin{proposition}\label{prop:levi-ls}
    Soit $M$ un facteur de Levi d'un élément de $\rc(X)$. Il existe une unique application
    \begin{equation}
      \label{eq:lappli2}
      \fc(M)\to \lc\Sc(X)
    \end{equation}
    telle que
\begin{enumerate}
\item elle coïncide sur $\pc(M)$ avec l'application \eqref{eq:lappli} ;
\item un élément de $\fc(M)$ est conjugué à son image ;
\item elle préserve les inclusions ;
\item elle est surjective ;
\item le groupe  $\Norm_G(M)$ agit transitivement sur les fibres de \eqref{eq:lappli2}. Le stabilisateur de $Q\in \fc(M)$ est $\Norm_Q(M)$
\item l'application \eqref{eq:lappli2} est bijective si et seulement $X$ est simple au sens de la définition \ref{def:simple}.
\end{enumerate}
  \end{proposition}

  \begin{preuve}
    Montrons d'abord l'unicité. Notons $Q\mapsto \tilde{Q}$ l'application \eqref{eq:lappli2}. Soit $Q\in \fc(M)$ et $P\in \pc(M)$ tel que $P\subset Q$. Par la propriété 1, $\tilde{P}$ est l'image de $P$ par l'application \eqref{eq:lappli}. Donc il existe un unique élément $g\in G/P$ tel que $gPg^{-1}=\tilde{P}$. Par conséquent, on $\tilde{P}\subset gQg^{-1}$. Par la propriété 3, on a $\tilde{P}\subset \tilde{Q}$ et par la propriété 2, le groupe  $\tilde{Q}$ est conjugué à $gQg^{-1}$. Comme   $\tilde{Q}$ et  $gQg^{-1}$ contiennent tous deux  $\tilde{P}$, on a nécessairement
    \begin{equation}
      \label{eq:cj}
      \tilde{Q}=gQg^{-1}.
    \end{equation}
 
Réciproquement, l'équation \eqref{eq:cj} permet de définir    l'application \eqref{eq:lappli2} : d'après le lemme \ref{lem:leg} ci-dessous, le membre de droite de \eqref{eq:cj} ne dépend pas du choix de  $P\in \pc(M)$ tel que $P\subset Q$. Les propriétés 1 à 3 sont alors évidentes.

Prouvons la surjectivité. Soit $Q\in \fc(X)$. Il existe $P\in \rc(X)$ tel que $P\subset Q$. D'après la surjectivité de  \eqref{eq:lappli} (cf. proposition \ref{prop:P-R}), il existe $g\in G$ tel que $gPg^{-1}\in \pc(M)$, donc $gQg^{-1}\in \fc(M)$ et l'image de $gQg^{-1}$ par  \eqref{eq:lappli2} est donc $Q$.

Prouvons l'assertion 5. Soit $Q_1$ et $Q_2$ deux éléments de $\fc(M)$ qui ont même image, notée $Q$, dans $\lc\Sc(X)$ par l'application \eqref{eq:lappli2}. Pour $i\in \{1,2\}$, soit $P_i\in \pc(M)$ et $g_i\in G$ tels que $P_i\subset Q_i$ et $g_i P_ig_i^{-1}\in \rc(X)$. Par construction même de l'application \eqref{eq:lappli2}, on a
\begin{equation}
  \label{eq:Q}
  Q=g_1Q_1 g_1^{-1}=g_2Q_2 g_2^{-1} \in \lc\Sc(X).
\end{equation}
On note $\ngo_i=\ngo_{P_i}$.  On a $X\in \ngo_1\cap \ngo_2$ et on écrit 
$$X=Y+U$$ 
avec $Y\in \mgo_Q\cap\ngo_1\cap\ngo_2$ et $U\in \ngo_Q$. Comme $X\in I_{g_iP_ig_i^{-1}}^G(0)$, on a aussi $Y\in I_{g_iP_ig_i^{-1}}^Q(0)$ (cf. lemme \ref{lem:proj-induite}). Donc 
$$ I_{g_1P_1g_1^{-1}}^Q(0)=I_{g_2P_2g_2^{-1}}^Q(0).$$
D'après le lemme \ref{lem:Richard2} (appliqué à $G=M_Q$ qui est un produit de groupes linéaires), les sous-groupes paraboliques    $g_iP_ig_i^{-1}\cap M_Q$   de $M_Q$ sont associés. Il existe donc $m\in {M_Q}$ tel qu'on ait 
$$ m  g_1M g_1^{-1}m^{-1}= g_2M g_2^{-1}$$
c'est-à-dire
$$x=g_2^{-1}mg_1\in \Norm_G(M).$$
On conclut en utilisant \eqref{eq:Q} d'où l'on tire
$$xQ_1x^{-1}= g_2^{-1}m  Q m^{-1}g_2=g_2^{-1}Q g_2=Q_2.$$
Réciproquement si $Q\in \fc(M)$ et si $x\in \Norm_G(M)$ alors pour tout $P\in \pc(M)$ tel que $P\subset Q$, on a $xPx^{-1}\in \pc(M)$ donc $xQx^{-1}$ et $Q$ ont même image par \eqref{eq:lappli2}.   Le stabilisateur de $Q$ sous l'action par conjugaison de $G$ est évidemment$Q$.

Prouvons enfin l'assertion 6. L'injectivité de \eqref{eq:lappli2} implique celle de \eqref{eq:lappli} et donc que $X$ vérifie la condition de l'assertion 6. Réciproquement, on sait (cf. proposition \ref{prop:P-R} assertion 4, que la condition sur $X$ implique qu'on a $\Norm_G(M)=M$. La bijectivité est alors une conséquence immédiate des assertions 4 et 5.
  \end{preuve}

\begin{lemme}\label{lem:leg}
Soit $M$ un facteur de Levi d'un élément de $\rc(X)$ et $Q\in \fc(M)$. 

Pour tout $i\in \{1,2\}$ soit  $P_i\in \pc(M)$ et $g_i\in G$ tels que 
$$P_i\subset Q \text{  et   }g_iP_ig_i^{-1}\in \rc(X).$$
 Alors $g_2$ appartient à l'ensemble
$$G_X g_1 M_Q P_2.$$
En particulier, on  a
$$g_2Qg^{-1}_2=g_1Qg_1^{-1}$$
\end{lemme}

\begin{preuve}  Soit $\oc$ l'orbite de $X$. Soit  $P_i\in \pc(M)$ et $g_i\in G$ tels que $P_i\subset Q$   et $g_iP_ig_i^{-1}\in \rc(X).$ Pour alléger un peu les notations, on pose $\ngo_{P_i}=\ngo_i$. Soit $Y_i=g_i^{-1}X_ig\in \ngo_i$. Alors $\oc\cap \ngo_i$ est un ouvert dense de $\ngo_i$ et c'est la $P_i$-orbite de $Y$. Écrivons 
$$Y_i=Z_i+U_i$$
avec $Z_i\in \mgo_Q\cap \ngo_i$ et $U_i\in \ngo_Q$. On a $Y_i\in \oc \cap \ngo_i$ et $\oc=I_{P_i}^G$. D'après le lemme \ref{lem:proj-induite}, l'induite  $I_{P_i\cap M_Q}^{M_Q}(0)$ est la $M_Q$-orbite de $Y_i$.  Les sous-groupe paraboliques $P_1\cap M_Q$ et $P_2\cap M_Q$ de $M_Q$ ont même facteur de Levi $M$. Il résulte du lemme \ref{lem:Richard2}, que les  orbites induites
\begin{equation}
  \label{eq:=orb}
  I_{P_1\cap M_Q}^{M_Q}(0)=  I_{P_2\cap M_Q}^{M_Q}(0)
\end{equation}
sont égales. Il existe donc $h\in M_Q$ tel que $hZ_1h^{-1}=Z_2$. Mais alors $hY_1h^{-1}$ et $Y_2$ appartient à $\ngo_2\cap \oc$. Il existe donc  $p_2\in P_2$ tel que 
$$p_2h Y_1 (p_2h)^{-1}= Y_2$$
d'où
$$(p_2hg_1^{-1})X (p_2hg^{-1}_1)^{-1}= g^{-1}_2Xg_2.$$
Ainsi $g_2 p_2 h g^{-1}_1 \in G_X$ d'où
$$g_2\in G_X g_1  h^{-1} p_2^{-1}\subset G_X g_1  M_Q P_2$$
ce qu'il fallait démontrer. On a donc
$$g_2 Qg_2^{-1}= x g_1Q g_1^{-1}x^{-1}$$
pour un certain $x\in G_X$. Or on a $g_1P_1g_1^{-1}\in \rc(X)$ d'où
$$G_X\subset  g_1P_1g_1^{-1}\subset g_1Q g_1^{-1}$$
d'où $x g_1Q g_1^{-1}x^{-1}= g_1Q g_1^{-1}$.
\end{preuve} 
  \end{paragr}

\section{Forme standard d'un endomorphisme nilpotent}\label{sec:standard}

\begin{paragr}[Nilpotent standard.] --- \label{S:surZ}
  Soit $n\geq 1$ et $r\geq 1$ des entiers. Soit $E$ un $\ZZ$-module libre de rang $n$  muni d'une décomposition en somme directe 
$$E= \bigoplus_{1\leq i \leq j\leq r} V_j^i$$
telle que le rang $d_j\geq 0$ de $V_j^i$ ne dépende que de $j$. Pour tous $1\leq i\leq j\leq r$ et $1\leq k\leq d_j$ soit $e^i_{k,j}\in V_j^i$ tel que la famille $(e_{k,j}^i)_{1\leq k\leq d_j}$ soit une base du $\ZZ$-module $V_j^i$. 
On obtient ainsi une base $e$ de $E$ qu'on ordonne de la façon suivante : $e^i_{k,j}< e^{i'}_{k',j'}$ si  l'une des conditions suivantes est satisfaite

\begin{itemize}
\item $i<i'$ ;
\item $i=i'$ et $j>j'$ ;
\item $i=i'$, $j=j'$ et $k<k'$.
\end{itemize}

De la sorte, on identifie $E$ à $\ZZ^n$ et le schéma en groupe $\gc=\Aut_\ZZ(E)$ à $GL(n)_\ZZ$. Soit $\bc$ et $\tc$ les sous-schémas en sous-groupes de Borel standard et en sous-tores maximaux standard.

Soit $X\in \End_\ZZ(E)$ défini par 
$$Xe_{k,j}^i=\left\lbrace
  \begin{array}{l}
    e_{k,j}^{i-1} \text{ si } i> 1 ;\\
0  \text{   si } i= 1.
  \end{array}\right.
$$
pour $1\leq i\leq j\leq r$ et $1\leq k\leq d_j$.
On observe que le $k$-ème noyau itéré s'écrit
$$\Ker(X^k)=\bigoplus_{1\leq j\leq r} \bigoplus_{1\leq i \leq \min(j,k)} V_j^i.$$

Soit $F$ un corps et $G=\mathcal{G}\times_\ZZ F$. Soit $B=\bc\times_\ZZ F$ et $T_0=\tc\times_\ZZ F$. L'endomorphisme $X$ s'identifie à un élément nilpotent de $\ggo(F)$ et tout élément nilpotent obtenu de cette façon est dit  \og standard\fg. Tout élément nilpotent de $\ggo(F)$ est conjugué sous $G(F)$ à un nilpotent standard.

 Le sous-groupe parabolique $P_0$ de $G$ qui stabilise le drapeau des noyaux itérés de $X$ est standard au sens où il contient $B$. Soit
$$M=M_{P_0}$$
l'unique facteur de Levi de $P_0$ qui contient $M$. Un sous-groupe parabolique, resp. de Levi, est dit semi-standard s'il appartient à $\fc(T_0)$, resp. $\lc(T_0)$, cf.  §\ref{S:fM}.

\begin{lemme}\label{lem:tous-std}
  On a l'inclusion $\lc\Sc(X)\subset \fc(T_0)$.
\end{lemme}

\begin{preuve}
Par la proposition \ref{prop:drap-R} ci-dessous, on a $\rc(X)\subset  \fc(T_0)$. Comme tout élément de $\lc\Sc(X)$ contient un élément de $\rc(X)$, le lemme en résulte.
\end{preuve}
\end{paragr}

\begin{paragr}[Groupe de Weyl et produit scalaire sur $a_{T_0}$.]---\label{S:Weyl} Soit $\Norm_G(T_0)/T_0$ le groupe de Weyl de $(G,T_0)$. Dans toute la suite, \emph{on l'identifie} au sous-groupe $W\subset \Norm_G(T_0)$ des matrices de permutation de la base $e$. Pour tout $L\in \lc(T_0)$, soit $W^L=W\cap L$. Le groupe $W$ agit naturellement sur $a_{T_0}$ et son dual. La base du $\ZZ$-module $E$ (cf. \ref{S:surZ}) permet d'identifier l'espace $a_{T_0}$ à l'espace euclidien $\RR^n$. Cette identification munit $a_{T_0}$ d'un produit scalaire invariant sous $W$. On note $\|\cdot\|$ la norme euclidienne sur $a_{T_0}$ qui s'en déduit. Dans la suite, pour tout $L\in \lc(T_0)$, on identifie $a_L$ à l'orthogonal de $a_{T_0}^L$ dans $a_{T_0}$ ce qui munit $a_L$ du produit scalaire induit. Des identifications et des notations analogues valent pour l'espace dual $a_{T_0}^*$ et ses sous-espaces. 
  \end{paragr}

\begin{paragr}[Élément $w_P$.] ---

\begin{lemme}\label{lem:r-std}
  Pour tout $P\in \fc(M)$, il existe un élément $w_P\in W$, unique à translation à gauche près par $W^{M_P}$  tel que 
$$w_P^{-1}Pw_P\in \lc\Sc(X)$$
soit l'image de $P$ par l'application \eqref{eq:lappli2}. 
Si $P\in \pc(M)$, la condition $w_P^{-1}Pw_P\in \rc(X)$ suffit pour que $w_P^{-1}Pw_P$ soit l'image de $P$ par  l'application \eqref{eq:lappli2}.
\end{lemme}
  
\begin{preuve}
Montrons d'abord le lemme lorsque $P\in \pc(M)$. D'après la proposition \ref{prop:levi-ls}, il existe un unique élément $\tilde{P}\in \rc(X)$ tel que $\tilde{P}$ soit conjugué à $P$. Donc il existe un unique élément $g\in P\back G$ tel que $g^{-1}Pg\in \rc(X)$. D'après le lemme \ref{lem:tous-std},  le sous-groupe parabolique $g^{-1}Pg$ est standard. En identifiant $g$ à un représentant dans $G$, on a $gT_0g^{-1}\subset P$. Or $T_0\subset P$. Ces deux sous-tores déployés sont donc $P$-conjugués : il existe $p\in P$ tel que $pg\in \Norm_G(T_0)$. Quitte à translater $p$ à gauche par un élément de $T_0$, on peut supposer que $pg\in W$. D'où l'existence de $w_P$. L'unicité résulte de l'égalité $\Norm_W(P)=W^M$.

Soit $Q\in \fc(M)$. Il existe $P\in \pc(M)$ tel que $P\subset Q$. Partant de $w_P\in W$ tel que $w_P^{-1}Pw_P\in \rc(X)$, on obtient   $w_P^{-1}Pw_P\in \rc(X)$ (cf. éq. \eqref{eq:cj} dans la démonstration de la proposition \ref{prop:levi-ls}). Cela donne l'existence. L'unicité résulte de $\Norm_W(Q)=W^{M_Q}$.

\end{preuve}
  
\end{paragr}

\begin{paragr}[Une description de $\rc(X)$ en terme d'algèbre linéaire.] ---  Soit   
$$J=\{j\in \NN^* \mid d_j\not=0)\}$$
et 
$$r=\max(J).$$ 
L'entier $r$ est encore l'indice de nilpotence de $X$ et $d_j$ est la multiplicité du bloc de taille $j$ dans la décomposition de Jordan de $X$. Le cardinal $|J|$ est l'entier noté $\inv(X)$, cf. \eqref{eq:invX}.

Soit  
$$\ec=\ec(X)$$  
l'ensemble (fini) des applications
$$\eps :
\begin{array}{ccc}
  \{k \in \NN \mid 1 \leq k \leq r\} \times  J & \to & \{0,1\}\\
(k,j) & \mapsto &\eps_{k,j}
\end{array}$$
qui vérifient pour tout $j\in J$
\begin{enumerate}
\item $\sum_{k=1}^r  \eps_{k,j}=j$ ;
\item pour tout $1\leq k\leq r$, l'application 
$$\begin{array}{{ccc}}
  J& \to & \{0,1\} \\
j&\mapsto& \eps_{k,j} 
\end{array}$$
est croissante.
\end{enumerate}

  La proposition suivante détermine explicitement l'ensemble $\rc(X)$.

\begin{proposition}\label{prop:drap-R}
Pour tout $\eps \in \ec$, soit  $E_\bullet(\eps)$ le drapeau de $E$ défini par $E_0(\eps)=(0)$ et pour $1\leq k \leq r$
$$E_k(\eps) =\bigoplus_{j\in J}\bigoplus_{i=1}^{\sum_{l=1}^k \eps_{l,j}} V_j^i.$$
L'application 
$$\eps \mapsto P_\eps=\stab_G(E_\bullet(\eps)\otimes_\ZZ F)$$
est une bijection de $\ec$ sur $\rc(X)$.
\end{proposition}

\begin{preuve}
Elle est donnée au paragraphe suivant \ref{S:preuve-drapR}.
\end{preuve}

\begin{remarque} \label{rq:particulier}Soit $\xi\in \ec$ défini par 
$$\xi_{k,j}=\left\lbrace
  \begin{array}{l}
    1 \text{ si } k\leq j\ ;\\
0 \text{ sinon.}
  \end{array}\right.$$
  Le drapeau correspondant $E_\bullet(\xi)$ est le drapeau des noyaux itérés de $X$. Le drapeau des images itérées de $X$ est le drapeau   $E_\bullet(\tilde{\xi})$  où $\tilde{\xi}\in \ec$ est défini par 
$$\tilde{\xi}_{k,j}=\left\lbrace
  \begin{array}{l}
    1 \text{ si } r-k\leq j\ ;\\
0 \text{ sinon.}
  \end{array}\right.$$
\end{remarque}
\end{paragr}

\begin{paragr}[Démonstration de la proposition \ref{prop:drap-R}.] --- \label{S:preuve-drapR}Le lemme \ref{lem:inj}ci-dessous montre  que l'application $\eps\mapsto P_\eps$ est injective. La surjectivité résulte du  lemme \ref{lem:surj}. Commençons par un lemme auxiliaire.

  \begin{lemme}\label{lem:transitif}
    Pour tout $\sigma$ dans le groupe symétrique $\SG_r$ et tout $\eps\in \ec$, l'application 
$$\sigma(\eps): (k,j)\mapsto \eps_{\sigma^{-1}(k),j}$$
 appartient à $\ec$. Cela définit une action \emph{transitive} de $\SG_r$ sur $\ec$.
  \end{lemme}

\begin{preuve}  
 Soit $\eps\in \ec$. On va montrer qu'il existe $\sigma\in \SG_r$ tel que $\xi=\sigma(\eps)$ où $\xi$ est défini dans la remarque \ref{rq:particulier}. Soit $s=\min(J)$. Par la condition 1, on a $\sum_{k=1}^r \eps_{k,s}=s$. Il y a donc exactement $s$ éléments de $\{1,\ldots,r\}$ tels que $\eps_{k,s}=1$. Quitte à remplacer $\eps$ par $\sigma(\eps)$ pour $\sigma\in \SG_r$, on peut et on va supposer qu'on a $\eps_{k,s}=1$ si et seulement si $1\leq k\leq s$. Soit $j_0\in J$ le plus grand élément de $J$ tel que pour tout $j\in J$ avec $j\leq j_0$ et  tout $1\leq k\leq r$ on ait 
\begin{equation}
  \label{eq:kj}
\eps_{k,j}=\left\lbrace
  \begin{array}{l}
    1 \text{ si } k\leq j\\
0 \text{ sinon }
  \end{array}\right.
\end{equation}
Si $j_0=r$, on a $\eps=\xi$. Sinon soit $j_1$ le plus petit élément de $J$ strictement plus grand que $j_0$.
On a $\eps_{k,j_1}=1$ pour tout  $k\leq j_0$ (par croissance de $j\mapsto \eps_{k,j}$). On a  donc
$$j_0 + \sum_{k=j_0+1}^r \eps_{k,j_1}= \sum_{k=1}^r \eps_{k,j_1}=  j_1.$$
Il existe donc exactement $j_1-j_0>0$ éléments de $\{j_0+1,j_0+2,\ldots, r\}$ tels que $\eps_{k,j_1}=1$. Il existe donc $\sigma\in \SG_r$ qui permute ces éléments tels que $\sigma(\eps)$ vérifie \eqref{eq:kj}  pour $j\leq j_1$ et $1\leq k\leq r$. De proche en proche, on obtient bien $\xi$. Ainsi l'action de $\SG_r$ sur $\ec$ est transitive.
\end{preuve}

\begin{lemme}\label{lem:inj}
  L'application $\eps \in \ec \mapsto P_\eps$ est une injection de $\ec$ dans $\rc(X)$.
\end{lemme}

\begin{preuve} L'application est clairement injective. Il s'agit de montrer que $P_\eps\in \rc(X)$ pour tout $\eps\in \ec$. Pour tout  $\eps\in \ec$, on a par construction $X\in \ngo_{P_\eps}$. Pour tout $\eps\in \ec$, on a 
$$\mathrm{rang}(E_k(\eps)/E_{k-1}(\eps))=\sum_{j\in J}  \eps_{k,j} d_j.$$
Il en résulte que l'action de $\SG_r$ sur $\ec$  préserve la suite ordonnée des rangs $\mathrm{rang}(E_k(\eps)/E_{k-1}(\eps))$. On déduit alors du lemme \ref{lem:Richard2} que $P_\eps\in \rc(X)$ si et seulement si   $P_{\sigma(\eps)}\in \rc(X)$ pour un élément $\sigma\in \SG_r$. Comme l'action de $\SG_r$ sur $\ec$ est transitive, il suffit de trouver un élément $\eps$ tel que $P_\eps\in \rc(X)$. Il suffit de prendre $\eps=\xi$ défini à la remarque \ref{rq:particulier} qui donne le stabilisateur des noyaux itérés.
\end{preuve}

\begin{lemme}\label{lem:surj}
  Tout drapeau de Richardson de $X$ est de la forme $E_\bullet(\eps)$  pour $\eps\in \ec$.
\end{lemme}

\begin{preuve} On raisonne par récurrence sur l'entier $r$. Lorsque $r=1$, l'endomorphisme nilpotent $X$ est nul, $\rc(X)=\{G\}$ et le résultat est évident  ce qui amorce la récurrence. On suppose désormais $r\geq 2$.

Soit $E_\bullet$ un drapeau de sous-$F$-espaces de $E\otimes_\ZZ F$ dont le stabilisateur est un sous-groupe parabolique de Richardson $P$ de $X$. On sait, par le lemme \ref{lem:Richard2}, décrire la suite ordonnée des dimensions du gradué de $E_\bullet$. En particulier, il existe $1\leq j_0\leq r$ tel que 
$$\dim(E_r/E_{r-1})=d_r+d_{r-1}+\ldots+d_{j_0}.$$

On va appliquer  l'hypothèse de récurrence au sous-$F$-espace vectoriel $E_{r-1}$ muni de l'endomorphisme nilpotent $Y$ donné par la restriction de $X$ à $F$. Soit $Q\subset G$ le sous-groupe parabolique maximal qui stabilise $E_{r-1}$. On a $P\subset Q$ ; il résulte du lemme \ref{lem:proj-induite} qu'on a $Y\in I_P^Q(0)$. Il résulte alors d'une nouvelle application du lemme \ref{lem:Richard2} que l'indice de nilpotence de $Y$ est $r'=r-1$ et que les multiplicités  $(\delta_j)_{1\leq j \leq r'}$ de ses blocs de Jordan vérifient
$$
\delta_{r-1}+\delta_{r-1}+\ldots+\delta_j=\left\lbrace
  \begin{array}{c}
    d_r+\ldots +d_{j+1} \text{ si } r> j\geq j_0 \\
 d_r+\ldots +d_{j} \text{ si } j<j_0.
  \end{array}
\right.
$$

On en déduit que la multiplicité du bloc de Jordan de taille $j$ de $Y$ est donnée par 
$$\delta_j=\left\lbrace
  \begin{array}{c}
    d_{j+1} \text{ si } j_0\leq j <r\\ 
d_j+d_{j+1} \text{ si } j=j_0-1 \\ 
d_j \text{ si } j<j_0-1 
  \end{array}\right.$$

Posons pour $j\geq 1$ et $1\leq i\leq j$ 

$$W_j^i=\left\lbrace
  \begin{array}{c}
    V_{j+1}^i  \text{ si } j\geq j_0\\ 
V_j^i\oplus V_{j+1}^i \text{ si } j=j_0-1 \\ 
V_j^i \text{ si } j<j_0-1 
  \end{array}\right.$$

\begin{lemme} \label{lem:egalite}On a 
  \begin{equation}
  \label{eq:egalite}
  (\bigoplus_{j\geq 1,1\leq i\leq j} W_j^i)\otimes_\ZZ F =  E_{r-1}.
\end{equation}
\end{lemme}

\begin{preuve}
Il suffit de prouver l'inclusion pour $j\geq 1$ et $1\leq i\leq j$
\begin{equation}
  \label{eq:inclusion}
   W_j^i\subset F.
\end{equation}
car les deux membres de l'égalité à prouver ont même dimension. On a $\Ima(X)\subset E_{r-1}$. Donc  pour $1\leq i<j$ on a $V_j^{i}=X(V_j^{i+1})\subset E_{r-1}$. En particulier, pour $j\geq j_0$ et $j\geq i$, on a $W_j^i=V_{j+1}^i\subset E_{r-1}$. On a terminé si $j_0=1$. Si $j_0>1$, on utilise l'inclusion 
\begin{equation}
  \label{eq:uneinclusion}
  \Ker(X^{j_0-1})\subset E_{r-1}
\end{equation}
autrement dit l'égalité 
$$\Ker(Y^{j_0-1})=\Ker(X^{j_0-1}).$$
Cette dernière résulte de l'inclusion $\subset$ évidente et de l'égalité des  dimensions (ces dimensions se calculent aisément en termes des multiplicités $d_j$ et $\delta_j$) 
$$\dim(\Ker(X^{j_0-1}))=\sum_{k=1}^{j_0-1}  \sum_{j\geq k} d_j =  \sum_{k=1}^{j_0-1}  \sum_{j\geq k} \delta_j=  \dim(\Ker(Y^{j_0-1})).$$
De \eqref{eq:uneinclusion}, on déduit $\Ker(X^{j_0-1})\subset F$ et donc 
$$V_j^{i}\subset E_{r-1}$$
pour $1\leq i\leq j_0-1$. Cela donne les inclusions \eqref{eq:inclusion} pour $j\leq j_0-1$. 
\end{preuve}

L'endomorphisme $Y$ de $E_{r-1}$ est alors standard pour la base extraite de $e$. L'hypothèse de récurrence entraîne l'existence de  $\xi\in \ec(Y)$ tel que pour tout $1\leq k\leq r-1$, on ait
$$E_k= (\bigoplus_{j\in J'}\bigoplus_{i=1}^{\sum_{l=1}^k \xi_{l,j}} W_j^i)\otimes_\ZZ F.$$
Ici on a $J'=\{j\geq 1 \mid \delta_j\not=0)\}$ avec $\max(J')=r-1$. L'application $\xi$ va de $\{k\in \NN \mid 1\leq k\leq r-1\}\times J'$ dans $\{0,1\}$. On a  $j\in J'$ si et seulement si l'une des trois conditions suivantes est réalisée :
\begin{itemize}
\item $j+1\in J$ et $j\geq j_0$ ;
\item $j\in J$ ou $j+1\in J$ et $j=j_0-1$ ;
\item $j\in J$ si $j<j_0-1$.
\end{itemize}
On pose alors pour tout $j\in J$

$$\eps_{k,j}=\left\lbrace
  \begin{array}{c}
    \xi_{k,j-1} \text{ si } j\geq j_0\\ 
\xi_{k,j} \text{ si } j< j_0\\
  \end{array}\right.$$
pour $k<r$. On a donc
$$\sum_{k=1}^{r-1}\eps_{k,j}=\left\lbrace
  \begin{array}{c}
    j-1 \text{ si } j\geq j_0\\ 
 j \text{ si } j< j_0\\
  \end{array}\right.$$
On pose alors 
$$\eps_{k,r}=\left\lbrace
  \begin{array}{c}
    1 \text{ si } j\geq j_0\\ 
 0 \text{ si } j< j_0\\
  \end{array}\right.$$
de sorte que $k\mapsto \eps_{k,r}$ est croissante et $\sum_{k=1}^{r}\eps_{k,j}=j$ pour tout $j\in J$.

On vérifie immédiatement qu'on a $\eps\in \ec$ et $E_\bullet=E_\bullet(\eps)\otimes_\ZZ F$.
\end{preuve}
\end{paragr}

\section{Un calcul d'orbite}\label{sec:orb}

\begin{paragr}[Décomposition de Levi du centralisateur.] \label{S:Levi-cent} --- On reprend les notations du §\ref{S:surZ}. Le $\ZZ$-module des endomorphismes de $E$ se décompose
  \begin{equation}
    \label{eq:decomp}
    \End_\ZZ(E)=\bigoplus \Hom_\ZZ(V_j^i,V_{j'}^{i'})
  \end{equation}
  où la somme porte sur les $1\leq i \leq j\leq r$ et $1\leq i' \leq j'\leq r$. On a donc des injections 
$$\Hom_\ZZ(V_j^i,V_{j'}^{i'})\hookrightarrow \End_\ZZ(E).$$
Soit $\ggo= \End_\ZZ(E)$ l'algèbre de Lie de $\gc=\Aut_\ZZ(E)$.   Soit $\ggo_X\subset \ggo$ la sous-algèbre des endomorphismes commutant à $X$. Le morphisme  de restriction
  \begin{equation}
    \label{eq:isom-endo}
    \ggo_X \to \prod_{j=1}^r \Hom_\ZZ(V_j^j, \Ker(X^j))
  \end{equation}
  est un isomorphisme. Le noyau de la flèche $\ggo_X \to \prod_{j=1}^r \Hom_\ZZ(V_j^j, V_j^j)$ obtenue par composition avec les projections  $\Hom_\ZZ(V_j^j, \Ker(X^j)) \to \Hom_\ZZ(V_j^j, V_j^j)$ est une sous-algèbre notée $\ngo_X$. L'inclusion 
$$\prod_{j=1}^r \Hom_\ZZ(V_j^j, V_j^j)\subset \prod_{j=1}^r \Hom_\ZZ(V_j^j, \Ker(X^j))$$
fournit via l'isomorphisme \eqref{eq:isom-endo} un supplémentaire $\mgo_X$ de $\ngo_X$ dans $\ggo_X$. C'est en fait une sous-algèbre de $\ggo_X$. La décomposition $\ggo_X=\mgo_X\oplus\ngo_X$ est une \og décomposition de Levi\fg{} de l'algèbre de Lie $\ggo_X$.

Le centralisateur $\gc_X$ de $X$ dans $\gc$ admet aussi une décomposition de Levi $M_XN_X$ où $M_X=\mgo_X \cap G$ est un schéma en groupes réductifs et $N_X=(Id_E+\ngo_X\cap G)$ un schéma en groupes unipotents. Il résulte de qui précède que $M_X$ est isomorphe au produit 
 $$\prod_{1\leq j \leq r}   \Aut_\ZZ(V_j^j).$$
  
\end{paragr}

\begin{paragr} Affublons chaque $\ZZ$-module d'un poids
  \begin{eqnarray*}
  p(V_j^i)&=& r+(r-1)+\ldots+ (r-(i-2))+ r-j+1\\
&=& (i-1)(2r-i+2)/2+ r-j+1.
\end{eqnarray*}

\begin{lemme}\label{lem:obs}
  \begin{enumerate}
  \item   On a $p(V_j^i)>p(V_{j'}^{i'})$ si et seulement si l'une des deux conditions est satisfaite
    \begin{enumerate}
    \item $i>i'$ ;
    \item $i=i'$ et $j<j'$.
    \end{enumerate}
\item Soit $1<i \leq j$ et $1< i'\leq j'$. Supposons qu'on a $p(V_j^{i})\geq p(V_{j'}^{i'})$. Alors on a 
$$p(V_j^{i-1})-p(V_{j'}^{i'-1})\geq p(V_j^{i})-p(V_{j'}^{i'}) .$$
\end{enumerate}
\end{lemme}

\begin{preuve} Prouvons l'assertion 1.
  Supposons $i>i'$. Alors on a
  \begin{eqnarray*}
    p(V_j^i)-p(V_{j'}^{i'})&\geq & (r-i'+1) + (j'-j)\\
&=& (r-j)+(j'-i')+1\\
&\geq & 1
\end{eqnarray*}
Si $i=i'$ alors $p(V_j^i)-p(V_{j'}^{i'})= j'-j$ d'où l'assertion 1. Prouvons l'assertion 2. On a la formule
$p(V_j^{i-1})-p(V_j^{i})= i-r+2.$ d'où $p(V_j^{i-1})-p(V_{j'}^{i'-1}) -(p(V_j^{i})-p(V_{j'}^{i'}))= i-i'.$
Si $p(V_j^{i})\geq p(V_{j'}^{i'})$, on a $i\geq i'$ d'où l'assertion 2.
\end{preuve}

\end{paragr}

\begin{paragr}[Filtration $\ngo^{\geq t}$.]  \label{S:filtr}---  Pour tout $t\in \ZZ$, soit 
$$\ngo^{\geq t}=\bigoplus \Hom_\ZZ(V_j^i,V_{j'}^{i'})$$
où la somme porte sur les $1\leq i \leq j\leq r$ et $1\leq i' \leq j'\leq r$ tels que $p(V_j^i)-p(V_{j'}^{i'})\geq t.$
On a donc $\ngo^{\geq t}\subset \ngo^{\geq t'}$ pour $t\geq t'$. Pour tous $t$ et $t'$ dans $\ZZ$, on a aussi $\ngo^{\geq t}\cdot \ngo^{\geq t'}\subset  \ngo^{\geq t+t'}$ et en particulier $[ \ngo^{\geq t},\ngo^{\geq t'} ]\subset  \ngo^{\geq t+t'}$ pour le crochet de Lie.
\end{paragr}

\begin{paragr}[Sous-groupe parabolique $R$.] --- \label{S:R}Soit $\rgo=\ngo^{\geq 0}$. C'est encore la sous-algèbre de $\End_\ZZ(E)$ qui stabilise le raffinement du drapeau des noyaux itérés qu'on obtient en  intercalant pour $0\leq i \leq r-1$ le drapeau 
  $$\Ker(X^i)\subset \Ker(X^i)+ \Ker(X^{i+1})\cap \Ima(X^{r-i})\subset \ldots \subset  \Ker(X^i)+ \Ker(X^{i+1})\cap \Ima(X) \subset \Ker(X^{i+1}).$$
Bien sûr, les inclusions ci-dessus ne sont pas nécessairement strictes. En tout cas, on a $\ggo_X\subset \ngo^{\geq 0}$. On a une décomposition de Levi $\rgo=\mgo \oplus\ngo$ où $\ngo=\ngo^{\geq 1}$ et 
$$\mgo= \bigoplus \Hom_\ZZ(V_j^i,V_{j}^{i})$$
où la somme porte sur les $1\leq i \leq j\leq r$. On observera qu'on a $\mgo_X\subset \mgo$ et $\ngo_X\subset \ngo$.

 Soit $R$ le sous-schéma en groupes de $\gc$ défini par $R=\rgo\cap \gc$. Il admet une décomposition de Levi $R=MN$ avec $M= \mgo\cap G$ et $N=\Id_E+\ngo$.  On notera qu'on a 
$$M_X\subset M \text{  et  }N_X\subset N.$$
Le schéma en groupes réductifs $M$ s'identifie au produit $\prod_{1\leq i \leq j \leq r} \Aut_\ZZ(V_j^i)$ de sorte que l'inclusion $M_X\subset M$ s'identifie à l'application
$$\prod_{1\leq j \leq r} \Aut_\ZZ(V_j^j)\to \prod_{1\leq i \leq j \leq r} \Aut_\ZZ(V_j^i)$$
qui $j$ par $j$ est donnée par le plongement diagonal $ \Aut_\ZZ(V_j^j)\to \prod_{1\leq i \leq j } \Aut_\ZZ(V_j^i)$ où l'on identifie $\Aut_\ZZ(V_j^j)\simeq \Aut_\ZZ(V_j^i)$ via l'isomorphisme $V_j^j\simeq V_j^i$ induit par $X^{j-i}$.
\end{paragr}

\begin{paragr}[Description de la $N$-orbite de $X$.] --- \label{S:o} Soit $t\in \ZZ$ et $\of^{\geq t}= \bigoplus \Hom_\ZZ(V_j^i,V_{j'}^{i'})$ où la somme porte sur les $1< i \leq j\leq r$ et $1\leq i' \leq j'\leq r$ tels que $p(V_j^{i-1})-p(V_{j'}^{i'})\geq t.$ Comme $p(V_j^{i})> p(V_j^{i-1})$, on a $\of^{\geq t}\subset \ngo^{\geq t+1}$. Soit $\of=\of^{\geq 1}$. Soit $N^{\geq t}=Id_E+\ngo^{\geq t}$ et $N_X^{\geq t}=N_X\cap N^{\geq t}$.

  \begin{proposition}\label{prop:iso}
L'application 
\begin{equation}
  \label{eq:n}
  \begin{array}{lll}
    N & \to & \ngo \\
n&\mapsto& n^{-1}Xn -X 
  \end{array}
\end{equation}
induit un isomorphisme de $\ZZ$-schéma
\begin{equation}
  \label{eq:n-iso}
  N_X\back N \to \of.
\end{equation}
  \end{proposition}

  \begin{preuve}
On va démontrer par récurrence sur $t\geq 1$ que l'application \eqref{eq:n} induit un isomorphisme $ N_X^{\geq t}\back N^{\geq t} \to \of^{\geq t}.$  Elle est clairement injective. Le cas $t$ assez grand est trivial et le cas $t=1$ donne la proposition. Montrons que le cas $t-1$ implique le cas $t$ : soit $Y\in \of^{\geq t}$. D'après le lemme \ref{lem:iso} ci-dessous, il existe $n\in N^{\geq t}$ tel que $Y-(n^{-1}Xn-X)\in \of^{\geq t+1}.$ On a donc $n(Y-(n^{-1}Xn-X)n^{-1}\in \of^{\geq t+1}.$
Par hypothèse de récurrence, il existe $n_1\in N^{\geq t}$ tel que 
$$n(Y-(n^{-1}Xn-X)n^{-1}=n_1^{-1}Xn_1-X$$
d'où $Y= (n_1n)^{-1}X n_1n -X$ comme voulu.
  \end{preuve}

\begin{lemme}\label{lem:iso}
  Soit $t\geq 1$. On a un isomorphisme de $\ZZ$-schéma en groupes
  \begin{equation}
    \label{eq:etude}
(N_X^{\geq t} N^{\geq t+1})\back N^{\geq t} \to \of^{\geq t}/\of^{\geq t+1}
\end{equation}
induit par $n\mapsto n^{-1}Xn -X$.
\end{lemme}

\begin{preuve}
Observons tout d'abord que l'application $\ngo^{\geq t}\to N^{\geq t}$ donnée $U\mapsto \Id_E +U$ induit un isomorphisme de groupes
\begin{equation}
  \label{eq:isom-gp}
   (\ngo_X^{\geq t} +\ngo^{\geq t+1})\back \ngo^{\geq t} \to    (N_X^{\geq t} N^{\geq t+1})\back N^{\geq t} .
 \end{equation}
 Soit $U\in \ngo^{\geq t} $ et $n=\Id_E +U$. On a
\begin{equation}
  \label{eq:nX}
  n^{-1}Xn-X= n^{-1} [X,U].
\end{equation}
Prouvons que $[X,U]\subset \of^{\geq t}$. Il suffit de considérer un élément $U\in  \Hom(V_j^i,V_{j'}^{i'})$ avec $p(V_j^i)-p(V_{j'}^{i'})\geq t\geq 1$. Si $i=j$ on a $UX=0$ sinon on a $UX \in  \Hom(V_j^{i+1},V_{j'}^{i'})\subset \of^{\geq t}.$
Si $i'=1$, on a $XU=0$. Sinon $i'>1$ et $XU\in  \Hom(V_j^{i},V_{j'}^{i'-1}).$
Comme $p(V_j^i)-p(V_{j'}^{i'})\geq 1$, on a $i\geq i'>1$. Par le lemme \ref{lem:obs} assertion 2, on a 
$$p(V_j^{i-1})-p(V_{j'}^{i'-1})\geq p(V_j^{i})-p(V_{j'}^{i'})\geq t,$$
d'où $XU\in  \of^{\geq t}$. 
On a donc $n^{-1} [X,U]\in \of^{\geq t}$ et par \eqref{eq:nX}, on a $ n^{-1}Xn-X\in \of^{\geq t}$. De la sorte, on obtient une application
$$\ngo^{\geq t} \to \of^{\geq t}/\of^{\geq t+1}$$
qui, à $U\in \ngo^{\geq t} $, associe la classe de $(\Id_E+U)^{-1}X (\Id_E+U)-X$. Par \eqref{eq:nX} et le fait que  $[X,U]\in \of^{\geq t}$, cette application coïncide avec celle induite par $U\mapsto [X,U]$. Elle se factorise en une application
\begin{equation}
  \label{eq:etude2}
  (\ngo_X^{\geq t} +\ngo^{\geq t+1})\back \ngo^{\geq t} \to  \of^{\geq t}/\of^{\geq t+1}.
\end{equation}
Via l'isomorphisme \eqref{eq:isom-gp}, elle coïncide avec l'application \eqref{eq:etude}. Prouvons que \eqref{eq:etude2} est un isomorphisme. On peut raisonner par récurrence sur $t$ (pour $t$ assez grand les deux membres sont nuls). En utilisant l'hypothèse de récurrence pour $t-1$, on voit que \eqref{eq:etude2} est injective. 
Soit $Y\in \of^{\geq t}/\of^{\geq t+1}$ ; on l'identifie à un élément de $\bigoplus \Hom_\ZZ(V_j^i,V_{j'}^{i'})$ où la somme porte sur les $1< i \leq j\leq r$ et $1\leq i' \leq j'\leq r$ tels que $p(V_j^{i-1})-p(V_{j'}^{i'})=t.$
Soit $p(Y)$ le plus petit entier parmi les entiers $p(V_{j'}^{i'})$ tels que la projection $Y'$ de $Y$ sur $\Hom_\ZZ(V_j^i,V_{j'}^{i'})$ soit non nulle. S'il n'existe pas de tels entiers, on a $Y=0$. Sinon $p(Y)\geq 1$. Soit $Y'$ la projection de $Y$ sur $\Hom_\ZZ(V_j^i,V_{j'}^{i'})$ avec  $p(V_{j'}^{i'})=p(Y')$. Il existe $U\in \Hom_\ZZ(V_j^{i-1},V_{j'}^{i'})\subset \ngo^{\geq t}$ tel que $Y'=UX$. On a donc $Y-[U,X]\in \of^{\geq t}$. On a aussi $XU\in \Hom_\ZZ(V_j^{i-1},V_{j'}^{i'+1})$. On a $p(V_j^{i-1})<p(V_j^i)$ donc $p(Y-[U,X])<p(Y)$. De proche en proche, on construit $U\in  \ngo^{\geq t}$ tel que $Y\in [X,U]+\of^{\geq t+1}$.
\end{preuve}
\end{paragr}

\section{Une famille orthogonale associée à un endomorphisme nilpotent}\label{sec:fam-ortho}

\begin{paragr}[Fonctions $H_P$.] --- \label{S:HP} Dans ce paragraphe, on reprend la situation du §\ref{S:surZ}. On suppose que  $F$ est un corps local de caractéristique $0$. Soit $|\cdot|$ la valeur absolue normalisée de $F$. Soit $K\subset G(F)$ le sous-groupe compact maximal défini ainsi :
  \begin{itemize}
  \item si $F=\CC$, c'est le groupe unitaire associée à la forme hermitienne pour laquelle la base $e$ est orthonormale ;
  \item si $F=\RR$, c'est  le groupe othogonal associé à la forme quadratique pour laquelle la base $e$ est orthonormale ;
  \item si $F$ est non-archimédien, on pose $K=\gc(\oc)$ où $\oc\subset F$ est l'anneau des entiers de $F$.
  \end{itemize}
Pour tout $P\in \fc(T_0)$, soit  l'application 
$$H_P: G(F)\to a_P$$
donnée par
\begin{equation}
  \label{eq:def-HP}
  \bg \chi, H_P(g)\bd = \log |\chi(p)|
\end{equation}
où l'on écrit, selon la décomposition d'Iwasawa, $g=pk$ avec $p\in P(F)$ et $k\in K$. Il s'ensuit que l'application $H_P$ est invariante à gauche par $N_P(F)$ et à droite par $K$. 
\end{paragr}

\begin{paragr}[Famille orthogonale.] ---  \label{S:fam-ortho}Soit $M\in \lc(T_0)$. Suivant Arthur (cf. \cite{localtrace} section 3), on dit qu'une famille de points $(Y_P)_{P\in \pc(M)}$ de $a_M$ est \emph{orthogonale} si elle vérifie la condition suivante : pour tout couple $(P,P')\in \pc(M)^2$ adjacents, le vecteur $Y_P-Y_{P'}$ appartient à la droite engendrée par la coracine associée à l'unique élément de
$$ (-\Sigma(P',A_M))\cap\Sigma(P,A_M).$$ 
Si l'on impose de plus que  $Y_P-Y_{P'}$ appartienne à la demi-droite engendrée par ce vecteur, on dit que la famille est \emph{orthogonale positive}. C'est le cas, par exemple, de la famille $(-H_P(g))_{P\in \pc(M)}$ pour tout $g\in G(F)$ (cf. \cite{dis_series}).

Lorsqu'on dispose d'une famille orthogonale $(Y_P)_{P\in \pc(M)}$, on définit pour tout $Q\in \fc(M)$ un point $Y_Q\in a_{M_Q}$ de la manière suivante : on choisit $P\in \pc(M)$ tel que $P\subset Q$ et, par définition, $Y_Q$ est l'image de $Y_P$ par la projection $a_M\to a_{M_Q}$ duale de la restriction $X^*(M_Q)\to X^*(M)$. Les propriétés d'orthogonalité font que le point $Y_Q$ ne dépend pas du choix de $P\subset Q$. Notons que pour $g\in G(F)$, le point  $H_Q(g)$ est bien la projection de $H_P(g)$. Soit $L\in \lc(M)$. La famille ainsi obtenue $(Y_Q)_{Q\in \pc(L)}$ est encore orthogonale.
\end{paragr}

\begin{paragr}[Fonction $R_P$.] --- \label{S:RP}Désormais $M$ désigne le facteur de Levi semi-standard du sous-groupe parabolique $P_0$ qui stabilise le drapeau des noyaux itérés de $X$. Pour $P\in \pc(M)$, on va introduire une légère variante de l'application $H_P$ : pour tout $g\in G(F)$ on pose
$$R_P(g)=H_P(w_P g)$$
où $w_P\in W$ vérifie $w_P^{-1}Pw_P\in \rc(X)$. D'après le lemme \ref{lem:r-std}, un tel élément non seulement existe mais est unique à translation à gauche près par un élément $w\in W^M$. Soit $w_Pg=mnk$ la décomposition d'Iwasawa de $w_Pg$. Pour $w\in W^M$, on a la décomposition d'Iwasawa $ww_P=wmnk$. D'où
$$\bg \chi, H_P(w w_Pg)\bd = \log |\chi(w m)|=\log |\chi(w)|+\log |\chi(m)|=\log |\chi(m)|=\bg \chi, H_P(w w_Pg)\bd,$$
puisque  $|\chi(w)|=1$. Ainsi $R_P$ est indépendante du choix de $w_P$.
\end{paragr}

\begin{paragr}[Orthogonalité de la famille $(-R_P(g))_{P\in \pc(M)}$.] ---\label{S:orth-famille} On montre le lemme suivant.

 \begin{lemme}\label{lem:ortho}
La famille $(-R_P(g))_{P\in \pc(M)}$ est orthogonale au sens d'Arthur.     
  \end{lemme}

  \begin{preuve}
    Il s'agit de vérifier que pour  $P$ et $P'$ adjacents, le vecteur 
$$-R_P(g)+R_{P'}(g)$$
appartient à la droite  engendrée par l'unique élément de $\Delta_P^\vee\cap (-\Delta_{P'}^\vee)$. Cette droite est la droite  $\ago_M^{M_Q}$ pour un unique sous-groupe parabolique $Q$ contenant $P$ et $P'$. Si l'on revient à la définition de $R_P$ et $R_{P'}$, on voit qu'il s'agit de prouver qu'on a 
$$H_Q(w_Pg)=H_{Q}(w_{P'}g).$$
D'après le lemme \ref{lem:leg}, on a 
$$w_{P'}^{-1} \in G_X w_P^{-1} M_Q P'.$$
Comme $w_P^{-1}Pw_P\in \rc(X)$, on a $G_X\subset  w_P^{-1}Pw_P$ d'où
$$w_{P'}^{-1} \in w_P^{-1}Q.$$ 
On a donc
\begin{equation}
  \label{eq:w1w2}
  w_Pw_{P'}^{-1}\in W^{M_Q}.
\end{equation}
On a donc $H_Q(w_Pw_{P'}^{-1})=0$ et 
$$H_{Q}(w_{P}g)=H_Q(w_Pw_{P'}^{-1}   )+H_Q(w_{P'}g)=H_Q(w_{P'}g),$$
ce qu'il fallait démontrer.
\end{preuve}

Soit $g\in G(F)$. Puisque la famille $(R_P(g))_{P\in \pc(M)}$ est orthogonale, on dispose pour tout $Q\in \fc(M)$ de points $R_Q(g)\in a_{M_Q}$. On vérifie qu'on a
\begin{equation}
  \label{eq:HQ-RQ}
  R_Q(g)=H_Q(w_Qg)
\end{equation}
où $w_Q\in W$ vérifie $w_Q^{-1}Qw_Q\in \lc\Sc(X)$ (cf. lemme \ref{lem:r-std}). Comme l'expression $H_Q(w_Qg)$ est indépendante de toute translation à droite de $w_Q$ par un élément de $W^{M_Q}$, elle ne dépend pas du choix de $w_Q$.
\end{paragr}

\begin{paragr}[Action du normalisateur de $M$ dans $W$.] --- \label{S:normM} Soit  $\Norm_W(M)$ le normalisateur de $M$ dans $W$. Pour tout $P\in \pc(M)$, le groupe $P^w=w^{-1}Pw$ appartient encore à $\pc(M)$. On vérifie aisément la formule suivante : pour tout $w\in \Norm_W(M)$, tout $P\in \pc(M)$ et tout $g\in G(F)$ on  a
  \begin{equation}
    \label{eq:normM}
    R_{P}(g)=w\cdot R_{P^w}(g)
  \end{equation}
  où l'on note $w\cdot \, $ l'action naturelle de $\Norm_W(M)$  sur l'espace $a_M$. Plus généralement pour tout $L\in \lc(M)$ et tout $Q\in \pc(L)$, on a 
\begin{equation}
    \label{eq:normMsurL}
    R_{Q}(g)=w\cdot R_{Q^w}(g).
  \end{equation}

\end{paragr}

\begin{paragr}[Action du centralisateur.]  --- On rappelle que $P_0\in \pc(M)$ est le sous-groupe parabolique qui stabilise le drapeau des noyaux itérés de $X$. Le centralisateur de $X$ dans $G$ possède une décomposition de Levi notée $M_XN_X$ où $N_X$ est le radical unipotent de $G_X$ ; la description des groupes $M_X$ et $N_X$ est donnée au §\ref{S:Levi-cent}  (modulo le changement de base à $F$).

  \begin{lemme}
    \label{lem:act-centralisateur}
Soit $h\in G_X(F)$. Pour tout $P\in \pc(M)$, on a 
$$R_P(h)=R_{P_0}(h).$$
En particulier, pour tout $g\in G(F)$ la famille orthogonale $(-R_P(hg))_{P\in \pc(M)}$ se déduit  de la famille  $(-R_P(h))_{P\in \pc(M)}$ par la translation par $-R_{P_0}(h)$.
  \end{lemme}

  \begin{preuve} Soit $P\in \pc(M)$. Par définition, pour tout $g\in G(F)$, on a $R_P(g)=H_P(w_Pg)$ où $w_P\in W$ vérifie $Q:=w_P^{-1}Pw_P\in \rc(X)$. Soit $M_Q$ le facteur de Levi standard de $Q$. On a $M_Q=w_P^{-1}Mw_P$. La conjugaison par $w_P$ induit un isomorphisme $a_{M_Q}\simeq a_M$ indépendant du choix de $w_P$ et pour lequel on a $R_P(g)=w_P\cdot H_Q(g).$
On a donc pour $h\in Q(F)$
\begin{equation}
  \label{eq:translation}
  R_P(hg)=w_P\cdot H_Q(hg)=w_P\cdot H_Q(h)+ w_P\cdot H_Q(g)= R_P(h)+R_P(g).
\end{equation}
Il s'ensuit que la fonction $g\mapsto R_P(g)$ est invariante à gauche par les éléments $q\in Q(F)$ tels que pour tout $\chi\in X^*(P)$, on a $|\chi(w_P qw_P^{-1})|=1$. Le centralisateur $G_X$ est inclus dans $Q$ (cf. remarque \ref{rq:centralisateur}) et admet une décomposition de Levi $M_X N_X$ comme on l'a dit plus haut (cf. aussi §\ref{S:Levi-cent}). La fonction $g\mapsto R_P(g)$ est invariante à gauche par le groupe unipotent $N_X(F)$. On est ramené à prouver le résultat pour un élément $h\in M_X(F)$.

  Par la proposition \ref{prop:drap-R}, il existe $\eps\in \ec$ tel que le sous-groupe parabolique $Q$ soit le stabilisateur du drapeau $E_\bullet(\eps)\otimes_\ZZ F$. Avec les notations de la section \ref{sec:standard}, le groupe $M_Q$ s'identifie alors au produit
  \begin{equation}
    \label{eq:MQ}
    \prod_{1\leq k \leq r}   \Aut_F( \oplus_{1\leq j\leq r}  (V_j^{\al_{k,j}})^{\oplus \eps_{k,j}})
  \end{equation}
  où $\al_{k,j}=\sum_{l=1}^k \eps_{l,j}$. Le groupe $M_X$ s'identifie à
  \begin{equation}
    \label{eq:MX}
    \prod_{1\leq j\leq r}   \Aut_F(V_j^j)
  \end{equation}
Ainsi tout élément  $h\in M_X$ s'identifie à un $r$-uplet $(h_j)_{1\leq j\leq r}$.   On a une flèche évidente 

$$
\prod_{1\leq j\leq r}   \Aut_F(V_j^j) \to \Aut_F( \oplus_{1\leq j\leq r}  (V_j^{\al_{k,j}})^{\oplus \eps_{k,j}})
$$
 qu'on obtient en utilisant les isomorphismes $V_j^i \simeq V_j^{i'}$ induits par les puissances de $X$. L'inclusion $M_X\subset M_Q$ s'identifie alors au plongement \og diagonal\fg{} de \eqref{eq:MX} dans \eqref{eq:MQ}. En utilisant comme base de $X^*(M_Q)$ les déterminants des facteurs de \eqref{eq:MQ}, on identifie $a_{M_Q}$ à $\RR^r$.

De la sorte, on a 
$$H_Q(h)= (\sum_{1\leq j\leq r}    \eps_{k,j}\log|\det(h_j)|)_{1\leq k\leq r}$$
pour tout $h\in M_X(\AAA)$. Rappelons que le sous-groupe parabolique $P_0\in \pc(M)$ est le stabilisateur du drapeau $E_\bullet(\xi)\otimes_\ZZ F$ des noyaux itérés où $\xi$ est défini dans la remarque \ref{rq:particulier}. On a pour $1<k\leq r$
$$E_k(\xi)=E_{k-1}(\xi)\oplus \bigoplus_{k\leq j\leq r} V_j^k.$$
Comme $P$ appartient aussi à $\pc(M)$, il existe une permutation $\sigma\in \mathfrak{S}_r$ telle que $P$ soit   le stabilisateur du drapeau $F_\bullet(\sigma)\otimes_\ZZ F$ défini récursivement par
$$F_k(\sigma)=F_{k-1}(\sigma)\oplus  \bigoplus_{\sigma(k)\leq j\leq r} V_j^{\sigma(k)}.$$
L'endomorphisme $w_P$ de $E$ envoie le drapeau $E_\bullet(\eps)$ sur $F_\bullet(\sigma)$. En particulier, pour tout $k$, on a égalité entre les dimensions des quotients à savoir
\begin{equation}
  \label{eq:dim-quotients}
  \sum_{1\leq j \leq r} \eps_{k,j}d_j= \sum_{ \sigma(k)\leq j\leq r } d_j.
\end{equation}
On en déduit également qu'on a
\begin{equation}
  \label{eq:RP-concret}
  R_P(h)=w_P\cdot H_Q(h)=(\sum_{1\leq j\leq r}    \eps_{\sigma^{-1}(k),j}\log|\det(h_j)|)_{1\leq k\leq r}
\end{equation}
D'après le lemme \ref{lem:transitif}, il existe une permutation $\tau \in  \mathfrak{S}_r$ tel que $\eps_{k,j}=\xi_{\tau(k),j}$ pour tout $1\leq k\leq k$. En utilisant la définition de $\xi$, on voit que les égalités \eqref{eq:dim-quotients} et \eqref{eq:RP-concret} deviennent
\begin{equation}
  \label{eq:dim-quotients2}
  \sum_{\tau(k)  \leq j \leq r} d_j= \sum_{ \sigma(k)\leq j\leq r } d_j
\end{equation}
et
\begin{equation}
  \label{eq:RP-concret2}
  R_P(h)=(\sum_{\tau\sigma^{-1}(k)\leq j\leq r}   \log|\det(h_j)|)_{1\leq k\leq r}.
\end{equation}
Lorsque $\tau$ et $\sigma$ sont égaux à la permutation identique, on trouve
\begin{equation}
  \label{eq:RP-concret3}
  R_{P_0}(h)=(\sum_{k\leq j\leq r}   \log|\det(h_j)|)_{1\leq k\leq r}.
\end{equation}
L'égalité $R_P(h)=R_{P_0}(h)$ résulte alors de l'égalité suivante pour tout $1\leq k \leq r$
\begin{equation}
  \label{eq:triviale} 
  \sum_{\tau\sigma^{-1}(k)\leq j\leq r}   \log|\det(h_j)|=\sum_{k\leq j\leq r}   \log|\det(h_j)|.
\end{equation}
 L'égalité \eqref{eq:dim-quotients2} donne, lorsqu'on remplace $k$ par $\sigma^{-1}(k)$ 
\begin{equation}
  \label{eq:dim-quotients3} 
  \sum_{\tau\sigma^{-1}(k)  \leq j \leq r} d_j= \sum_{ k\leq j\leq r } d_j.
\end{equation}
Donc $d_j=0$ pour $k\leq j <\tau\sigma^{-1}(k)$ ou pour $\tau\sigma^{-1}(k)\leq  j <k$. Pour un tel $j$, on a $h_j=1$ et le terme   $\log|\det(h_j)|$ est nul. L'égalité \eqref{eq:triviale} est donc évidente.

Enfin l'assertion sur la famille orthogonale est alors une conséquence immédiate de \eqref{eq:translation} et de l'égalité $R_P(h)=R_{P_0}(h)$ pour tout $P\in \pc(M)$. 
  \end{preuve}
\end{paragr}

\begin{paragr}[Un calcul pour des sous-groupes paraboliques adjacents.] --- \label{S:non-pos} On rappelle que deux sous-groupes paraboliques $P_1$ et $P_2$ de $\pc(M)$ sont adjacents s'il existe un sous-groupe parabolique $Q$ minimal pour l'inclusion dans $\fc(M)-\pc(M)$  tel que $P_1\subset Q$, $P_2\subset Q$ et $P_1\cap P_2\cap M_Q= M$.

 Soit $P_1$ et $P_2$ deux éléments adjacents de $\pc(M)$. Soit $\tilde{P}_1$ et $\tilde{P}_2$ les éléments de $\rc(X)$ images respectives de $P_1$ et $P_2$ par l'application \eqref{eq:lappli}.

\begin{lemme}\label{lem:adj} Il existe $\eps\in \ec(X)$ et une transposition $\tau\in \mathfrak{S}_r $ tels $\tilde{P}_1$ et $\tilde{P}_2$  soient les stabilisateurs respectifs des drapeaux $E_\bullet(\eps)\otimes_\ZZ F$ et $E_\bullet(\tau(\eps))\otimes_\ZZ F$ (où $\tau(\eps)$ est défini au lemme \ref{lem:transitif}).
  \end{lemme}

  \begin{preuve}
   Soit $Q$ le plus petit sous-groupe parabolique qui contient à la fois $P_1$ et $P_2$. Comme $P_1$ et $P_2$ sont adjacents, $Q$ est aussi minimal parmi les sous-groupes paraboliques qui contiennent strictement $P_1$ (resp. $P_2$). Soit $w_1$ et $w_2$ des éléments de $W$ tels que $w_i^{-1}P_iw_i =\tilde{P}_i$ (cf. lemme \ref{lem:r-std}. Alors $\tilde{Q}=w_1^{-1}Q w_1\in   \lc\Sc(X)$. En fait, on a $w_1w_2^{-1}\in W^{M_Q}$ (cf. éq. \eqref{eq:w1w2} démonstration du lemme \ref{lem:ortho}), on a donc aussi  $\tilde{Q}=w_2^{-1}Q w_2$. Le sous-groupe parabolique $\tilde{Q}$ contient donc  à la fois $\tilde{P}_1$ et $\tilde{P}_2$. Il est également minimal parmi les sous-groupes paraboliques qui contiennent strictement  $\tilde{P}_1$ (resp.  $\tilde{P}_2$) . Par la proposition \ref{prop:drap-R}, il existe $\eps$ et $\eps'$ dans $\ec(X)$ tels que les sous-groupes paraboliques $\tilde{P}_1$ et $\tilde{P}_2$ soient les stabilisateurs respectifs des drapeaux  $E_\bullet(\eps)\otimes_\ZZ F$ et $E_\bullet(\eps')\otimes_\ZZ F$. Par minimalité de $\tilde{Q}$ parmi les sous-groupes paraboliques contenant $\tilde{P}_1$, il existe un unique $1\leq k <r$ tel que le sous-groupe parabolique $\tilde{Q}$ est le stabilisateur du drapeau associé à 
$$(0)=E_0(\eps)\subsetneq \ldots \subsetneq E_{k-1}(\eps)\subsetneq  E_{k+1}(\eps) \subsetneq \ldots \subsetneq E_r(\eps)=E$$
qui se déduit de $E_\bullet(\eps)$ par suppression du sous-espace $E_k(\eps)$. Observons que les groupes $\tilde{P}_1$ et $w_1^{-1}P_2w_1$ appartiennent tous deux à $\pc(w_1^{-1}Mw_1)$, sont inclus dans $\tilde{Q}$ et sont adjacents. Par ailleurs $\tilde{P}_2$ est conjugué par $w_2^{-1}w_1\in W^{M_{\tilde{Q}}}$ à  $w_1^{-1}P_2w_1$. Il s'ensuit que $\tilde{Q}$ est aussi le stabilisateur du drapeau associé à $E_\bullet(\eps')$ privé de $E_k(\eps')$ pour le \emph{même indice} $k$.   Pour tout $1\leq l \leq r$ tel que $l\not=k$ on a donc  $E_l(\eps)=E_{l}(\eps')$ d'où l'on déduit pour tout $j\in J$
$$\sum_{m=1}^l \eps_{l,j}=\sum_{m=1}^l \eps_{l,j}'$$
On obtient alors les conditions suivantes
\begin{equation}
  \label{eq:cond1}
  \eps_{l,j}=\eps_{l,j}'
\end{equation}
pour $l\not=k$ et
\begin{equation}
  \label{eq:cond2}
  \eps_{k,j}+\eps_{k+1,j} =\eps_{k,j}'+\eps_{k+1,j}'. 
\end{equation}
Il reste à vérifier qu'on a $\eps=\tau(\eps')$ pour $\tau$ la transposition qui échange $k$ et $k+1$ c'est-à-dire qu'on a 
\begin{equation}
  \label{eq:laconclusion}
  \left\lbrace
    \begin{array}[l]{l}
       \eps_{k,j}=\eps_{k+1,j}'\\ \eps_{k+1,j}= \eps_{k,j}'
    \end{array}\right.
\end{equation}
Notons d'ailleurs que les deux égalités ci-dessus sont équivalentes par \eqref{eq:cond2}.

Considérons d'abord le cas où  $\eps=\eps'$. Il s'agit de voir qu'on a $\eps_{k,j}=\eps_{k+1,j}$ pour tout $j\in J$. Supposons le contraire par exemple qu'il existe $j_0\in J$ tel que $0=\eps_{k,j_0}$ et $\eps_{k+1,j_0}=1.$ La croissance de $j\mapsto \eps_{k,j}$ implique alors qu'on a   $\eps_{k,j}\leq \eps_{k+1,j}$ pour tout $j\in J$ et l'inégalité est stricte pour $j=j_0$. On a alors
$$
\mathrm{rang}(E_k(\eps)/E_{k-1}(\eps))=\sum_{j\in J} \eps_{k,j}d_j < \sum_{j\in J} \eps_{k+1,j}d_j= 
\mathrm{rang}(E_{k+1}(\eps)/E_{k}(\eps)).
$$
On en déduit en particulier que le normalisateur de $M_{\tilde{P}_1}$ dans $W^{\tilde{Q}}$ est égal à $W^{M_{\tilde{P}_1}}$. Or ce groupe contient évidemment $w_2^{-1}w_1$. Donc $w_2^{-1}w_1\in W^{M_{\tilde{P}_1}}$ d'où 
$$P_1=w_1\tilde{P}_1w_1=w_2\tilde{P}_1w_2= w_2\tilde{P}_2w_2=P_2$$
puisque $\eps=\eps'$ entraîne $\tilde{P}_1=\tilde{P}_2$. Donc $P_1=P_2$ ce qui n'est pas par hypothèse.

Considérons ensuite le cas  $\eps\not=\eps'$. Il existe donc $j_0$ tel que $\eps_{k,j_0}\not=\eps_{k,j_0}'$. Quitte à échanger les rôles de $\eps$ et $\eps'$, on suppose qu'on a  $\eps_{k,j_0}<\eps_{k,j_0}'$ ce qui signifie $\eps_{k,j_0}=0$ et $\eps_{k,j_0}'=1$. En utilisant \eqref{eq:cond2} pour $j=j_0$, on  a aussi $\eps_{k+1,j_0}= 1+\eps_{k+1,j_0}'$ ce qui force  $\eps_{k+1,j_0}=1$ et $\eps_{k+1,j_0}'=0$. On a donc prouvé \eqref{eq:laconclusion} pour $j=j_0$. Par les propriétés de croissance de $j\mapsto \eps_{k,j}$, on a donc pour tout $j<j_0$, $\eps_{k,j}=\eps_{k+1,j}'=0$ et pour tout $j>j_0$, $\eps_{k,j}'=\eps_{k+1,j}=1$ ce qui donne  \eqref{eq:laconclusion} pour $j\not=j_0$.

  \end{preuve}

Reprenons les notations du lemme  \ref{lem:adj} ci-dessus. Les sous-groupes paraboliques $\tilde{P}_1$ et 
$\tilde{P}_2$ sont les stabilisateurs  respectifs des drapeaux   $E_\bullet(\eps)\otimes_\ZZ F$ et  $E_\bullet(\tau(\eps))\otimes_\ZZ F$ où $\tau$ est la transposition qui échange $k$ et $k+1$. Comme on l'a vu au cours de la preuve du lemme  \ref{lem:adj}, quitte à échanger les rôles de $\eps$ et $\tau(\eps)$, on peut et on va supposer  qu'on a   $\eps_{k,j}\leq \eps_{k+1,j}$ pour tout $j\in J$.
Soit
$$J_1=\{j\in J \mid \eps_{k,j}=\eps_{k+1,j}=1\}$$
et
$$J_2=\{j\in J \mid \eps_{k,j}=0 \text{  et  } \eps_{k+1,j}=1\}.$$
On a toujours $r\in J_1$ et l'ensemble $J_2$ est vide si et seulement si $\tilde{P}_1=\tilde{P}_2$. 
 On a  alors $E_k(\eps)\subset E_k(\tau(\eps))$ et les formules suivantes pour les dimensions
 \begin{equation}
   \label{eq:r1}
   r_1:=\mathrm{rang}(E_k(\eps)/E_{k-1}(\eps))= \mathrm{rang}(E_{k+1}(\eps)/E_k(\tau(\eps)))=\sum_{j\in J_1} d_j>0
 \end{equation}
 et
 \begin{equation}
   \label{eq:r2}
   r_2:=\mathrm{rang}(E_k(\tau(\eps))/E_k(\eps))=\sum_{j\in J_2} d_j\geq 0.
 \end{equation}
 
Le raffinement  du drapeau $E_\bullet(\eps)$
\begin{equation}
  \label{eq:drap-raffine}
  (0)\subset E_1(\eps) \subsetneq \ldots  \subsetneq E_{k}(\eps) \subset E_k(\tau(\eps))\subsetneq E_{k+1}(\eps)\subsetneq \ldots  \subsetneq E_r(\eps)
\end{equation}
donne par extention des scalaires un drapeau de sous-espaces dont $\tilde{P}=\tilde{P_1}\cap \tilde{P}_2$ est le stabilisateur. Soit $Q$ le sous-groupe parabolique qui stabilise le drapeau
$$(0)\subset E_1(\eps) \subsetneq \ldots  E_{k-1}(\eps) \subsetneq E_{k+1}(\eps)\subsetneq \ldots E_r(\eps).$$
Par construction des modules $E_k$ comme somme de  $V_j^i$, le drapeau 
$$ E_{k-1}(\eps) \subsetneq E_{k}(\eps) \subset E_k(\tau(\eps))\subsetneq E_{k+1}(\eps)$$
vient avec un scindage $E_{k}(\eps)= E_{k-1}(\eps) \oplus W_1$, $ E_k(\tau(\eps))=E_{k}(\eps)\oplus W_2$ et  $E_{k+1}(\eps)=E_k(\tau(\eps))\oplus W_3$ où les $W_1,\ldots,W_3$ sont eux aussi des sommes de $V_J^i$.
En outre, $X$ induit un isomorphisme de $W_3$ sur $W_1$. On a $\mathrm{rang}(W_1)=\mathrm{rang}(W_3)=r_1$ et $\mathrm{rang}(W_2)=r_2$. Le sous-groupe parabolique ${\tilde{P}}$ est semi-standard. Son facteur de Levi standard $M_{\tilde{P}}$ est muni d'une projection  sur 
$$\Aut_F(W_1)\times \Aut_F(W_2)\times \Aut_F(W_3).$$
On la note $m\mapsto (m_1,m_2,m_3)$.

\begin{lemme}\label{lem:non-pos}
Soit  $g=mnk$ avec $m\in M_{\tilde{P}}(F)$, $n\in N_{\tilde{P}}(F)$ et $k\in K$. Avec les notations ci-dessus, 
$$-R_{P_1}(g)+R_{P_2}(g)= \log|\det(m_1)^{-1}\det(m_3)| \al^\vee$$
où $\al^\vee$ est l'unique élément de $\Delta_{P_1}^\vee\cap (-\Delta_{P_2}^\vee)$.
 \end{lemme}

 \begin{remarque}
   On voit que, sans hypothèse supplémentaire sur $g$, la famille $(-R_P(g))_{P\in \pc(M)}$ n'est pas nécessairement positive.
 \end{remarque}

\begin{preuve} D'après le lemme \ref{lem:ortho}, on sait que l'égalité est vraie à un scalaire près. Le point est d'évaluer ce scalaire. On reprend les notations précédentes, en particulier celles du lemme \ref{lem:adj}. La fonction $g\mapsto R_{P_i}(g)$ est invariante à gauche par les éléments $h\in \tilde{P}_i(F)$ tels que $|\chi(h)|=1$ pour tout $\chi\in X^*(\tilde{P}_i)$. Elle est donc invariante à gauche par $N_{\tilde{P}}$. Elle est par ailleurs invariante à droite par $K$. On 
 On est ramené à prouver l'égalité pour $n$ et $k$ triviaux. On a également $R_{P_i}=w_i\cdot H_{\tilde{P}_i}$.
On est donc ramener à prouver pour tout $m\in M_P(F)$ l'égalité
\begin{equation}
  \label{eq:verif1}
  -H_{\tilde{P}_1}(m)+(w_1^{-1}w_2)\cdot H_{\tilde{P}_2}(m)= \log|\det(m_1)^{-1}\det(m_3)| \be^\vee
\end{equation}
où $\beta^\vee$ est l'unique élément de $\Delta_{\tilde{P}_1}^{\tilde{Q},\vee}$.  Le groupe $M_{\tilde{P}_1}$ est muni naturellement d'une projection sur $\Aut_F(W_1)$ et on va composer cette projection avec le déterminant. Cela donne un caractère $\chi$ de $X^*(\tilde{P}_1)$ pour lequel $\bg \chi , \be^\vee\bd=1$ et $\bg \chi , H_{Q_1}(m)\bd =\log|\det(m_1)|$. Comme \eqref{eq:verif1} est au moins vraie à une constante près,  il suffit de la vérifier sur le caractère $\chi$.

Pour continuer, il nous faut comprendre l'action de $w_1^{-1}w_2$. On a déjà remarqué dans la preuve du lemme \ref{lem:adj} que le sous-groupe  $w^{-1}_1P_2w_1$ appartient à $\pc(M_{\tilde{P}_1})$, qu'il est contenu dans $\tilde{Q}$ et qu'il est adjacent à $\tilde{P}_1$. Ainsi $w^{-1}_1P_2w_1$ est nécessairement le stabilisateur du drapeau associé à 

$$ (0)\subset E_1(\eps) \subsetneq \ldots  \subsetneq E_{k-1}(\eps) \subset E_k(\eps)\oplus W_2\oplus W_3 \subsetneq E_{k+1}(\eps)\subsetneq \ldots  \subsetneq E_r(\eps)$$

Le groupe $\tilde{P}_2$ est le stabilisateur du drapeau associé à 
$$ (0)\subset E_1(\eps) \subsetneq \ldots  \subsetneq E_{k-1}(\eps) \subset E_k(\eps)\oplus W_1\oplus W_2 \subsetneq E_{k+1}(\eps)\subsetneq \ldots  \subsetneq E_r(\eps).$$
L'élément $w\in W^{M_{\tilde{Q}}}$ d'ordre $2$ qui échange $W_1$ et $W_3$ (en induisant $X$ sur $W_3$) conjugue  $\tilde{P}_2$ en  $w_1^{-1}P_2 w_1$ tout comme d'ailleurs  $w_1^{-1}w_2$. Ces éléments sont donc égaux à un élément de $W^{M_{\tilde{P}_2}}$ près qui, de toute façon, agit trivialement sur $a_{\tilde{P}_2}$. Le caractère $\chi\circ (w_1^{-1}w_2)$ est alors celui obtenu par composition du déterminant avec la projection $M_P\to \Aut_F(W_3)$. Le lemme devient évident.

\end{preuve}

La décomposition \eqref{eq:decomp} induit une projection $\ggo \to  \Hom(W_3,W_1)\otimes_\ZZ F$ dont on note $U\mapsto U_{1,3}$ la restriction à $\ngo_P$. Soit $\det(U_{1,3})$ le déterminant de la matrice de $U_{1,3}$ dans les bases de $W_1$ et $W_3$ extraites de  la base $e$ de $E$ fixée au § \ref{S:surZ}. Observons qu'on a
\begin{equation}
  \label{eq:decomp-n1-n2}
  \ngo_{\tilde{P}_1}\cap\ngo_{\tilde{P}_2}= \Hom(W_3,W_1)\oplus \ngo_{\tilde{Q}}.
\end{equation}

Le lemme suivant fournit une autre interprétation de la différence $-R_{P_1}(g)+R_{P_2}(g)$.

\begin{lemme}\label{lem:Y}
  Soit $Y\in \ggo(F)$ dans l'orbite de $X$ sous $G(F)$. 
  \begin{enumerate}
  \item Il existe $k\in K$ tel que 
$$U=k Y k^{-1} \in \ngo_{\tilde{P}_1}(F)\cap \ngo_{\tilde{P}_2}(F).$$
L'élément $|\det(U_{1,3})| \in \RR $ ne dépend pas du choix de $k$.
\item Pour tout $g\in G(F)$ tel que $Y=g^{-1}X g$, on a 
$$-R_{P_1}(g)+R_{P_2}(g)= \log|\det(U_{1,3})|\al^\vee,$$
où $\al^\vee$ est l'unique élément de $\Delta_{P_1}^\vee\cap (-\Delta_{P_2}^\vee)$.
\end{enumerate}
\end{lemme}

\begin{preuve}
  Afin d'alléger les notations durant la démonstration,  on pose $\ngo_i=\ngo_{\tilde{P}_i}$ pour $i=1,2$.
On a $Y=g^{-1}X g$ pour un certain $g\in G(F)$. Par décomposition d'Iwasawa, on a $g=pk$ avec $p\in \tilde{P}(F)$. Il s'ensuit que $U=  k Y k^{-1}$ appartient à la $\tilde{P}(F)$-orbite de $X$. Or $\tilde{P}=\tilde{P}_1 \cap \tilde{P}_2$ et la $\tilde{P}_i$-orbite de $X$ est incluse dans $\ngo_i$. On a donc $U\in  \ngo_{1}(F)\cap \ngo_2(F)$. Comme $G_X\subset \tilde{P}$, la $\tilde{P}$-orbite de $X$ est l'intersection des orbites de $X$ sous $\tilde{P}_1$ et $\tilde{P}_2$. C'est aussi l'intersection de $\ngo_1\cap \ngo_2$ avec la $G$-orbite de $X$. Donc si $k\in K$, vérifie 
$$kg^{-1}X gk^{-1}\in \ngo_{1}(F)\cap \ngo_2(F)$$
on en déduit qu'il existe $p\in \tilde{P}(F)$ tel que 
$$kg^{-1}X gk^{-1}=p^{-1}X p$$
donc $g \in G_X(F) pk\subset \tilde{P}(F)k$. Alors $k$ est nécessairement la composante sur $K$ de la décomposition d'Iwasawa de $g$ relative au sous-groupe parabolique $ \tilde{P}$. Il s'ensuit que $k$ est bien défini à une translation à gauche près par un élément de $K\cap  \tilde{P}(F)$.

 Pour tout élément $p\in  \tilde{P}$ de projection $(m_1,m_3)$ sur $\Aut_F(W_1)\oplus\Aut_F(W_3)$ et tout $U \in \ngo_1\cap \ngo_2$, on a 
$$(p^{-1}U p)_{1,3}= m_1^{-1}U_{1,3}m_3$$
et donc
\begin{equation}
  \label{eq:equiv}
  \det( (p^{-1}U p)_{1,3})= \det(m_1)^{-1}\det(U_{1,3})\det(m_3).
\end{equation}

Si, de plus, $p\in K\cap   \tilde{P}(F)$, on a $|\det(m_1)|=|\det(m_3)|=1$ d'où l'indépendance vis-à-vis du choix de $k$ dans l'assertion $1$.

Prouvons l'assertion 2. Soit $g=pk$ la décomposition d'Iwasawa de $g$ selon $G(F)=  \tilde{P}(F)K$. Dans ce cas, on a $U=p^{-1}Xp$. Il vient d'après \eqref{eq:equiv} et l'égalité $\det(X_{1,3})=1$
$$\det(U_{1,3})=\det((p^{-1}X p)_{1,3}))= \det(m_1)^{-1}\det(X_{1,3})\det(m_3)= \det(m_1)^{-1}\det(m_3) .$$
L'assertion 2 résulte alors du lemme \ref{lem:non-pos}.
\end{preuve}

\end{paragr}

\section{Des $(G,M)$-familles associées à un endomorphisme nilpotent}\label{sec:GMfam}

\begin{paragr}[Rappel sur les $(G,M)$-familles.] ---\label{S:GMfam}
On reprend les notations de la section \ref{sec:Richard}. Soit $M$ un sous-groupe de Levi de $G$. À la suite d'Arthur (cf. \cite{trace_inv}), on appelle $(G,M)$-famille toute famille $(c_P(\la))_{P\in \pc(M)}$ de fonctions holomorphes de la variable $\la\in a_{M,\CC}^*$ qui vérifient la condition de recollement suivante : pour couple $(P,P')$ d'éléments adjacents de $\pc(M)$, on a $c_P(\la)=c_{P'}(\la)$ sur l'hyperplan défini par $\bg \la, \al^\vee \bd=0$ où $\al^\vee$ est l'unique élément de $\Delta_P^{\vee}\cap (-\Delta_{P'}^\vee)$.

Soit $L\in \lc(M)$. On associe à une $(G,M)$-famille $c$ la $(G,L)$-famille $d$ définie ainsi :  pour $Q\in \pc(L)$ et  $\la\in a_{L,\CC}^*$, on pose $d_Q(\la)=c_P(\la)$ pour n'importe quel $P\in \pc(M)$ qui vérifie $P\subset Q$ (le résultat est indépendant du choix de $P$).

Soit $Q\in \pc(L)$. On associe à une $(G,M)$-famille $c$ la $(L,M)$-famille $c^Q$ définie ainsi : pour tout $P\in \pc^L(M)$  et  $\la\in a_{M,\CC}^*$, on pose $c_P^Q(\la)=c_{PN_Q}(\la)$.  
\end{paragr}

\begin{paragr}[Orbite nilpotente et $(G,M)$-famille : la construction d'Arthur.] --- \label{S:ctr-Arthur} Désormais  $F$ est un corps local. On utilise les notations du §\ref{S:HP}. Soit $P$ un sous-groupe parabolique semi-standard quelconque de $G$ de facteur de Levi noté $M$. Soit $\of\subset \mgo$ une orbite nilpotente sous l'action de $M$. Soit $\ago_M$ l'algèbre de Lie du centre $A_M$ de $M$ et $\ago_M^{\reg}$ l'ouvert de $\ago_M$ défini par la condition 
$$\prod_{\al \in \Sigma(G,A_M)}\al \not=0.$$
Pour alléger les notations, on pose $N=N_{P}$ et $\ngo=\ngo_{P}$. 
Dans la suite de ce paragraphe, on va rappeler des constructions dues à Arthur (cf. \cite{wei_or} §3), pour tout $Z\in \ago_M^{\reg}$, tout $U\in \of(F)$ et tout $Y\in \ngo$, il existe un unique 
$$n=n(Z,U,Y)\in N$$ 
défini par la condition
\begin{equation}
  \label{eq:fonda}
  n^{-1}(Z+U)n=Z+U+Y.
\end{equation}

Soit $\al\in \Sigma(G,A_M)$. Soit $L\in \lc(M)$ tel que $ \Sigma(L,A_M)=\{\pm\al\}$. On suppose que $P\cap L$ est le sous-groupe parabolique de $L$ tel que $\Sigma(P\cap L,A_M)= \{-\al\}$. Soit $P'\in\pc(M)$ tel que $P'\cap L$ soit opposé à $P\cap L$. Pour tout  $Y\in \ngo\cap  \lgo$, on  a $n(Z,U,Y) \in N\cap L$ pour tout $Z\in  \ago_M^{\reg}$ et $U\in \of$. En particulier, $H_{P'}(n(Z,U,Y))$ ne dépend pas du choix de $P'$. 

Il existe un unique réel positif $\rho(\al,\of)$ tel que si l'on pose 
$$r_\al(\la,\of,Z)= |\al(Z)|^{ \rho(\al,\of) \bg \la, \al^\vee\bd}$$ 
pour $\la\in a_{M,\CC}^*$,  la limite
\begin{equation}
  \label{eq:lalimiteexiste}
  \lim_{Z \in \ago_M^{\reg} \to 0} r_\al(\la,\of,Z) \exp(-\bg \la, H_{P'}(n(Z,U,Y))\bd)
\end{equation}
existe et définit une fonction non identiquement nulle des variables $U\in \of(F) $ et  $Y\in \ngo(F)\cap \lgo(F)$. Cette propriété caractérise  $\rho(\al,\of)$. On a aussi $\rho(\al,\of)=\rho(-\al,\of)$.

À la suite d'Arthur (\emph{ibid.}), on introduit la famille suivante, qui est une  $(G,M)$-famille,
\begin{equation}
  \label{eq:wP}
  w_{P'}(\la,Z,U,Y)= \big(\prod_{\al \in \Sigma(P',A_M)\cap    \Sigma(\bar{P},A_M) }  r_\al(\la,\of,Z) \big) \exp(-\bg \la, H_{P'}(n(Z,U,Y))\bd)
\end{equation}
pour $P'\in \pc(M)$, $Z\in  \ago_M^{\reg}(F)$, $U\in \of(F)$ et $Y\in \ngo(F)$. On a noté $\bar{P}\in \pc(M)$ le sous-grope parabolique opposé à $P$. Cette famille dépend du sous-groupe parabolique \og point-base\fg{} $P$ mais, pour ne pas alourdir encore les notations, on ne le fait pas figurer dans la notation. 

Soit $I_{P}^G(\of)$ l'orbite induite de $\of$ au sens de Lusztig-Spaltenstein. Pour tout $U\in \of(F)$ et tout $Y\in  \ngo(F)$ tels que $U+Y\in I_{P}^G(\of)$ la limite
$$ w_{P'}(\la,U,Y):= \lim_{Z \in \ago_M^{\reg} \to 0}  w_{P'}(\la,Z,U,Y)$$
existe et est non-nulle (cf. \cite{wei_or} lemme 4.1). Cela définit une autre $(G,M)$-famille  
\begin{equation}
  \label{eq:wP'}
  (w_{P'}(\la,U,Y))_{P'\in \pc(M)}.
\end{equation}
\end{paragr}

\begin{paragr}\label{S:wP'} Dans le lemme suivant, on considère l'orbite nulle dans $\mgo$. On a forcément  $U=0$ et on note simplement 
 $$(w_{P'}(\la,Z,Y))_{P'\in \pc(M)} $$
 et  
 \begin{equation}
   \label{eq:wPbis}
   (w_{P'}(\la,Y))_{P'\in \pc(M)}
 \end{equation}
les  $(G,M)$-familles \eqref{eq:wP} et  \eqref{eq:wP'} pour $U=0$ et tout $Y\in \ngo(F)\cap I_{P}^G(0)$.

Soit $P_1$ un sous-groupe parabolique contenant $P$ et soit $M_1=M_{P_1}$, $N_1=N_{P_1}$. Soit $\of_1=I_{P\cap M_1}^{M_1}(0)$ : c'est une orbite nilpotente dans $\mgo_1$ pour l'action de $M_1$. Relativement au sous-groupe parabolique \og point-base\fg{} $P_1$, on dispose de la $(G,M_1)$-famille $(w_{P_1'}(\la,U,Y))_{P_1'\in \pc(M_1)}$ définie pour $U\in \of_1(F)$ et  $Y\in \ngo_1(F)$ tels que $U+Y\in I_{P_1}^G(\of_1)$. Supposons de plus qu'on a $U\in \ngo(F)$. Alors $U+Y\in \ngo(F)$ et par transitivité de l'induction de Lusztig-Spaltenstein, on a $U+Y\in I_{P}^G(0)$. Donc on dispose aussi de la $(G,M)$-famille  $(w_{P'}(\la,U+Y))_{P'\in \pc(M)}$ pour le \og point-base\fg{} $P$.

\begin{lemme}  \label{lem:descente-GM}Sous les hypothèses ci-dessus, la $(G,M_1)$-famille déduite de  la $(G,M)$-famille $(w_{P'}(\la,U+Y))_{P'\in \pc(M)}$ coïncide avec la $(G,M_1)$-famille  $(w_{P_1'}(\la,U,Y))_{P_1'\in \pc(M_1)}$. En d'autres termes, on a
  \begin{equation}
    \label{eq:descente}
    w_{P'}(\la,U+Y)=w_{P_1'}(\la,U,Y),
  \end{equation}
  pour tout $\la  \in a_{M_1,\CC}^* $ et tous sous-groupes paraboliques $P'\in \pc(M)$ et $P_1'\in \pc(M_1)$ tels que $P'\subset P_1'$.
\end{lemme}

\begin{preuve}
  Soit $Z\in \ago_M^{\reg}$.  Pour tout $U\in \of_1 \cap \ngo$ et $Y\in \ngo_1$, soit $n\in N$ tel que 
  $$n^{-1}Zn= Z+U+Y.$$
Par définition (cf. éq. \eqref{eq:wP}), on a pour $P'\in \pc(M)$ et $\la\in a_{M,\CC}^*$
\begin{equation}
  \label{eq:defbrute}
   w_{P'}(\la,Z,0,U+Y)= \big(\prod_{\al \in \Sigma(P',A_M)\cap    \Sigma(\bar{P},A_M) }  r_\al(\la,(0),Z) \big) \exp(-\bg \la, H_{P'}(n)\bd).
\end{equation}
On peut spécialiser cette égalité en $\la\in  a_{M_1,\CC}^*$.  Écrivons $n=m_1n_1$ avec $m_1\in M_1\cap N$ et $n_1\in N_1$. Pour $\la\in  a_{M_1,\CC}^* $ et $P_1'\supset P'$, on a 
\begin{equation}
  \label{eq:egalite-acc}
  \bg \la, H_{P'}(n)\bd=\bg \la, H_{P_1'}(n)\bd= \bg \la, H_{P_1'}(n_1)\bd
\end{equation}
et pour $\al\in \Sigma(M_1,A_M)$
\begin{equation}
  \label{eq:vaut1}
   r_\al(\la,(0),Z)=1.
\end{equation}
  Soit $U'=m_1^{-1}Zm_1 -Z\in \mgo_1\cap \ngo$. On a alors
\begin{eqnarray*}
  n^{-1}Zn&=& n_1^{-1}(U'+Z)n_1\\
&=& Z+ U' +  n_1^{-1}(U'+Z)n_1 -(U'+Z)\\
&\in & Z+U' +\ngo_1.
\end{eqnarray*}
Il s'ensuit qu'on a $U=U'$ et 
$$n_1^{-1}(Z+U) n_1 = Z+U+Y.$$
Cette égalité caractérise $n_1$ lorsque $Z$ est dans l'ouvert $\ago_M'\subset \ago_M$ défini par 
$$\prod_{\al\in \Sigma(N_1,A_M)}\al\not=0.$$
On obtient un morphisme 
$$(Z,U,Y)\in \ago_M'\times (\of_1\cap \ngo)\times \ngo_1 \mapsto n_1(Z,U,Y)\in N_1.$$
Observons que l'ouvert $\ago_M'$ contient  l'ouvert $\ago_{M}^{\reg}$ et la partie localement fermée $\ago_{M_1}^{\reg}$. En tenant compte de \eqref{eq:defbrute}, \eqref{eq:egalite-acc} et \eqref{eq:vaut1}, on obtient que   pour $P'\in \pc(M)$ et tous $Z\in \ago_{M}^{\reg}(F)$,  $U\in \of_1(F) \cap \ngo(F)$, $Y\in \ngo_1(F)$ et $\la\in a_{M_1,\CC}^*$ on a 
\begin{equation}
  \label{eq:complique}
 w_{P'}(\la,Z,0,U+Y)= \big(\prod_{\al \in \Sigma(P'_1,A_M)\cap    \Sigma(\bar{P}_1,A_M) }  r_\al(\la,(0),Z) \big) \exp(-\bg \la, H_{P_1'}(n_1(Z,U,Y))\bd).
\end{equation}
où $P'_1\in \pc(M_1)$ est tel que $P'\subset P_1'$.
D'après Arthur (cf. \cite{wei_or}, lemmes 4.1 et 4.2 ainsi que leur preuve), l'expression  $w_{P'}(\la,Z,0,U+Y)$, \emph{a priori} définie seulement pour $Z\in \ago_{M}^{\reg}(F)$, se prolonge par continuité en une fonction non nulle sur  un ouvert de  $  \ago_{M}(F) \times (\of_1(F) \cap \ngo(F))\times \ngo_1(F)$ qui contient les triplets $(0,U,Y)$ tels que $Y+U\in I_{P_1}^G(\of_1)=I_{P}^G(0)$. Le membre de droite de \eqref{eq:complique} est, lui, bien défini dès que $Z\in \ago_M'(F)$, en particulier si $Z\in \ago_{M_1}^{\reg}$. Considérons le cas où $P_1$ et $P_1'$ sont adjacents. Soit $\beta\in \Sigma(P_1',A_L)\cap -\Sigma(P_1,A_L)$. Pour $Z\in \ago_M'$, le membre de droite de \eqref{eq:complique} se réduit à 
$$|\beta(Z) |^{\rho(\beta) \bg \la, \beta^\vee\bd}  \exp(-\bg \la, H_{P_1'}(n_1(Z,U,Y))\bd)$$
avec
$$\rho(\beta)=\sum_{\al \in \Sigma(G,A_M),\, \al_{|A_L }=\beta} \rho(\al,0).$$
La caractérisation du réel $\rho(\beta,\of_1)$ par l'existence et la non trivialité de la limite \eqref{eq:lalimiteexiste} implique qu'on a $\rho(\beta,\of_1)=\rho(\beta)$. En faisant varier les paraboliques $P$ et $P'$, on en déduit que cette égalité est vraie pour tout $\beta$. On en déduit qu'on a pour $Z\in A_{M_1}^{\reg}(F)$
$$\prod_{\al \in \Sigma(P'_1,A_M)\cap    \Sigma(\bar{P}_1,A_M) }  r_\al(\la,(0),Z) =\prod_{\al \in \Sigma(P'_1,A_{M_1})\cap    \Sigma(\bar{P}_1,A_{M_1}) }  r_\al(\la,\of_1,Z)$$
Maintenant le membre de droite de  \eqref{eq:complique} s'interprète pour $Z\in A_{M_1}^{\reg}(F)$ comme $w_{P_1'}(\la,Z,U,Y)$. Prenant la limite quand $Z\to 0$ dans chacun des membres de \eqref{eq:complique}, on obtient l'égalité cherchée \eqref{eq:descente}. 

\end{preuve}
  
\end{paragr}

\begin{paragr}[Une $(G,M)$-famille associée à la famille orthogonale $(R_P(g))_{P\in \pc(M)}$.]  --- Les notations sont celles  de la section \ref{sec:fam-ortho} et des paragraphes ci-dessus. En particulier, on rappelle que $X$ est un élément nilpotent de $\ggo(F)$ \og standard\fg{}. Le groupe $M$ est ici le facteur de Levi du stabilisateur $P_0$ du drapeau des noyaux itérés de $X$. Pour tout $g\in G(F)$, on dispose de la famille orthogonale $(R_{P}(g))_{P\in \pc(M)}$. Fixons un sous-groupe parabolique \og point-base\fg{} $P_1\in \pc(M)$. Les propriétés des familles orthogonales font que la famille
$$(\exp(\bg \la, R_{P_1}(g)-R_{P}(g)\bd))_{P\in \pc(M)}$$
est une $(G,M)$-famille. Rappelons que puisque $P_1\in \pc(M)$, l'induite $I_{P_1}^G(0)$ n'est autre que l'orbite de $X$ sous $G(F)$. Pour tout $Y\in I_{P_1}^G(0)\cap \ngo_1(F)$, on dispose aussi de la $(G,M)$-famille  $(w_{P}(\la,Y))_{P\in \pc(M)}$, cf. §\ref{S:wP'}, définie relativement au même point-base $P_1$. La proposition suivante relie ces deux $(G,M)$-familles. Pour tout indice $i$,  on abrège  $\ngo_{P_i}$ en $\ngo_i$.

\begin{proposition}\label{lem:comparaisonGM}
  Pour tous $g\in G(F)$,  $Y \in \ngo_1(F)\cap I_{P_1}^G(0)$, et $k\in K$ tels que 
\begin{equation}
  \label{eq:egality}
  k^{-1}Yk= g^{-1}Xg
\end{equation}
on a pour tout $P_2\in \pc(M)$
\begin{equation}
  \label{eq:aprouver}
   w_{P_2}(\la,Y)=\exp(\bg \la, R_{P_1}(g)-R_{P_2}(g)\bd),
 \end{equation}
où, l'on rappelle que $P_1$ est le  \og  point-base \fg{} de la famille $(w_P)_{P\in \pc(M)}$.  

En particulier, le membre de gauche est invariant par conjugaison de $Y$ par un élément de $P(F)\cap K$.
\end{proposition}

\begin{preuve}
La dernière assertion est une conséquence immédiate de l'invariance de la formule \eqref{eq:aprouver} par toute translation à droite de $g$ par un élément de $K$. 

Prouvons \eqref{eq:aprouver}.
Le cas $P_1=P_2$ est évident. Soit  $\mathrm{dist}(P_1,P_2)$, la distance de $P_1$ à $P_2$, c'est-à-dire la longueur minimale d'une chaîne de sous-groupes paraboliques adjacents de $\pc(M)$ qui relie $P_1$ à $P_2$. On va d'abord montrer par récurrence sur $\mathrm{dist}(P_1,P_2)$ qu'il suffit de prouver la proposition dans le cas où $P_1$  et $P_2$ sont adjacents, cas qu'on admet provisoirement. On suppose   $\mathrm{dist}(P_1,P_2)\geq 2$ et soit $P_3\in \pc(M)$ adjacent à $P_2$ tel que
$$\mathrm{dist}(P_1,P_2)=\mathrm{dist}(P_1,P_3)+1.$$

Pour $i\in \{1,2,3\}$,  on introduit l'espace topologique
$$\mathcal{S}_i=\{(k,H)\in (K\cap P_i(F))\back K \times \ggo(F) \mid  kHk^{-1} \in \ago_M(F)\oplus \ngo_i(F)    \}$$
Cette ensemble est muni d'une application continue $p_i: \mathcal{S}_i \to \ago_M(F)  \times \ggo(F)$ définie ainsi : l'image de $((K\cap P_i(F))k, H)$ est donnée par $(Z,H)$ tel que $kHk^{-1} \in Z\oplus \ngo_i(F)$. Soit   $\mathcal{S}_i^{\reg}$ l'image réciproque de l'ouvert $\ago_M^{\reg}(F)\oplus \ngo_i(F)$. D'après Arthur (\cite{wei_or} preuve du lemme 4.2), il existe une unique application continue 
\begin{equation}
  \label{eq:fcontinue}
  f: \mathcal{S}_1^{\reg} \to \mathcal{S}_3^{\reg}
\end{equation}
telle que $p_3\circ f=p_1$. On peut décrire cette application ainsi : soit $((K\cap P_1(F))k, H)\in \mathcal{S}_1^{\reg}$. Soit $Z\in \ago_{M}^{\reg}(F)$ et $n_1\in N_1(F)$ les éléments uniquement définis par la condition : 
$$kHk^{-1}=n^{-1}_1Zn_1.$$
Soit $n_1=mn_3 k_3$ la décomposition d'Iwasawa de $n_1$ selon $G(F)=M(F)N_3(F) K$. L'élément $k_3$ est bien défini modulo $K\cap P_3(F)$. On a $f((K\cap P_1(F))k, H))= ((K\cap P_3(F))k_3,H)$. On vérifie aisément que cette application est indépendante du représentant $k$ choisi.

Soit $Y\in \ngo_1(F)$ et $Z\in \ago_{M}^{\reg}$. Le couple $(K\cap P_1(F), Z+Y)$ appartient à  $\mathcal{S}_1^{\reg}$. Soit $((K\cap P_3(F))k_3, Z+Y)$ son image par $f$ où $k_3\in K$. Soit $Y_3= k_3(Z+Y)k_3^{-1}-Z$ ; c'est un élément de $\ngo_3(F)$. On a alors l'identité suivante (cf. \cite{wei_or} preuve du lemme 4.1)
\begin{equation}
  \label{eq:produit}
  w_{P_2}(\la, Z,Y) = \tilde{w}_{P_2}(\la,Z,Y_3)\cdot w_{P_3}(\la, Z,Y)
\end{equation}
où l'on note  $(\tilde{w}_{P'}(\la,Z,Y_3))_{P'\in \pc(M)}$ la $(G,M)$-famille analogue à  $(w_{P'}(\la,Z,Y))_{P'\in \pc(M)}$ mais définie relativement au point-base $P_3$. D'après Arthur (cf. \cite{wei_or} lemme 4.2 et sa preuve), l'application \eqref{eq:fcontinue} a un prolongement continu à un ouvert de $\Sc_1$ qui contient 
$\{K\cap P_1(F)\}\times \ngo_1(F)\cap I_{P_1}^G(0)$. De plus, l'image par ce prolongement d'un couple $(K\cap P_1(F),Y)$ dans cet ensemble  est $(Kk,Y)$ où $k$ est un élément quelconque de $K$ qui vérifie $kYk^{-1}\in \ngo_3(F)$. Lorsque, de plus,  $Y\in \ngo_1(F)\cap I_{P_1}^G(0)$, on peut prendre la limite quand $Z\to 0$ dans \eqref{eq:produit} pour obtenir

\begin{equation}
  \label{eq:produit2}
  w_{P_2}(\la,Y) = \tilde{w}_{P_2}(\la,Y_3)\cdot w_{P_3}(\la, Y)
\end{equation}
où $Y_3$ est un élément  de $\ngo_3(F)$ conjugué sous $K$ à $Y$. Supposons, de plus, que $Y$ est $K$-conjugué à $g^{-1}Xg$. Il en est donc de même de  $g^{-1}Xg$. En particulier, on a aussi $Y_3\in \ngo_3(F)\cap I_{P_3}^G(0)$. Par hypothèse de récurrence, on a  alors
\begin{equation}
  \label{eq:fait1}
  w_{P_3}(\la,Y)=\exp(\bg \la, R_{P_1}(g)-R_{P_3}(g)\bd)
\end{equation}
et  
\begin{equation}
  \label{eq:fait2}
  \tilde{w}_{P_2}(\la,Y_3)=  \exp(\bg \la, R_{P_3}(g)-R_{P_2}(g)\bd).
\end{equation}
En combinant \eqref{eq:produit2}, \eqref{eq:fait1} et \eqref{eq:fait2}, on obtient bien la formule cherchée \eqref{eq:aprouver} pour le couple $(P_1,P_3)$.

Il nous reste à prouver \eqref{eq:aprouver} lorsque $P_1$ et $P_2$ sont adjacents. On reprend les notations utilisées au §\ref{S:non-pos}. On utilise les éléments $\tilde{P}_1$ et $\tilde{P}_2$ de $\rc(X)$ images respectives de $P_1$ et $P_2$ par \eqref{eq:lappli}. Pour $i=1,2$, soit $w_i\in W$ tel que $w_i^{-1}P_i w_i= \tilde{P_i}$. On rappelle que $Q$ est le plus petit sous-groupe parabolique qui contient $P_1$ et $P_2$. Soit $\tilde{Q}$ l'image de $Q$ par \eqref{eq:lappli2} et soit $\tilde{P}=\tilde{P}_1\cap \tilde{P}_2$ : ce dernier groupe est un sous-groupe parabolique.

Soit $Y\in \ngo_1(F)\cap I_{P_1}^G(0)$, $k\in K$ et $g\in G(F)$ tels que   $k^{-1}Yk= g^{-1}Xg$. Le second membre de \eqref{eq:aprouver} est invariant lorsqu'on translate  $g$ par un élément de $K$. Par décomposition d'Iwasawa, on est ramené au cas où $g\in \tilde{P}(F)$, ce qu'on suppose désormais. On a alors $\Ad(g^{-1})X\in \ngo_{\tilde{P}_1(F)} \cap \ngo_{\tilde{P}_2}(F)$ et  $\Ad(w_1g^{-1})X\in \ngo_1(F)\cap I_{P_1}^G(0)$. Comme $Y$ vérifie la même propriété, il existe $p_1\in P_1(F)$ tel que $\Ad(p_1w_1g^{-1})X=Y=\Ad(kg^{-1})X  $. On en déduit 
$$g w_1^{-1} p_1^{-1}k g^{-1}\in G_X(F)\subset \tilde{P_1}(F)$$
d'où l'on tire $w_1^{-1}k \in \tilde{P_1}(F)\cap K$. Quitte à translater $g$ à droite par l'inverse de cet élément, on peut et va supposer qu'il est trivial c'est-à-dire qu'on a $k=w_1$.
Soit $Z\in \ago_M^{\reg}(F)$ et $n_1\in N_1(F)$ tel que 
$$\Ad(n^{-1}_1)Z=Z+Y.$$
Soit $Y=Y_Q+Y^Q$ la décomposition selon $\qgo=\mgo_Q\oplus\ngo_Q$. On écrit $n_1=m_Q n_Q$ avec $m_Q\in M_Q(F)\cap N_1(F)$ et $n_Q\in N_Q(F)$. Alors on a 
$$Y^Q=\Ad(n^{-1}_Q)( \Ad(m_Q^{-1})Z)-   \Ad(m_Q^{-1})Z$$
et 
\begin{equation}
  \label{eq:mQ}
  Z+Y_Q= \Ad(m_Q^{-1})Z.
\end{equation}
Ces éléments sont les projections respectives de $Z+Y$ sur $\ngo_Q$ et $\mgo_Q$. Soit $\al$ l'unique élément de $\Sigma(P_2,A_M)\cap \Sigma(\bar{P}_1,A_M)$. On pose $\rho(\al)=\rho(\al,(0))$ où $(0)$ est ici l'orbite nulle de $\mgo_1$.  Par ailleurs, on a
\begin{eqnarray*}
  w_{P_2}(\la,Z,Y)&=&|\al(Z)|^{\rho(\al) \bg \la, \al^\vee\bd } \exp(-\bg \la, H_{P_2}(n_1)\bd)\\
&=& |\al(Z)|^{\rho(\al) \bg \la, \al^\vee\bd } \exp(-\bg \la, H_{P_2}(m_Q)\bd)
\end{eqnarray*}
 Soit $\tilde{Z}=\Ad(w_1^{-1})Z$ et $\tilde{Y}= \Ad(w_1^{-1})Y_Q$. Soit $\tilde{m}_Q = w_1^{-1}m_Q w_1$. Alors $\tilde{Z}\in \ago_{\tilde{M}}^{\reg}(F)$ où l'on pose $\tilde{M}=M_{\tilde{P}_1}$ et $\tilde{Y}\in \mgo_{\tilde{Q}}(F)\cap \ngo_{\tilde{P}_1}(F)$. En outre $\tilde{m}_Q \in M_{\tilde{Q}}(F)\cap N_{\tilde{P}_1}(F)$. On a évidemment
\begin{equation}
  \label{eq:mQ2}
  \tilde{Z}+\tilde{Y}= \Ad(\tilde{m}_Q^{-1})\tilde{Z}
\end{equation}
et cette égalité caractérise $\tilde{m}_Q$ dans $M_{\tilde{Q}}(F)\cap N_{\tilde{P}_1}(F)$ en fonction de $\tilde{Z}$ et $\tilde{Y}$. Soit $\tilde{P}_1'\in \pc(\tilde{M})$ l'unique sous-groupe parabolique inclus dans $\tilde{Q}$ tel que $\tilde{M}\cap \tilde{P}_1= \tilde{M}\cap \tilde{P}_1'$ ; on a   $\tilde{P}_1'=w_1^{-1}P_2w_1 $. On obtient 
$$ H_{P_2}(m_Q)=H_{P_2}(w_1 \tilde{m}_Q)=w_1\cdot H_{\tilde{P}_1'}( \tilde{m}_Q).$$

Écrivons $X=X_{\tilde{Q}}+X^{\tilde{Q}}$ selon la décomposition $\tilde{Q}=\mgo_{\tilde{Q}}\oplus\ngo_{\tilde{Q}}$. Soit $m\in \tilde{P}(F)\cap M_{\tilde{Q}}(F)$ et $n\in N_{\tilde{Q}}(F)$ tels $g=mn$. Par projection sur $\mgo_{\tilde{Q}}$ la relation $k^{-1}Yk= g^{-1}Xg$ implique  
$$\tilde{Y}=  \Ad(m^{-1})X_Q.$$
Par le lemme \ref{lem:non-pos} et avec ses notations (à l'exception notable de  $\al$ là-bas qui est devenu $-\al$ ici), on a 
$$R_{P_1}(g)-R_{P_2}(g)= \log|\det(m_1)^{-1}\det(m_3)|\al^\vee.$$
Pour $i=1,2$, on a $R_{P_i}(g)=R_{P_i}(m)= w_i\cdot H_{\tilde{P_i}}(m)$. Il s'ensuit qu'on est ramené à prouver l'énoncé suivant : soit $m\in  \tilde{P}(F)\cap M_{\tilde{Q}}(F)$ et  $\beta$ l'unique élément de $\Delta_{\tilde{P}_1'}^{\tilde{Q}}$. Il existe un unique $\rho\in \RR$ tel pour tout $\la \in a_{\tilde{M},\CC}^*$, la limite ci-dessous existe et est non identiquement nulle :
\begin{equation}
  \label{unelimite}
  \lim_{Z \to 0} |\be(Z)|^{\rho  \bg \la, \be^\vee\bd } \exp(-\bg \la, H_{\tilde{P_1'}}(n)\bd).
\end{equation}
Dans l'expression ci-dessus, on suppose qu'on a $Z\in \ago_{\tilde{M}}^{\reg}(F)$ et  $n\in M_{\tilde{Q}}(F)\cap N_{\tilde{P}_1}(F)$ est l'unique élément qui vérifie 
$$\Ad(n^{-1})Z=Z+  \Ad(m^{-1})X_Q.$$
De plus, cette limite est égale à $|\det(m_1)^{-1}\det(m_3)|^{\bg \la, \be^\vee  \bd }$  où $m_i$ est la projection de $m$ sur $\Aut_F(W_i)$. Il est clair que tout se passe dans le \og bloc\fg{} $\Aut_F(W_1\oplus W_2 \oplus W_3)$ de $M_{\tilde{Q}}$ c'est-à-dire que l'expression \eqref{unelimite} ne dépend de $Z$, $m$ etc. que via leur projection sur  $\Aut_F(W_1\oplus W_2 \oplus W_3)$. On a va donc supposer $M_{\tilde{Q}}=\Aut_F(W)$ où $W=W_1\oplus W_2 \oplus W_3$. On écrit matriciellement sous forme de matrices $3\times 3$ par blocs (éventuellement réduites à $2\times 2$ si $\dim(W_2)=0$). On a (ici on suppose implicitement que $\tilde{P}_1$ est le stabilisateur de $W_1$ ; l'autre situation où l'on échangerait le rôle de $\tilde{P}_1$ et $\tilde{P}_2$ conduit au même résultat) 
$$\tilde{P}_1= \{\left(
  \begin{array}{ccc}
    * & * & * \\ 0 & * & * \\ 0 & * & *
  \end{array}\right)\}  \text{   et  } \tilde{P}_1'= \{\left(
  \begin{array}{ccc}
    * & 0 & 0 \\ * & * & * \\ * & * & *
  \end{array}\right)\},
$$
$$X_Q=\left(
  \begin{array}{ccc}
    0 & 0 & I \\ 0 & 0 & 0 \\ 0 & 0 & 0
  \end{array}\right)
\text{   et    } m=\left(
  \begin{array}{ccc}
    m_1 & * & * \\ 0 & m_2 & * \\ 0 & 0 & m_3
  \end{array}\right)
$$
de sorte qu'on a 
$$\Ad(m^{-1})X_Q= \left(
  \begin{array}{ccc}
    0 & 0 & m_1^{-1}m_3 \\ 0 & 0 & 0 \\ 0 & 0 & 0
  \end{array}\right).$$

Un élément $Z\in \ago_{\tilde{M}}^{\reg}(F)$ s'écrit  
$$Z=  \left(
  \begin{array}{ccc}
    aI & 0 & 0 \\ 0 & bI & 0 \\ 0 & 0 & bI
  \end{array}\right),$$
où $a$ et $b$ sont deux éléments de $F$ tels que $a\not=b$. Le caractère $\beta$ est donné par $\beta(Z)=b-a$. On calcule facilement $n$ on trouve : 
$$n= \left(
  \begin{array}{ccc}
    I & 0 & -\beta(Z)^{-1} m_1^{-1}m_3 \\ 0 & I & 0 \\ 0 & 0 & I
  \end{array}\right).$$

Clairement, il suffit d'obtenir le résultat pour l'élément $\la$ particulier suivant : on va prendre $\la=\varpi$ où $\varpi\in X^*(\tilde{M})$ est le caractère défini par $\mathrm{diag}(m_1,m_2,m_3)\mapsto \det(m_2)\det(m_3)$. On retiendra qu'on a $\bg \varpi , \beta^\vee\bd =1.$ 

Rappelons qu'on note $r_1=\rang(W_1)=\rang(W_3)$ et $r_2=\rang(W_2)$. Prenons la puissance extérieure  $\bigwedge^{r_1+r_2}(W)$ qu'on munit de la base déduite  de la base de $W$ extraire de $e$ (cf.  §\ref{S:surZ}). Par extension des scalaires, on obtient un $F$-espace vectoriel  sur lequel le groupe $M_{\tilde{Q}}$ agit naturellement ;  on  note $\bigwedge^{r_1+r_2}$ cette représentation. Soit  $\|\cdot \|$ la norme sur $\bigwedge^{r_1+r_2}W\otimes_\ZZ F$ définie ainsi : si $F$ est réel (resp. complexe, resp. non-archimédienne), la norme est donnée par la racine carrée des valeurs absolues au carré (resp. la somme des valeurs absolues, resp. le maximum des valeurs absolues) des coefficients dans la base fixée. Par construction, la norme est invariante sous l'action du groupe $K$. Soit $e_\varpi$ l'unique élément de la base qui appartient à $\bigwedge^{r_1+r_2}(W_2\oplus W_3)$. On a 
$$\exp(\bg- \varpi, H_{\tilde{P}_1'}(n)\bd)= \| \bigwedge^{r_1+r_2}(n^{-1})  e_\varpi \|$$
La norme $ \| \bigwedge^{r_1+r_2}(n^{-1})  e_\varpi \|$ se calcule aisément en terme des déterminant extraits de taille $r_1+r_2$ de la matrice
$$n= \left(
  \begin{array}{cc}
     0 & \beta(Z)^{-1} m_1^{-1}m_3 \\ I & 0 \\ 0 & I
  \end{array}\right).
$$
En particulier, on a 
$$ \|\bigwedge^{r_1+r_2}(n^{-1})  e_\varpi \|= |\beta(Z)|^{-r_1} |\det(m_1^{-1}m_3)|+|\beta(Z)|^{-r_1+1} \phi(Z,m) $$
où  $\phi(Z,m)$ tend vers $0$ quand $Z$ tend vers $0$.
Quand $Z$ tend vers $0$, l'expression
\begin{equation}
  \label{eq:expr-norme}
   |\be(Z)|^{r_1} \exp(\bg- \varpi, H_{\tilde{P}_1'}(n)\bd)
 \end{equation}
a donc une limite qui vaut  $|\det(m_1^{-1}m_3)|$ ce qu'il fallait démontrer.
\end{preuve}

\end{paragr}

\section{Intégrale orbitale nilpotente : théorie locale}\label{sec:IOloc}

\begin{paragr}[Normalisation de mesures.] --- \label{S:norm}Dans cette section, $F$ est un corps local de caractéristique $0$. Les autres notations sont empruntées aux  sections précédentes, en particulier les section \ref{sec:standard} et  \ref{sec:fam-ortho}. 
  
Expliquons comment nous allons allons normaliser les mesures de Haar sur les différents groupes que nous rencontrerons. On munit le groupe additif $\RR$ de la mesure de Lebesgue usuelle notée $dx$, le groupeadditif  $\CC$ de la mesure $dz$ qui est le double de la mesure de Lebesgue. Si $F$ est non-archimédien, le groupe additif $F$ est muni de la mesure qui donne le volume $1$ à l'anneau des entiers $\oc$ de $F$. Les groupes multiplicatifs $\RR^\times$ et $\CC^\times$ sont respectivement munis des mesures $dx/|x|$ et $dz/|z|$ (rappelons que $|z|$ est la valeur absolue normalisée donc égale au carré de la valeur absolue usuelle). Le groupe multiplicatif $F^\times$ est muni de le mesure de Haar qui donne le volume $1$ à $\oc^\times$.

Soit $N$ un schéma en groupes sur $\ZZ$ qui admet une suite de composition $N_i$ dont les quotients sont isomorphes à $\mathbb{G}_{m,\ZZ}$. Chaque quotient $N_{i+1}(F)/N_i(F)$ est muni de la mesure de Haar choisie sur $F$. Chaque $N_i(F)$ est muni de la mesure de Haar telle que la mesure quotient sur  $N_{i+1}(F)/N_i(F)$ soit celle qu'on vient de choisir. En particulier, cela munit $N(F)$ d'une mesure de Haar qui ne dépend pas de la suite choisie. Cette normalisation s'applique en particulier aux radicaux unipotents des sous-groupes paraboliques semi-standard des sous-groupes de Levi semi-standard de $G$. Elle s'applique aussi au radical unipotent du centralisateur $G_X$ de $X$. Enfin elle s'applique aux algèbres de Lie des schémas en groupes sur $\ZZ$.

On a fixé au \ref{S:surZ} une $\ZZ$-base de $E$. On en déduit un isomorphisme sur $\ZZ$ entre le schéma en tores $\tc$ et $\mathbb{G}_{m,\ZZ}^n$. Sur $(F^\times)^n$, on met la mesure de Haar produit. Par transport, on en déduit une mesure de Haar sur $T_0(F)$.  Pour tout $M\in \lc(T_0)$, on munit $M(F)\cap K$ de la mesure de Haar qui donne le volume total égal à $1$.

Soit $dm$ la mesure de Haar sur $M(F)$ telle que pour toute fonction $f$ continue à support compact sur $M(F)$, on ait
  \begin{equation}
    \label{eq:MPK}
    \int_{M(F)} f(m) \, dm =  \int_{T_0(F)} \int_{N_B(F)\cap M(F)}\int_{K\cap M(F)} f(tnk) \,dt dn dk.
  \end{equation}
On en déduit la relation suivante, pour tout $P\in \pc(T_0)$ et  toute fonction $f$ continue à support compact sur $G(F)$
\begin{equation}
    \label{eq:GPK}
    \int_{G(F)} f(g) \, dg =  \int_{M_P(F)} \int_{N_P(F)}\int_{K} f(mnk) \,dm dn dk.
  \end{equation}

Fixons ensuite une mesure de Haar sur le centralisateur $G_X(F)$ de $X$ dans $G(F)$. On a déjà fixé une mesure sur le radical unipotent $N_X(F)$. On demande que la mesure quotient donne sur le facteur de Levi 
$M_X(F)$ (cf. §\ref{S:Levi-cent}) la mesure de Haar suivante : regardons comme $M_X$ comme schéma en groupes sur $\ZZ$. Comme tel, il s'identifie au produit $\prod_{1\leq j \leq r}   \Aut_\ZZ(V_j^j)$ et sur chaque facteur on utilise la mesure qu'on normalise comme ci-dessus. On a alors une décomposition d'Iwasawa sous la forme 
$$M_X(F)=T_X(F)\cdot (M_X(F)\cap N_B(F))\cdot( M_X(F)\cap K),$$
où $T_X=T\cap G_X,$ qui est compatible aux mesures suivantes :  le groupe unipotent $M_X(F)\cap N_B(F)$ est normalisé comme précédemment de même que le groupe compact  $M_X(F)\cap K$. Le sous-schéma en tore qui centralise $X$ est isomorphe sur $\ZZ$ à $\mathbb{G}_{m,\ZZ}^r$ d'où une mesure de Haar sur $T_X(F)$ par transport de la mesure produit.

On met sur l'espace homogène $G_X(F)\back G(F)$ la mesure quotient $dg^*$ qui est invariante à droite.
\end{paragr}

\begin{paragr}[Fonction zêta.] --- \label{S:zetaloc}Soit $\mathbf{1}$ la fonction sur $F$ suivante :
  \begin{itemize}
  \item  $\mathbf{1}(x)=\exp(-\pi x^2)$ si $F=\RR$ ;
  \item  $\mathbf{1}(z)=\exp(-\pi|z|)$ si $F=\CC$ ;
  \item $\mathbf{1}$ est la fonction caractéristique de $\oc$ si $F$ non-archimédien.
  \end{itemize}
Soit $\Sc(F)$ et $\Sc(\ggo(F))$ l'espace des fonctions complexes de Schwartz-Bruhat respectivement sur $F$ et $\ggo(F)$. Si $F$ est non-archimédien, c'est simplement l'espace des fonctions lisses à support compact.

Pour toute fonction $f\in \Sc(F)$, soit
$$Z(f,s)=\int_{F^\times} f(x) |x|^s \, dx.$$
L'intégrale est définie et holomorphe en $s$ pour $s\in \CC$ de partie réelle strictement positive. Elle admet un prolongement analytique à $\CC$. Lorsque $f=\mathbf{1}$, on note simplement $Z(s)=Z(\mathbf{1},s)$ cette fonction. Explicitement $Z(s)$ vaut (on note $\Gamma$ la fonction d'Euler)
 \begin{itemize}
  \item  $\pi^{-s/2}\Gamma(s/2)$ si $F=\RR$ ;
  \item  $(2\pi)^{1-s}\Gamma(s)$ si $F=\CC$ ;
  \item $\frac1{1-q^{-s}}$ où $q$ est le cardinal du corps résiduel de $F$.
  \end{itemize}

On identifie $\ggo(F)$ à $F^{n^2}$ via la base $e$ du § \ref{S:surZ} ; soit $\mathbf{1}_{i,j}$ la fonction sur $\ggo(F)$ obtenue par composition de $\mathbf{1}$ avec la fonction coordonnée d'indice $(i,j)$.  Soit $\mathbf{1}\in \Sc(\ggo(F))$ la fonction obtenue comme produit des fonctions $\mathbf{1}_{i,j}$  pour $1\leq i,j \leq n$. Ce procédé donne pour tous $F$-espaces $V_1$ et $V_2$ implicitement munis d'une base une fonction encore notée  $\mathbf{1}$ sur $\Hom_F(V_1,V_2)$, par exemple pour des sous-espaces de l'espace $E$ engendrés par des parties de sa base $e$.

Soit 
\begin{equation}
  \label{eq:Zd}
  Z_n(s)=\int_{G(F)} \mathbf{1}(g) |\det(g)|^s \, dg.
\end{equation}
Un calcul élémentaire utilisant la décomposition d'Iwasawa des mesures montre que $Z_n(s)$ est définie par une intégrale convergente pour $\Re(s)>n-1$ et qu'elle y vaut
$$Z_n(s)=Z(s)Z(s-1)\ldots Z(s-n+1).$$
En particulier, cette égalité donne le prolongement analytique à $\CC$. Par convention, on pose $Z_0(s)=1$.
\end{paragr}

\begin{paragr}[Comparaison d'intégrales.] ---  Le but de cette section est de prouver le lemme suivant.

\begin{lemme}\label{lem:Howe} Il existe $c_{X}>0$ tel que pour tout $P\in \rc(X)\cup \pc(M)$ et toute fonction $f\in \Sc(\ggo(F))$ on a 
\begin{equation}
  \label{eq:GXLie}
 \int_{G_X(F)\back G(F)} f(\Ad(g^{-1})X) \, dg^*=c_X\cdot \int_{\ngo_P(F)} \int_{K} f(k^{-1}Uk) \, dU dk.
\end{equation}
De plus, on a 
\begin{equation}
  \label{eq:formule-c}
  c_{X}=\prod_{j=1}^r  \prod_{i=1}^{j-1} Z_{d_j}(d_i+\ldots+d_{j})
\end{equation}
\end{lemme}

\begin{preuve}
Tout élément de $\pc(M)$ est conjugué sous $W$ à un élément de $\rc(X)$ (cf. lemme \ref{lem:r-std}). Il est clair alors qu'il suffit de prouver le lemme pour $P\in \rc(X)$, cas qu'on suppose désormais.

  Prouvons d'abord l'existence de la constante $c_X$. Tout d'abord les deux membres sont définis ; pour celui de droite c'est immédiat et pour celui de gauche, c'est un résultat général de Rao (cf. \cite{rao}, cf. aussi proposition 5 de \cite{Howe_germe}).  Soit $dp^*$ la  mesure   sur $G_X(F)\back P(F)$ relativement invariante à droite pour le caractère $\delta_P$ normalisée ainsi : c'est le quotient de la mesure sur $P(F)$ définie par 

$$\int_{P(F)}\phi(p)\, dp=\int_{M_P(F)}\int_{N_P(F)} \phi(mn) \, dmdn$$
par la mesure sur $G_X(F)$. Pour toute fonction lisse à support compact sur $G(F)$, on a avec nos choix de mesure
\begin{equation}
  \label{eq:GXG}
 \int_{G_X(F)\back G(F)} f(g) \, dg^*= \int_{G_X(F)\back P(F)} \int_{K} f(pk) \, dp^*dk.
\end{equation}
pour une certaine constante $c'>0$. Soit $V$ la $P(F)$-orbite de $X$. On sait que $V$ est un ouvert dense de $\ngo_P(F)$. On munit $V$ de la mesure induite par celle sur $\ngo_P(F)$. Il suffit de prouver le résultat suivant : il existe $c'>0$ telle que pour toute fonction $\phi$ complexe continue sur $V$, on a 

$$ \int_{G_X(F)\back P(F)} \phi(p^{-1}Xp) dp^* = c'\cdot \int_{V} \phi(Y)dY.
$$ 
Les deux membres définissent des mesures relativement invariantes sur lespace homogène $V$ pour le caractère $\delta_P$. L'existence de $c'$ résulte de ce qu'une telle mesure est unique à un scalaire près.

Pour terminer la preuve du lemme, il suffit d'évaluer chaque membre en la fonction $f=\mathbf{1}$, qui est $K$-invariante par conjugaison. On vérifie que l'intégrale du second membre vaut alors $1$ de sorte qu'on a 

$$c_X=\int_{G_X(F)\back G(F)} \mathbf{1}(g^{-1}Xg) \, dg^*.
$$
Le lemme résulte alors du lemme \ref{lem:calcul-IO} ci-dessous.
\end{preuve}

\begin{lemme}\label{lem:calcul-IO}
On a
$$  \int_{G_X(F)\back G(F)} \mathbf{1}(\Ad(g^{-1})X) \, dg^*= \prod_{j=1}^r  \prod_{i=1}^{j-1} Z_{d_j}(d_i+\ldots+d_{j}).$$ 
\end{lemme}

\begin{preuve}
  Soit $R$ le sous-groupe parabolique du §\ref{S:R} et $R=M_R N_R$ sa décomposition de Levi. Rappelons que le centralisateur $G_X$ possède une décomposition de Levi $M_X N_X$ avec $N_X\subset N_R$ et $M_X\subset M_R$ (cf. §§\ref{S:Levi-cent} et \ref{S:R}). Pour toute fonction $f$ continue à support compact sur $G_X(F)\back G(F)$, on a 
$$\int_{G_X(F)\back G(F)} f(g) \, dr^*= \int_{M_X(F)\back M_R(F)} \exp(-\bg 2\rho_R,H_R(m)\bd)\int_{N_X(F)\back N_R(F)} \int_K f(nmk)\,dk\,  dn^*\, dm^*,$$
où les quotients sont munis des mesures quotients. On en déduit que l'intégrale à calculer vaut 
$$  \int_{M_X(F)\back M_R(F)} \exp(-\bg 2\rho_R,H_R(m)\bd)\int_{N_X(F)\back N_R(F)} \mathbf{1}(\Ad(m^{-1}n^{-1})X )\,  dn^*\, dm^*.$$
On a défini au §\ref{S:o} un sous-$\ZZ$-module $\of$ de $\End_\ZZ(E)$. On met sur le $F$-espace vectoriel $\of(F)$ la mesure de Haar qui donne le volume $1$ à $\of(\oc)$. Son image par la translation par $X$ est l'orbite de $X$ sous $N_R(F)$.
En utilisant le lemme \ref{lem:int-o} ci-dessous, on est ramené à calculer
$$  \int_{M_X(F)\back M_R(F)} \exp(-\bg 2\rho_R,H_R(m)\bd)\int_{\of(F)} \mathbf{1}(\Ad(m^{-1})(X+Y))\,  dn^*\, dm^*.$$
Avec les notations des §§\ref{S:filtr} et \ref{S:o}, on a 
$$\Ad(m^{-1})X\in \bigoplus_{2\leq i \leq j\leq r} \Hom_F(V_j^i,V_j^{i-1})$$
et cet espace est en somme directe avec $\of(F)$. Il s'ensuit qu'on a 
$$ \mathbf{1}(\Ad(m^{-1})(X+Y))= \mathbf{1}(\Ad(m^{-1})X)\cdot \mathbf{1}(\Ad(m^{-1})Y).$$
On identifie $M$ à $\prod_{1\leq i \leq j \leq r} \Aut_F(V_j^i)$. On écrit $m=(m_j^i)$ selon cette décomposition. On utilise le lemme \ref{lem:calcul-jacob} ci-dessous pour obtenir que l'intégrale à calculer est égale à 
$$ \int_{M_X(F)\back M_R(F)}  \prod_{j=1}^r \prod_{i=1}^{j-1}  \mathbf{1}_{}( (m_j^{i})^{-1}m_j^{i+1}) |\det(m_j^{i})^{-1}\det(m_j^{i+1})|^{d_i+\ldots+d_{j}} \, dm  $$
qui est égale par un changement de variables évident à 
$$ \prod_{j=1}^r  \prod_{i=1}^{j-1} Z_{d_j}(d_i+\ldots+d_{j}).$$
\end{preuve}

Les deux lemmes suivants ont été utilisés dans la preuve du lemme \ref{lem:calcul-IO}.

\begin{lemme}\label{lem:int-o}
  Pour toute fonction $f$ continue à support compact sur $\of(F)$, on a 
$$\int_{N_X(F)\back N(F)} f(n^{-1}Xn-X) \, dn^* = \int_{\of(F)} f(Y)\,dY.$$
\end{lemme}

\begin{preuve}
C'est une conséquence immédiate de nos choix de mesure et du dévissage dans la preuve de la proposition \ref{prop:iso}.
\end{preuve}

\begin{lemme}
  \label{lem:calcul-jacob}Pour tout $m\in M_R(F)$, on a 
  \begin{equation}
    \label{eq:calcul-jac}
     \exp(-\bg 2\rho_R,H_R(m)\bd)\cdot \int_{\of(F)} \mathbf{1}(m^{-1}Ym)\, dY= \prod_{j=1}^r \prod_{i=1}^{j-1} |\det(m_j^{i})^{-1}\det(m_j^{i+1})|^{d_i+\ldots+d_{j}}
   \end{equation}
 \end{lemme}

 \begin{preuve}
   On a tout d'abord
$$ \exp(-\bg 2\rho_R,H_R(m)\bd)=\prod_{1\leq i \leq j \leq r}   |\det(m_j^i)|^{\rho_j^i}$$
où $\rho_j^i=\sum_{1\leq l \leq r} (\rho_j^{i,+}(l) - \rho_j^{i,-}(l)) d_l$ avec
$$\rho_j^{i,+}(l)= |\{ 1 \leq  k \leq l | p(V_l^k)<p(V_j^i) \}|$$
et
$$\rho_j^{i,-}(l)=|\{ 1 \leq  k \leq l | p(V_l^k)>p(V_j^i) \}|.$$
On a 
$$ \int_{\of(F)} \mathbf{1}(m^{-1}Ym)\, dY=  \prod_{1\leq i \leq j \leq r}   |\det(m_j^i)|^{\sigma_j^i}$$
où $\sigma_j^i=\sum_{1\leq l \leq r}( \sigma_j^{i,+}(l)- \sigma_j^{i,-}(l)) d_l$ avec
$$\sigma_j^{i,+}(l)= |\{ 1 \leq  k \leq l |  p(V_l^{k-1})>p(V_j^i) \}|$$
et
$$\sigma_j^{i,-}(l)=|\{ 1 \leq  k \leq l |  p(V_l^k)<p(V_j^{i-1})\}|.$$
On a ensuite

\begin{eqnarray*}
a_j^{i,+}(l)&:=&  \rho_j^{i,+}(l)- \sigma_j^{i,-}(l)\\
&=&  |\{ 1 \leq  k \leq l | p(V_j^{i-1})\leq p(V_l^k)<p(V_j^i) \}|\\
&=& \left\lbrace
  \begin{array}{l}
    1 \text{ si } i=1 \text{ et } j<l \text{ ou si } i\geq 2 \text{ et } l \geq i-1 ; \\
0 \text{ sinon.} 
  \end{array}\right.
\end{eqnarray*}

On a aussi
\begin{eqnarray*}
a_j^{i,-}(l)&:=& \rho_j^{i,-}(l)- \sigma_j^{i,+}(l)\\
&=&  |\{ 1 \leq  k \leq l | p(V_j^{k-1})\leq p(V_l^i)<p(V_j^k) \}|\\
&=& \left\lbrace
  \begin{array}{l}
    1 \text{ si } i< l \text{ ou si } l=i \text{ et } i<j ;\\
0 \text{ sinon.} 
  \end{array}\right.
\end{eqnarray*}
Soit $1\leq i \leq j$. On obtient alors,  avec $\delta_{i,j}$ le symbole de Kronecker,
\begin{eqnarray*}
  \sum_{k=1}^i \sum_{1\leq l \leq r} (a_j^{k,+}(l)-a_j^{k,-}(l))d_l &=& \sum_{l>j}d_l +  \sum_{k=2}^i \sum_{l\geq k-1}d_l -\sum_{k=1}^i \sum_{l\geq k} d_l + \delta_{i,j}d_j \\
&=& -\sum_{l=1}^{j} d_l +\sum_{k=2}^i d_{k-1} + \delta_{i,j}d_j\\
&=& -\sum_{l=1}^{j} d_l +\sum_{k=1}^{i-1} d_{k} + \delta_{i,j}d_j\\
&=& -\sum_{l=i}^{j} d_l + \delta_{i,j}d_j\\
&=& \left\lbrace
  \begin{array}{l}
    -(d_i+\ldots+d_{j}) \text{ si } 1\leq i <j \\
0 \text{ si } i=j.
  \end{array}\right.
\end{eqnarray*}
Le lemme s'en déduit.
\end{preuve}
\end{paragr}

\begin{paragr}[Une autre comparaison d'intégrales.] --- Dans ce paragraphe, on reprend les notations et les hypothèses en vigueur dans le §\ref{S:non-pos}. En particulier, $P_1$ et $P_2$ sont deux éléments adjacents de $\pc(M)$ et $Q$ le plus petit sous-groupe parabolique qui les contient tous deux. Les éléments  $\tilde{P}_1$, $\tilde{P}_2$ de $\rc(X)$ et $\tilde{Q}$ de $\lc\Sc(X)$ sont les images respectives de $P_1$, $P_2$ et $Q$ par l'application \eqref{eq:lappli2}. On a une décomposition, cf. \eqref{eq:decomp-n1-n2},
$$  \ngo_{\tilde{P}_1}\cap\ngo_{\tilde{P}_2}= \Hom_F(W_3,W_1)\oplus \ngo_{\tilde{Q}}$$
dont on note $U\mapsto U_{1,3}$ la première projection. Soit $I_{1,3}$ la variété sur $F$ des $F$-isomorphismes de $W_3\otimes_\ZZ F$ sur $W_1\otimes_\ZZ F$. C'est un torseur sous le groupe des $F$-automorphismes de $W_1\otimes_\ZZ F$. Le groupe des $F$-points de ce dernier est muni de la mesure de Haar normalisée selon les conventions de \ref{S:norm}. On en déduit une mesure invariante sur $I_{1,3}(F)$. Les modules $W_1$ et $W_3$ sont équipés de bases extraites de la base $e$ de $E$. On a alors un morphisme déterminant
$$\det : \Hom_F(W_3,W_1) \to \mathbb{G}_{a,F}$$
et $I_{1,3}$ est le lieu où ce morphisme ne s'annule pas. Par le choix de la base $e$, on a des inclusions naturelles
$$I_{1,3}\subset \Hom_F(W_3,W_1) \subset \End_F(E) =\ggo.$$

Rappelons enfin qu'on a défini des entiers $r_1>0$ et $r_2\geq 0$ en \eqref{eq:r1} et \eqref{eq:r2}.

\begin{lemme}\label{lem:alaRao}
Il existe une constante $c_{P_1,P_2}(X)>0$ telle que pour toute fonction $f\in \Sc(\ggo(F))$ on ait
\begin{equation}
  \label{eq:Rao}
 \int_{G_X(F)\back G(F)} f(g^{-1}Xg) \, dg^*=c_{P_1,P_2}(X)\cdot \int_{I_{1,3}(F)} \int_{\ngo_{\tilde{Q}}(F)} \int_{K} f(k^{-1}(x+ U)k)|\det(x)|^{r_1+r_2} \, dk\, dU \, dx.
\end{equation}
De plus, on a 
$$c_X=c_{P_1,P_2}(X)\cdot Z_{r_1}(r_1+r_2),$$
où $c_X$ est la constante du lemme \ref{lem:Howe}. 
\end{lemme}

\begin{preuve}
  La preuve est similaire à celle du lemme \ref{lem:Howe}. On note $\tilde{P}=\tilde{P}_1\cap \tilde{P}_2$. On traite d'abord l'existence de la constante de proportionnalité. Soit $V$ la $\tilde{P}(F)$-orbite de $X$ : c'est un ouvert dense de $\ngo_{\tilde{P}_1}(F)\cap\ngo_{\tilde{P}_2}(F)$ (cf. preuve du lemme \ref{lem:Y}). Comme dans la preuve du lemme \ref{lem:Howe}, on est prouve d'abord qu'il existe une constante $c>0$ telle que  pour toute fonction complexe continue $\phi$ sur $V$, on a 

$$ \int_{G_X(F)\back \tilde{P}(F)} \phi(p^{-1}Xp) dp^* = c \cdot \int_{V} \phi(U)|\det(U_{1,3})|^{r_2}  \, dU.
$$ 
où à gauche on utilise la mesure quotient. L'argument est le même : ce sont deux mesures relativement invariante pour le caractère  $\delta_{\tilde{P}}$. Par l'unicité d'une mesure invariante sur $I_{1,3}(F)$, il existe une constante $c'>0$ telle qu'on ait pour tout $\phi$ comme ci-dessus
$$\int_{V} \phi(U)|\det(U_{1,3})|^{r_2} \, dU= c'\cdot \int_{I_{1,3}(F)} \int_{\ngo_{\tilde{Q}}(F)}  \phi(x+U)|\det(x)|^{r_1+r_2} \, dx.
$$ 
Cela implique l'existence de $c_{P_1,P_2}(X)$. Pour déterminer cette constante, on évalue l'intégrale du second membre de \eqref{eq:Rao} en la fonction $f=\mathbf{1}$. On trouve aisément
$$\int_{I_{1,3}(F)} \int_{\ngo_{\tilde{Q}}(F)} \int_{K} \mathbf{1}(k^{-1}(x+ U)k)|\det(x)|^{r_1+r_2} \, dk\, dU \, dx=Z_{r_1}(r_1+r_2).$$
On utilise le lemme \ref{lem:calcul-IO} pour conclure.
\end{preuve}

\end{paragr}

\section{Intégrale orbitale pondérée nilpotente : théorie locale}\label{sec:IOPloc}

\begin{paragr}[Mesure sur $a_M^G$.] ---  \label{S:Haar-aM} Soit  $M\in \lc(T_0)$. Les espaces $a_M$ et $a_M^G$ sont munis d'un produit scalaire (cf. \ref{S:Weyl}). Il sont alors munis de la mesure de Haar qui donnent le covolume $1$ à un réseau engendré par une base orthonormale.
\end{paragr}

\begin{paragr}[Valeurs  associées à une $(G,M)$-famille.] --- \label{S:vM}
Soit $L\in \lc(M)$ et $Q\in \fc(L)$. Pour tout $P\in \pc^Q(L)$, soit la fraction rationnelle en $\la\in a_{L,\CC}^*$
\begin{equation}
  \label{eq:thetaP}
  \theta_P^Q(\lambda)=\vol(a_L^{M_Q}/\ZZ(\Delta_P^{Q,\vee}))\prod_{\al \in \Delta_P^Q} \frac1{\bg \lambda, \al^\vee\bd}.
\end{equation}
 Soit $v=(v_P(\la))_{P\in \pc(M)}$ une $(G,M)$-famille. On note $(v_P(\la))_{P\in \pc(L)}$ la  $(G,L)$-famille associée (cf. rappel en \ref{S:GMfam}). L'expression
$$\sum_{P\in \pc^Q(L)}v_P(\la) \theta_P^Q(\Lambda),$$
qui dépend de $\la\in a_{L,\CC}^*$, est en fait holomorphe en $\la=0$ (cf. \cite{trace_inv} lemme 6.2) et on note $v_L^Q$ sa valeur en $0$. On a  d'ailleurs pour tout $\Lambda\in a_{L,\CC}^*$  assez général et $k=\dim(a_L^Q)$,  (cf. \cite{trace_inv} éq. (6.5) p.37)
\begin{equation}
  \label{eq:derivee}
   v_{L}^{Q}=\frac1{k!} \sum_{P\in \pc^Q(L)} \frac{d^k}{dt^k} (v_P(t\Lambda))_{|t=0}  \theta_P^Q(\Lambda).
 \end{equation}

On a rappelé au §\ref{S:GMfam}comment on déduit une $(M_Q,M)$-famille $v^Q$ de toute  $(G,M)$-famille. Le procédé \eqref{eq:derivee} appliquée à cette famille $v^Q$ et au couple $(M_Q,L)$ (en lieu et place de $(Q,L)$) donne une valeur qui n'est autre que $v_L^Q$.
\end{paragr}

\begin{paragr}[Poids associé à un nilpotent standard.] --- On continue avec les notations des sections \ref{sec:standard}, \ref{sec:fam-ortho}, \ref{sec:GMfam} et \ref{sec:IOloc}. En particulier, $F$ est un corps local de caractéristique $0$, $X\in \ggo(F)$ est un élément nilpotent \og standard\fg{} et $M$ est facteur de Levi standard du sous-groupe parabolique $P_0$ qui stabilise le drapeau des noyaux itérés de $X$. On applique la construction précédente à la $(G,M)$-famille définie par
$$v_{P,X}(\la,g)=\exp(-\bg \la,R_P(g)\bd )$$
pour $P\in \pc(M)$. C'est bien une $(G,M)$-famille car la famille $(R_P(g))_{P\in \pc(M)}$ est orthogonale (cf. lemme \ref{lem:ortho}). On note $v_{L,X}^Q(g)$ la valeur associée en $0$ à cette famille et au couple $(L,Q)$. Lorsqu'on a $Q\in \pc(L)$, cette fonction vaut identiquement $1$.

\begin{lemme}
  \label{lem:inv-GX-loc} La fonction sur $G(F)$
$$g\mapsto v^Q_{L,X}(g)$$
est invariante à droite par $K$ et à gauche par le centralisateur $G_X(F)$ de $X$ dans $G(F)$. De plus, pour tout $w\in \Norm_W(M)$ et tout $g\in G(F)$, on  a
$$v^Q_{L,X}(g)= v^{Q^w}_{L^w,X}(g)$$
avec $L^w=w^{-1}Lw$ et $Q^w=w^{-1}Qw$.
\end{lemme}

\begin{preuve}
  L'invariance à droite par $K$ est évidente. Soit $g\in G(F)$ et  $h\in G_X(F)$. Il s'agit de montrer 
  \begin{equation}
    \label{eq:une-inv}
    v^Q_{L,X}(hg)=v^Q_{L,X}(g).
  \end{equation}
D'après le lemme \ref{lem:act-centralisateur},  on pour tout  $P\in \pc(M)$
$$R_P(hg)=R_{P_0}(h)+R_P(g).$$
On a donc pour $\la\in a_{L,\CC}^*$
$$\sum_{P\in \pc^{Q}(L)} \exp(\bg \lambda, -R_{P}(hg) \bd ) \theta_P^Q(\lambda)$$
$$=  \exp(-\bg \lambda, R_{P_0}(h) \bd ) \sum_{P\in \pc^{Q}(L)} \exp(\bg \lambda, -R_{P}(g) \bd ) \theta_P^Q(\lambda).$$
On déduit l'égalité \eqref{eq:une-inv} en prenant la limite pour $\lambda\to 0$. La denière égalité du lemme est une conséquence immédiate des définitions et de la formule \eqref{eq:normMsurL}.
\end{preuve}
\end{paragr}

\begin{paragr}[Intégrale orbitale pondérée nilpotente locale.] --- Les choix de mesures sont ceux de la section \ref{sec:IOloc}. Soit $L\in \lc(M)$ et  tout $Q\in \fc(L)$. Pour toute fonction $f \in \Sc(\ggo(F))$, on définit
\begin{equation}
    \label{eq:JMX-locale}
    J_{L,X}^Q(f)=\int_{G_X(F)\back G(F)} f(\Ad(g^{-1})X) v_{L,X}^Q(g)\, dg^*
  \end{equation}

  \begin{remarque}\label{eq:invariance}
    D'après le lemme \ref{lem:inv-GX-loc}, on a  $J_{L,X}^Q(f)=J_{L^w,X}^{Q^w}(f)$ pour tout $w\in \Norm_W(M)$.
  \end{remarque}

  \begin{proposition} \label{prop:cv-locale}Pour toute fonction $f \in \Sc(\ggo(F))$, l'intégrale \eqref{eq:JMX-locale} converge absolument.
\end{proposition}

\begin{preuve}
Si $Q\in \pc(L)$ on a $v_{L,X}^Q(g)=1$ pour tout $g\in G(F)$ et le résultat est un cas particulier du résultat général de Rao (cf. \cite{rao})

 Prenons un élément de $\pc(M)$ par exemple $P_0$. Soit $g\in G(F)$. La $(G,M)$-famille définie pour $P\in \pc(M)$ par
$$v_P(\la,g)=\exp(\bg \la,R_{P_0}-R_P(g)\bd )$$
donne les mêmes poids. Il résulte de \eqref{eq:derivee} qu'on a 
$$v_{L,X}^Q(g)= \frac1{k!} \sum_{P\in \pc^Q(M)} (\bg \Lambda,R_{P_0}-R_P(g)\bd)^k   \theta_P^Q(\Lambda).
$$
où $\Lambda\in a_{L,\CC}^*$  est assez général et $k=\dim(a_L^Q)$. On en déduit la majoration suivante : il existe $C>0$ tel pour tout $g\in G(F)$ on a la majoration suivante
\begin{equation}
  \label{eq:une-maj-adj}
  |v_{L,X}^{Q}(g)|\leq C\cdot \sum_{(P_1,P_2)\in \pc(M)^{\text{adj}}} \|R_{P_1}(g)-R_{P_2}(g)\|^k
\end{equation}
où $\pc(M)^{\text{adj}}\subset \pc(M)^2$ est l'ensemble formé de couples de paraboliques $(P_1,P_2)$ qui sont adjacents. La norme $\|\cdot\|$ est la norme euclidienne sur $a_M$. 

Soit $(P_1,P_2)\in \pc(M)^{\text{adj}}$. On est ramené à démontrer la convergence de l'intégrale orbitale où le poids $v_{L,X}^{Q}(g)$ est remplacé par $\|R_{P_1}(g)-R_{P_2}(g)\|^k$.

On reprend sans plus de commentaire les notations et les hypothèses du §\ref{S:non-pos}. D'après le lemme \ref{lem:Y}, pour tout  $k\in K$ et $U\in \ngo_{\tilde{P}_1}(F)\cap \ngo_{\tilde{P}_2}(F)$ tel que $kUk^{-1}=g^{-1}X g$, on a 
$$ \|R_{P_1}(g)-R_{P_2}(g)\|= |\log|\det(U_{1,3})| |  \|\al^\vee\|.$$
 D'après le lemme \ref{lem:alaRao}, on a, pour une certaine constante $c>0$ qu'il est inutile d'expliciter ici,

 $$\int_{G_X(F)\back G(F)} f(\Ad(g^{-1})X)  \|R_{P_1}(g)-R_{P_2}(g)\|^k \, dg^*$$
$$=c\cdot \int_{I_{1,3}(F)} \int_{\ngo_{\tilde{Q}}(F)} \int_{K} f(k^{-1}(x+U)k)|\log|\det(x)| |^k |\det(x)|^{r_1+r_2} \, dk\, dU \, dx.$$
On est ramené à prouver la convergence de cette dernière intégrale. Par une réduction immédiate, il suffit d'obtenir la convergence de 
$$\int_{GL(r_1,F)} f(x)|\log|\det(x)| |^k |\det(x)|^{r_1+r_2} \, dx$$ 
pour tout $f\in\Sc(\mathfrak{gl}(r_1,F))$. En utilisant la décomposition d'Iwasawa, on voit qu'il suffit de démontrer que l'intégrale suivante converge
$$\int_{T_{r_1}(F)} f(t)       |\log|\det(t)| |^k |\det(t)|^{r_1+r_2} |t_1|^{-r_1+1}|t_2|^{-r_2+2}\dots |t_{r_1-1}|        \,dt    $$
où $T_{r_1}$ est le tore diagonal de $GL(r_1)$ et $f\in\Sc(\mathfrak{t}_{r_1}(F))$. Quitte à majorer $f$ par un produit de fonctions des coordonnées $t_i$, on est ramené à étudier l'intégrale suivante en dimension $1$
$$\int_{F^\times} f(t) |\log|\det(t)| |^k |t|^{k'} \,dt$$
avec $k\geq 0$ et  $r_2+1\leq k'\leq r_2+r_1$. La convergence est alors évidente.
\end{preuve}

\end{paragr}

\begin{paragr}[Intégrale orbitale nilpotente d'Arthur.] --- \label{S:Arthur-loc}Soit $L\in \lc(T_0)$ et $\of\subset \lgo$ une $L$-orbite nilpotente. Soit $Q\in \pc(L)$. Au §\ref{S:ctr-Arthur} éq. \eqref{eq:wP'}, on a rappelé la construction d'Arthur d'une $(G,L)$-famille $(w_{Q'}(U,Y))_{Q'\in \pc(L)}$ définie pour tout $U\in \of(F)$ et tout $Y\in \ngo_Q(F)$ tels que $U+Y\in I_Q^G(\of)$. Soit $w_L(U,Y)$ la valeur associée à cette $(G,L)$-famille (cf. \ref{S:vM}).

Dans \cite{wei_or} Arthur a défini des intégrales orbitales unipotentes. La définition que met en avant Arthur (\emph{loc. cit.} formule (6.5) p. 254) ne semble pas très commode en pratique : dans sa formulation, elle invoque une limite d'intégrales orbitales pondérées \og usuelles\fg{}. En réalité, la construction d'Arthur repose fondamentalement sur l'utilisation des \og poids \fg{} $w_L(U,Y)$. Il est même possible d'extraire de  \emph{loc. cit.} une définition directe des intégrales orbitales pondérées nilpotentes : c'est le point de vue que nous adoptons.
 
Pour toute fonction $f\in \Sc(\ggo(F))$, on pose
\begin{equation}
  \label{eq:Jof}
  J_L^G(\of,f)=\int_{\ngo_P(F)} \int_{K} f(\Ad(k^{-1})U) \cdot w_L(U^L,U_Q)  \, dU dk,
\end{equation}
où
\begin{itemize}
\item $P\in \fc^Q(T_0)$ est un sous-groupe parabolique tel que $\of=I_{P\cap L}^L(0)$ (on a donc $I_Q^G(\of)=I_P^G(0)$);
\item on décompose $U\in \ngo_P$ en $U^L+U_Q$ selon $\ngo_P=(\lgo\cap \ngo_P) \oplus \ngo_Q$ ;
\item le poids  $w_L(U^L,U_Q)$ est bien défini sur l'ouvert dense $\ngo_P(F)\cap (\of(F)\oplus \ngo_Q(F))\cap I_P^G(0)$.
\end{itemize}
\emph{A priori}, cette intégrale dépend des choix auxiliaires $P$ et $Q$. Mais comme le suggère la notation, il n'en est rien. C'est l'objet de la proposition suivante, qui montre également la convergence.

\begin{theoreme}\label{thm:egalite-IOP} Soit $L\in \lc(T_0)$ et $\of\subset \lgo$ une $L$-orbite nilpotente. 

  \begin{enumerate}
  \item  Pour tout $f\in \Sc(\ggo(F))$, l'intégrale  $J_L^G(\of,f)$ définie par \eqref{eq:Jof} converge absolument  et ne dépend ni du choix de $P$ ni de celui de $Q\in \pc(L)$.
  \item Soit $X$ le représentant standard de l'orbite $I_L^G(\of)$ (au sens de la section \ref{sec:standard}) et soit $M\in \lc(T_0)$ attaché à $X$ (cf. \ref{S:surZ}). Il existe un élément $w\in W$ tel qu'on ait $wMw^{-1}\subset L$ et $\of=I_{wMw^{-1}}^L(0)$. On a 
  \begin{equation}
    \label{eq:egalite-IOP}
     J_{L^w,X}^G(f)=c_X\cdot J_{L}^G(\of,f)
   \end{equation}
   où $c_X$ est la constante du lemme \ref{lem:Howe}.
 \end{enumerate}
\end{theoreme}

\begin{remarque}
  Un élément $w\in W$ n'est pas unique dans l'assertion 2. Cependant différents choix de $w$ donnent des éléments $L^w$ qui sont conjugués sous le groupe $\Norm_{W}(M)$. En particulier, le membre de gauche dans \eqref{eq:egalite-IOP} ne dépend pas du choix de $w$ (cf. remarque \ref{eq:invariance}).
\end{remarque}

 \begin{preuve} L'assertion 1 est en fait une conséquence de l'assertion 2 et de la convergence donnée par la proposition \ref{prop:cv-locale}. Prouvons l'assertion 2. Soit $M'=M_P$. On a alors $\of=I_{M'}^L(0)$. Par transitivité de l'induction, on a $X\in I_L^G(I_{M'}^L(0))=I_{M'}^G(0)$. Or $X\in I_M^G(0)$. Il s'ensuit qu'il existe $w\in W$ tel que $M'=wMw^{-1}$. Posons $P_1=P^w$, $L_1=L^w$, $Q_1=Q^w$ et $\of_1=w^{-1}\of w$. Par définition, on a
   \begin{eqnarray*}
      J_L^G(\of,f)&=&\int_{\ngo_P(F)} \int_{K} f(\Ad(k^{-1})U) \cdot w_L(U^L,U_Q)  \, dU dk\\
&=&\int_{\ngo_{P_1}(F)} \int_{K} f(\Ad(k^{-1}w)U) \cdot w_L((\Ad(w)U)^L,(\Ad(w)U)_Q)  \, dU dk\\
&=&\int_{\ngo_{P_1}(F)} \int_{K} f(\Ad(k^{-1}w)U) \cdot w_{L_1}(U^{L_1},U_{Q_1}) \, dU dk
\end{eqnarray*}
La dernière égalité résulte de 
$$w_L((\Ad(w)U)^L,(\Ad(w)U)_Q)=w_{L_1}(U^{L_1},U_{Q_1})$$
pour tout $U\in \ngo_{P_1}(F)$.

 Le poids $v_{L,X}^{Q}(g)$ pour $g\in G(F)$ se déduit de la $(G,M)$-famille 
$$(\exp(\bg \la,R_{P_1}-R_{P'}(g)\bd )_{P'\in \pc(M)}.$$
 Soit $U\in \ngo_{P_1}(F)\cap (\of_1(F)\oplus \ngo_{Q_1}(F))\cap I_{P_1}^G(0)$ et $k\in K$ tel que $k^{-1}Uk=g^{-1}Xg$. D'après la proposition \ref{lem:comparaisonGM}, cette $(G,M)$-famille n'est autre que la famille $(w_{P'}(\la,U))_{P'\in \pc(M)}$ définie relativement au sous-groupe parabolique point-base $P_1$. En particulier, les $(G,L_1)$-familles qui s'en déduisent ainsi que leur valeur obtenue par le procédé du §\ref{S:vM} sont égales. Or d'après le lemme \ref{lem:descente-GM}, la $(G,L_1)$-famille déduite de $(w_{P'}(\la,U))_{P'\in \pc(M)}$ n'est autre que $(w_{Q'}(U^{L_1},U_{Q_1}))_{Q'\in \pc(L_1)}$. Il en résulte qu'on a 
$$v_{L_1,X}^{Q_1}(g)=w_{L_1}(U^{L_1},U_{Q_1}).$$
L'assertion 2 (ainsi que la convergence) est alors une conséquence du lemme \ref{lem:Howe}.
    \end{preuve}
\end{paragr}

\begin{paragr}[Formule de descente.] --- \label{S:descente} Au paragraphe \ref{S:Arthur-loc} ci-dessus, on a rappelé la définition d'Arthur des intégrales orbitales pondérées nilpotentes. Le poids et la définition des intégrales se généralisent aux sous-groupes de Levi de $G$. Pour tous sous-groupes de Levi $L_1\subset L_2$ et toute fonction $f\in \Sc(\lgo_2(F))$, pour toute $L_1$-orbite nilpotente $\of\subset\lgo_1$, on dispose de l'intégrale $J_{L_1}^{L_2}(\of,f)$ construite à l'aide du poids $w_{L_1}^{L_2}$ relatif à $L_2$. On a la formule de descente suivante

\begin{lemme} \label{lem:descente} Soit $L\in \lc(M)$. Pour tout $Q\in \pc(L)$ et tout $f\in \Sc(\ggo(F))$, on a 
$$J^{Q}_{M,X}(f)=c_X\cdot  J_M^L((0),f_Q)$$
où 
\begin{itemize}
\item l'on pose
$$f_Q(U)=\int_{\ngo_Q(F)} \int_K f(\Ad(k^{-1})(U+V)) \, dV\ ;$$
\item $(0)$ est l'orbite nulle dans $\mgo$ ;
\item la constante $c_X$ est celle définie dans le lemme \ref{lem:Howe}.
\end{itemize}
\end{lemme}

\begin{preuve}
Soit $P_0\in \pc^L(M)$ et $\ngo_0=\ngo_{P_0}$.   Par le lemme \ref{lem:Howe}, on  a 
$$J^{Q}_{M,X}(f)=c_X \cdot \int_{(\ngo_0\oplus\ngo_Q)(F)} f(\Ad(k^{-1})U) v^Q_M(U) \, dU $$
où l'on a défini pour $U$ dans l'ouvert $(\ngo_0\oplus\ngo_Q)(F)\cap I_{P_0N_Q}^G(0)$ un poids $v^Q_M(U)$ de la façon suivante : c'est la valeur associée à la $(L,M)$-famille définie pour $P\in \pc^L(M)$ par
$$
w^Q_P(\la,U)= \exp(\bg \la, R_{P_0N_Q}(g)-R_{PN_Q}(g)\bd)
$$
où $g\in G(F)$ doit satisfaire l'égalité $\Ad(g^{-1})X=U$. Bien sûr, cette $(L,M)$-famille est indépendante du choix de $g$. Soit $w_0\in W$ et $Y_0\in \ngo_0(F)$ tels que $w_0^{-1}P_0N_Qw_0\in \rc(X)$ et $\Ad(w_0)X\in Y_0\oplus \ngo_Q(F)$. Il existe $p\in P_0(F)$ et $n\in N_Q(F)$ tels que 

$$ \Ad (w_0^{-1} p n)^{-1}X= U.$$
En particulier, en écrivant $U=U_0+U_Q$ selon $\ngo_0\oplus\ngo_Q$, on  a
$$\Ad(p^{-1}) Y_0 = U_0.$$
 Soit  $\rc^L(Y_0)$  l'ensemble défini en \eqref{eq:sgp-Ric} relativement au groupe $L$ et à l'élément $Y_0$. Soit $P\in \pc^L(M)$.  Il existe  $w_P\in W^L$  tel que $w_P^{-1}P w_P\in \rc^L(Y_0)$ (la classe de $w_P$ dans $W^M\back W^L$ est uniquement définie ; le point-clef est l'inclusion $\rc^L(Y_0)\subset \lc^L(T_0)$, les arguments sont ceux de la démonstration du lemme \ref{lem:r-std}). Pour un tel $w_p$, on a $\Ad(w_Pw_0)X\in \ngo_P\oplus\ngo_Q$ et donc $(w_Pw_0)^{-1} PN_Q w_Pw_0 \in \rc(X)$. En se reportant à la définition \eqref{eq:HQ-RQ}, on obtient
\begin{eqnarray*}
   R_{P_0N_Q}(g)-R_{PN_Q}(g)&=&  R_{P_0N_Q}(w_0^{-1} p n )-R_{PN_Q}(w_0^{-1} p n )\\
&=& H_{P_0N_Q}(pn)-H_{PN_Q}(w_P p n )\\
&=& H_{P_0N_Q}(p)-H_{PN_Q}(w_P p ).
\end{eqnarray*}
Maintenant il résulte de la proposition \ref{lem:comparaisonGM} qu'on a 
$$v_M^Q(U)= w_{M}^L(U_0).$$
Le lemme en découle immédiatement.
\end{preuve}

\end{paragr}

\section{Un calcul d'intégrale orbitale pondérée nilpotente locale}\label{sec:uncalculIOPloc}

\begin{paragr} Dans cette section, les hypothèses et les notations sont celles de \ref{sec:IOloc} et \ref{sec:IOPloc}. En particulier $F$ est un corps local de caractéristique $0$ et $X\in \ggo(F)$ un élément nilpotent \og standard\fg{}. On suppose dans toute cette section qu'on a $\inv(X)=1$ (cf. \eqref{eq:invX}). On dit encore suivant \cite{scfhn} que $X$ est régulier par blocs. Soit $M$ le facteur de Levi standard de $P_0$, le stabilisateur du drapeau des noyaux itérés de $X$.

 Le but de cette section est le calcul des intégrales nilpotentes d'Arthur $J_L^G(I_M^L(0),\mathbf{1})$ pour $L\in \lc(M)$ où $\mathbf{1}\in \Sc(\ggo(F))$ est la fonction définie au § \ref{S:zetaloc}. En particulier si $F$ est non-archimédien, cette fonction est la fonction caractéristique de $\ggo(\oc)$. La formule obtenue peut être vue comme un analogue sur les corps locaux de caractéristique $0$ du résultat principal de \cite{scfhn}.
\end{paragr}

\begin{paragr}[Formulation du résultat.] --- Sous notre hypothèse $\inv(X)=1$, la seule multiplicité non nulle d'un bloc de Jordan de $X$ est $d_r$ pour $r$ l'indice de nilpotence de $X$. On  a  donc $rd_r=n$ et $M\simeq GL(d_r)^r$ et le groupe $W(M)=\Norm_W(M)/W^M$ s'identifie au groupe symétrique en $d$ variables ; il agit librement et transitivement sur $\pc(M)$ par $(w,P)\mapsto P^w=w^{-1}Pw$. On introduit la $(G,M)$-famille $(\mathcal{J}_{P,X}(\la))_{P\in \pc(M)}$ de fonctions sur $a_{M,\CC}^{G,*}$ holomorphes au voisinage de $0$ définie ainsi : pour tout $P\in \pc(M)$ 
$$\mathcal{J}_{P,X}(\la)= \prod_{\varpi^\vee \in \hat{\Delta}^\vee_{P_0}} \frac{Z_{d_r}(d_r+\bg w\cdot \la ,\varpi^\vee\bd)}{Z_{d_r}(d_r)}.$$
où 

\begin{itemize}
\item $w$ est l'unique élément de $W(M)$ tel que $P_0^w=P$ ;
\item $w\cdot \la$ désigne l'action naturelle de $W(M)$ sur $a_{M,\CC}^{G,*}$ ;
\item $\hat{\Delta}^\vee_{P_0}$ est l'ensemble des copoids $P_0$-dominants fondamentaux ;
\item la fonction $Z_{d_r}$ est celle définie en \eqref{eq:Zd}.
\end{itemize}

  \begin{theoreme}\label{thm:calcul-IOPlocale}
  Pour tout $L\in \lc(M)$, on a 
$$J_L^G(I_M^L(0),\mathbf{1})= \mathcal{J}_{L}$$
où $\mathcal{J}_{L}$ est la valeur associée à la $(G,L)$-famille déduite de $(\mathcal{J}_{P,X}(\la))_{P\in \pc(M)}$.
  \end{theoreme}

\begin{preuve} Soit $Q_1\in \pc(L)$ et $P_1\in \pc^Q(M)$. D'après \eqref{eq:Jof}, on a pour notre fonction particulière
 $$J_L^G(\of,f)=\int_{\ngo_{1}(F)} \mathbf{1}(U) \cdot w_L(U',U'')  \, dU 
$$
où l'on pose $\ngo_1=\ngo_{P_1}$ et $U=U'+U''$ selon $\ngo_1=(\lgo\cap \ngo_1 )\oplus \ngo_{Q_1}$. D'après la proposition \ref{lem:comparaisonGM} et le lemme \ref{lem:descente-GM}, le poids  $w_L(U',U'')$ est donné par la valeur associée à la $(G,L)$-famille déduite de la $(G,M)$-famille 
$$(\exp(\bg \la, R_{P_1}(g)-R_P(g)\bd))_{P\in \pc(M)},$$
où $g$ est un élément de $G(F)$ tel que $\Ad(g^{-1})X$ appartient à la $K$-orbite de $U$. Il est facile de voir que pour tout $P\in \pc(M)$, on a $R_P(G_X(F))\subset a_G$. En particulier, comme il est loisible, on se restreint à $\la\in a_{M,\CC}^{G,*}$, on voit que $(\exp(\bg \la, -R_P(g)\bd))$ ne dépend pas du choix de $g$. Il est facile de voir que le poids $w_L(U',U'')$ est aussi égal à la valeur associée à la $(G,L)$-famille déduite de la $(G,M)$-famille 
$$(\exp(\bg \la, -R_P(g)\bd))_{P\in \pc(M)}.$$
En particulier, pour un point $\Lambda\in a_{L,\CC}^{G,*}$ assez général (c'est-à-dire hors des hyperplans singuliers des $\theta_Q$) on a (cf. \eqref{eq:derivee})
$$w_L(U',U'')=\frac{1}{k!}\sum_{Q\in \pc(L)}  (\bg \Lambda ,  -R_Q(g)\bd)^k \theta_Q(\Lambda).$$
où $k=\dim(a_L^G)$. On a donc formellement à l'aide du lemme \ref{lem:Howe} (la constante $c_X$ est celle qui apparaît dans ce lemme) 
\begin{eqnarray}\label{eq:laserie}
  J_L^G(\of,f)&=& \frac{1}{k!\cdot c_X}\sum_{Q\in \pc(L)} \theta_Q(\Lambda)  \int_{G_X(F)\back G(F)}  \mathbf{1}(\Ad(g^{-1})X)  (\bg \Lambda ,  -R_Q(g)\bd)^k \, dg\\
\nonumber &=&  \frac{1}{k!\cdot c_X}\sum_{Q\in \pc(L)} \theta_Q(\Lambda)  \int_{G_X(F)\back G(F)}  \mathbf{1}(\Ad(g^{-1})X)  (\bg \Lambda ,  -R_Q(g)\bd)^k \, dg\\
\nonumber &=& \frac1{k!\cdot c_X} \sum_{Q\in \pc^Q(L)} \theta_Q(\Lambda)  \frac{d^k}{dt^k} (\mathcal{I}_{Q,X}(t\Lambda))_{|t=0},
\end{eqnarray}
où l'on pose pour tout $Q\in \fc(M)$ et tout $\la\in  a_{M,\CC}^{G,*}$
$$\mathcal{I}_{Q,X}(\la)= \int_{G_X(F)\back G(F)}  \mathbf{1}(\Ad(g^{-1})X) \exp(- \bg \la, R_Q(g)\bd \,dg.  
$$
On a, du moins formellement,  
\begin{equation}
  \label{eq:IP=IQ}
  \mathcal{I}_{P,X}(\la)=\mathcal{I}_{Q,X}(\la)
\end{equation}
si $\la \in a_{M_Q,\CC}^{G,*}$ et si $P\subset Q$ (cf.  §\ref{S:orth-famille}). Notons également que \eqref{eq:normMsurL} implique qu'on a pour tout $\la\in a_{M,\CC}^{G,*}$ et tout $w\in W(M)$
\begin{equation}
  \label{eq:IsousW}
  \mathcal{I}_{P_0^w,X}(\la)=\mathcal{I}_{P_0,X}(w\cdot \la).
\end{equation}
Un calcul réminiscent de celui effectué lors de la démonstration du lemme \ref{lem:calcul-IO} montre qu'on a
\begin{eqnarray}\label{eq:sousformezeta}
  \mathcal{I}_{P_0,X}(\la)&=& \prod_{\varpi^\vee \in \hat{\Delta}^\vee_{P_0}} \int_{GL(d_r,F)} \mathbf{1}_{d_r}(x)  |\det(x)|^{d_r+\bg  \la ,\varpi^\vee\bd}dx \\
\nonumber &=&\prod_{\varpi^\vee \in \hat{\Delta}^\vee_{P_0}} Z_{d_r}(d_r+\bg  \la ,\varpi^\vee\bd).
\end{eqnarray}
où  $\mathbf{1}_{d_r}$ est la fonction  $\mathbf{1}$ relative à $GL(d_r)$. On voit donc que, sur l'ouvert formé des $\la\in a_{M,\CC}^{G,*}$ tels que $\bg  \la ,\varpi^\vee\bd>-1$ pour tout $\varpi^\vee \in \hat{\Delta}^\vee_{P_0}$, l'intégrale définissant $\mathcal{I}_{P_0,X}(\la)$. Vu \eqref{eq:IsousW} et \eqref{eq:IP=IQ}, on en conlut que toutes les intégrales $I_Q(\la)$ convergent pour $\la$ dans un voisinage de $0$. De même, la série d'égalité \eqref{eq:laserie} est valide pour $\Lambda$ dans un voisinage de $0$. Maintenant, en combinant \eqref{eq:IsousW} et \eqref{eq:sousformezeta} ainsi que les formules explicites $|\hat{\Delta}^\vee_{P_0}|=r-1$,
$$c_X= Z_{d_r}(d_r)^{r-1}$$
on obtient  
$$c_X^{-1}\mathcal{I}_{P,X}(\la)=\mathcal{J}_{P,X}(\la)$$
pour tout $P\in \pc(M)$ et $\la$ au voisinage de $0$. Le théorème est alors une conséquence immédiate de \eqref{eq:laserie}.
\end{preuve}

\end{paragr}

\section{Intégrale orbitale pondérée nilpotente : théorie globale}\label{sec:IOPglob}

\begin{paragr}[Notations et choix des mesures de Haar.] ---\label{S:notations-Haar}
  Soit $F$ un corps global de caractéristique $0$. On reprend les notations de \ref{sec:standard} relativement à ce corps. Soit $V$ l'ensemble des places de $F$. Pour toute place $v\in V$, soit $F_v$ le complété du corps $F$ en $v$. Pour toute place non-archimédienne, soit $\oc_v\subset F_v$ l'anneau des entiers. Les notations des sections \ref{sec:fam-ortho} à \ref{sec:IOPloc} s'appliquent relativement au corps local $F$. En général, on affuble un objet relatif à $F_v$ d'un indice $v$. Par exemple, $K_v$ est le sous-groupe compact maximal de $G(F_v)$ défini au §\ref{S:HP}. Soit $K=\prod_{v\in V} K_v$. 

Soit $\AAA$ l'anneau des adèles de $F$. On note simplement $|\cdot|$ la valeur absolue adélique, qui est le produit des valeurs absolues normalisées locales $|\cdot|_v$  sur les corps $F_v$. En particulier, la valeur absolue adélique est triviale sur $F$.

Chaque groupe $G(F_v)$ est muni de la mesure de Haar décrite au §\ref{S:norm}. De même, chaque groupe $K_v$ est muni d'une mesure de Haar qui est le produit des mesures de Haar choisie sur les facteurs $K_v$. De façon générale, les groupes de points adéliques rencontrés sont munis de la mesure de Haar produit des mesures de Haar sur leurs facteurs choisies conformément au §\ref{S:norm}. L'espace homogène $G_X(\AAA)\back G(\AAA)$ est muni de la mesure quotient notée $dg^*$. 
\end{paragr}

\begin{paragr}[Fonction zêta globale.] --- Soit $\Sc(\ggo(\AAA))$ l'espace des fonctions complexes de Schwartz-Bruhat  $\ggo(\AAA)$. Pour tout $f\in \Sc(\ggo(\AAA))$, soit
\begin{equation}
  \label{eq:Zd-globale}
  Z_n(f,s)=\int_{G(\AAA)} \mathbf{1}(g) |\det(g)|^s \, dg.
\end{equation}
On note simplement $Z(f,s)=Z_1(f,s)$ dans le cas particulier  où le rang de $G$ vaut $1$. Il est bien connu depuis la thèse de Tate que l'intégrale qui définit $Z(f,s)$ est absolument convergente pour $s\in \CC$ de partie réelle $\Re(s)>1$ et se prolonge analytiquement à $\CC$. Soit $\mathbf{1}_v$ la fonction définie sur $\ggo(F_v)$ définie au \ref{S:zetaloc}. 

Soit $\mathbf{1}=\otimes_{v\in V} \mathbf{1}_v \in \Sc(\ggo(\AAA))$. On pose  $Z_n(s)=Z_n(\mathbf{1},f)$ et $Z(s)=Z_1(s)$. Comme dans le cas local, on montre que pour $s\in \CC$ de partie réelle $\Re(s)>n$, l'intégrale qui définit $Z_n(s)$  est absolument convergente et qu'on a 
  $$Z_n(s)=Z(s)Z(s-1)\ldots Z(s-n+1)$$
ce qui donne le prolongement analytique à $\CC$. Si l'on note $Z_{n,v}$ les fonctions zêta locales, on apour $\Re(s)>n$ l'égalité avec le produit convergent
$$Z_n(s)=\prod_{v\in V} Z_{n,v}(s).$$
\end{paragr}

\begin{paragr}[Intégrales orbitales nilpotentes globales.] --- Le théorème suivant montre que seules certaines intégrales  orbitales nilpotentes globales vont converger. Cela explique pourquoi nous allons dans la suite imposer la condition d'être simple à l'endomorphisme nilpotent $X$ (cf. définition \ref{def:simple}. 
 
  \begin{theoreme}\label{thm:cv}
  L'intégrale orbitale globale
  \begin{equation}
    \label{eq:IO-globale}
    \int_{G_X(\AAA)\back G(\AAA)} f(g^{-1}Xg)\, dg^*
  \end{equation}
  est convergente pour toute fonction $f\in \Sc(\ggo(\AAA))$ si et seulement si $X$ est simple au sens de la définition \ref{def:simple}.

De plus, si $X$ est simple, on a 
\begin{equation}
  \label{eq:GXLie-glob}
  \int_{G_X(\AAA)\back G(\AAA)} f(g^{-1}Xg)\, dg^*=c_X\cdot \int_{\ngo_P(F)} \int_{K} f(k^{-1}Uk) \, dU dk.
\end{equation}
où l'on pose 
\begin{equation}
  \label{eq:formule-c-globale}
  c_{X}=\prod_{j=1}^r  \prod_{i=1}^{j-1} Z_{d_j}(d_i+\ldots+d_{j}).
\end{equation}
\end{theoreme}

\begin{preuve}
Traitons la condition nécessaire. Il suffit de prendre $f=\mathbf{1}$. Par nos choix des mesures, l'intégrale orbitale est alors un produit d'intégrales orbitales locales qui se calculent à l'aide du lemme \ref{lem:calcul-IO}. On obtient donc que l'intégrale \eqref{eq:IO-globale} converge pour $f=\mathbf{1}$ si et seulement si le produit suivant
$$\prod_{v\in V }\prod_{j=1}^r  \prod_{i=1}^{j-1} Z_{d_j,v}(d_i+\ldots+d_{j})$$
est convergent. Si c'est le cas, ce produit vaut  $c_{X}$. Or ce produit converge si et seulement si chacun des produits ci-dessous converge
 $$\prod_{v\in V } Z_{d_j,v}(d_i+\ldots+d_{j}).$$
pour $1\leq i <j \leq r$. Ce dernier est convergent si et seulement $d_j=0$  ou 
$$d_i+\ldots+d_{j}>d_j.$$
Si $r=1$, il n'y a aucune condition à vérifier. Si $r>1$, on a doit avoir $d_{r-1}>0$ pour $j=r$ et $i=r-1$. De fil en aiguille, on obtient comme condition nécessaire et suffisante $d_j>0$ pour tout $1\leq j \leq r$ : c'est bien la condition de simplicité pour $X$.

Supposons $X$ simple. Alors on vient de voir que la constante $c_X$ est finie et est égale à l'intégrale  \eqref{eq:IO-globale} pour $f=\mathbf{1}$. Pour conclure, il suffit d'invoquer le lemme \ref{lem:Howe} et d'observer que l'intégrale dans le membre de droite de \eqref{eq:GXLie-glob} est convergente pour toute $f\in \Sc(\ggo(\AAA))$  et vaut $1$ pour $f=\mathbf{1}$.

\end{preuve}

\end{paragr}

\begin{paragr}[Poids global.] --- Les notations $T_0,B,M$ sont celles du §\ref{S:surZ}. Soit $T\in a_{T_0}$. Pour tout sous-groupe de Borel $B'\in \fc(T_0)$, il existe un unique élément $w\in W$ tel que $wB w^{-1}=B'$. On pose alors 
$$T_{B'}=w\cdot T.$$
La famille $(T_{B'})_{B'\in \pc(T_0)}$ est orthogonale. On en déduit la famille orthogonale $(T_P)_{P\in \pc(M)}$ (cf. §\ref{S:fam-ortho}). Soit $g=(g_v)_{v\in V }\in G(\AAA)$. Pour tout $P\in \fc(T_0)$, on définit une application 
$$H_P: G(\AAA)\to a_P$$
par 
$$H_P(g)=\sum_{v\in V} H_{P,v}(g_v),$$
où les fonctions $H_{P,v}$ sont celles relatives au corps local $F_v$ définies au §\ref{S:HP} et où la somme est en fait à support finie. Pour tout $P\in \pc(M)$, on introduit la fonction $R_P$ définie pour $g\in G(\AAA)$ par
$$R_P(g)=H_P(w_Pg)$$
où $w_P\in W$ vérifie $w_P^{-1}Pw_P\in \rc(X)$. C'est le pendant évident des fonctions $R_{P,v}$ locales du §\ref{S:RP}. On a bien sûr la relation suivante :
 
\begin{equation}
  \label{eq:RP-RPv}
  R_P(g)=\sum_{v\in V} R_{P,v}(g).
\end{equation}

Il résulte de cette égalité et du lemme \ref{lem:ortho} que la famille définie pour $\la\in a_{M,\CC}^*$ et $P\in \pc(M)$ par
$$v_{P,X}(\la,g,T)=\exp(\bg \la,T_P-R_P(g)\bd )$$
 est une $(G,M)$-famille. Pour tout $L\in \lc(M)$ et $Q\in \fc(L)$, on note $v_{L,X}^Q(g,T)$ le \og poids\fg{} qui s'en déduit (cf. \ref{S:vM}). Lorsque $T=0$, on note simplement $v_{P,X}(\la,g)=v_{P,X}(\la,g,0)$ et 
$$v_{L,X}^Q(g)=v_{L,X}^Q(g,0).
$$

\begin{lemme}
  \label{lem:inv-GX} La fonction sur $G(\AAA)$
$$g\mapsto v^Q_{L,X}(g,T)$$
est invariante à droite par $K$ et à gauche par le centralisateur $G_X(\AAA)$ de $X$ dans $G(\AAA)$.
\end{lemme}

\begin{preuve}
  Elle est similaire à la preuve du lemme \ref{lem:inv-GX-loc} compte tenu de la relation \eqref{eq:RP-RPv}.

\end{preuve}
  
\end{paragr}

\begin{paragr}[Intégrale orbitale pondérée globale pour $X$ nilpotent simple.] --- Soit $L\in \lc(M)$, $Q\in \fc(L)$ et $T\in a_{T_0}$. On suppose que $X$ est simple au sens de la définition \ref{def:simple}. Pour tout $f\in \Sc(\ggo(\AAA))$, on introduit l'intégrale orbitale pondérée globale
 \begin{equation}
     \label{eq:JMXG}
     J_{L,X}^{Q,T}(f):=\int_{G_X(\AAA)\back G(\AAA)} f(\Ad(g^{-1})X)   v^Q_{L,X}(g,T) \, dg^*.
   \end{equation}
Lorsque $T=0$, on pose
$$ J_{L,X}^{Q}(f)=J_{L,X}^{Q,T=0}(f)
$$

  \begin{theoreme}\label{thm:cv-pondere}
Supposons $X$ simple.  Soit $f\in \Sc(\ggo(\AAA))$. L'intégrale \eqref{eq:JMXG} est absolument convergente. En outre, l'application de source $a_{T_0}$
$$T\mapsto J_{L,X}^{Q,T}(f)$$
est polynomiale en $T$.
  \end{theoreme}

  \begin{preuve}
La $(G,M)$-famille $(\exp(\bg \Lambda, T_{P}-R_{P}(g) \bd ))_{P\in \pc(M)}$ est le produit des deux $(G,M)$-familles $(\exp(\bg \Lambda, T_{P}\bd ))_{P\in \pc(M)}$ et $(\exp(\bg \Lambda, -R_{P}(g) \bd ))_{P\in \pc(M)}$. D'après Arthur (cf. \cite{inv_loc}, corollaire 7.4), on a 
  $$v_{L,X}^Q(g,T)=\sum_{(L_1,L_2)\in \lc(L)^2} d_{L}^{M_Q}(L_1,L_2)   v_{L,X}^{Q_1}(1,T) v_{L,X}^{Q_2}(g,0).$$

Les coefficients  $ d_{L}^{M_Q}(L_1,L_2)$ ont peu d'importance ici et on ne rappellera pas  leur définition. L'expression ci-dessus dépend d'un certain choix qui fixe une section de l'application 
$$(Q_1,Q_2)\in \fc(L)^2 \mapsto (M_{Q_1},M_{Q_2})\in \lc(L)^2$$
définie sur la partie formée des couples $(L_1,L_2)\in \lc(L)^2$ tels que  $d_{L}^{M_Q}(L_1,L_2)\not=0$. Par abus, on a noté $(Q_1,Q_2)$ l'image de $(L_1,L_2)$ par cette section.
 Pour tout $\Lambda\in a_{L,\CC}^*$  assez général et $k=\dim(a_L^Q)$, on a (cf. \eqref{eq:derivee})
$$ v_{L,X}^{Q}(1,T)=\frac1{k!} \sum_{P\in \pc^Q(L)} \bg \Lambda,T_P\bd^k \theta_P^Q(\Lambda)$$
qui est clairement une fonction polynomiale de la variable $T$. 

Pour démontrer le théorème, il suffit de prouver la convergence de l'intégrale $J_{M,X}^{Q}(f)$ pour tout $Q\in \fc(M)$. Cette convergence, \emph{mutatis mutandis}, se démontre comme dans le cas local (cf. preuve de la proposition \ref{prop:cv-locale}). On se ramène tout d'abord à prouver la convergence de l'intégrale

$$\int_{G_X(\AAA)\back G(\AAA)} f(\Ad(g^{-1})X)  \|R_{P_1}(g)-R_{P_2}(g)\|^k \, dg^*$$
pour tout entier $k\geq 0$ et tout couple de sous-groupes paraboliques $(P_1,P_2)\in \pc(M)^2$ qui sont adjacents. À ce stade, comme dans la preuve de la proposition \ref{prop:cv-locale}, on reprend les notations du  §\ref{S:non-pos}. En particulier, on a introduit un entier naturel $r_2$ (cf. \eqref{eq:r2}). Il est important de remarquer que l'hypothèse selon laquelle $X$ est simple entraîne 
$$r_2>0.$$
En effet, on a $r_2=\sum_{j\in J_2} d_j$ pour un certain sous-ensemble $J_2\subset J$. On a observé au §\ref{S:non-pos} que $J_2$ est vide si et seulement si $\tilde{P}_1=\tilde{P}_2$. Or comme $X$ est simple, l'application  \eqref{eq:lappli} est injective (cf. proposition \ref{prop:P-R}), et l'égalité  $\tilde{P}_1=\tilde{P}_2$ entraîne $P_1=P_2$ ce qui n'est pas. Donc $J_2$ n'est pas vide ; on a donc $r_2>0$ puisque, sous l'hypothèse $X$ simple, on a $d_j>0$ pour tout $j$. La principale conséquence de $r_2>0$ est que le produit
$$\prod_{v\in V}Z_{r_1,v}(r_1+r_2)$$
converge et vaut $Z_{r_1}(r_1+r_2)$. On va pouvoir appliquer une variante adélique du lemme \ref{lem:alaRao}. Comme conséquence du lemme \ref{lem:Y} et de la formule \eqref{eq:RP-RPv}, on a l'égalité 
$$ \|R_{P_1}(g)-R_{P_2}(g)\|= |\log|\det(U_{1,3})| |  \|\al^\vee\|.$$
pour tout  $k\in K$ et $U\in \ngo_{\tilde{P}_1}(F)\cap \ngo_{\tilde{P}_2}(F)$ tel que $kUk^{-1}=g^{-1}X g$ (le vecteur $\al^\vee$ est défini au lemme \ref{lem:Y}). Le lemme \ref{lem:alaRao} implique alors qu'on a 

  $$
\int_{G_X(\AAA)\back G(\AAA)} f(\Ad(g^{-1})X)  \|R_{P_1}(g)-R_{P_2}(g)\|^k \, dg^*
$$
$$
=c\cdot \int_{I_{1,3}(\AAA)} \int_{\ngo_{\tilde{Q}}(\AAA)} \int_{K} f(k^{-1}(x+U)k)|\log|\det(x)| |^k |\det(x)|^{r_1+r_2} \, dk\, dU \, dx.
$$
où $I_{1,3}(\AAA)$ est muni du produit des mesures locales sur $I_{1,3}(F_v)$ et où
\begin{equation}
  \label{eq:cestc}
  c=  \|\al^\vee\|\cdot c_X\cdot Z_{r_1}(r_1+r_2)^{-1}
\end{equation}
et $c_X$ est la constante définie en \eqref{eq:formule-c-globale}.
On est donc ramené à démontrer la convergence de l'intégrale ci-dessus. Comme dans la preuve de la proposition \ref{prop:cv-locale}, on est ramené à démontrer la convergence de l'intégrale en dimension $1$

$$\int_{\AAA^\times} f(t) |\log|\det(t)| |^k |t|^{k'} \,dt$$
avec $k\geq 0$ et  $2\leq r_2+1\leq k'\leq r_2+r_1$ pour tout $f\in \Sc(\AAA)$. Cette dernière convergence est aiséee à obtenir. Cela conclut.
  \end{preuve}

\section{Asymptotique d'intégrales orbitales nilpotentes tronquées}\label{sec:asymp}

\begin{paragr}[Hypothèses.] Dans toute cette section, $F$ est un corps global. Les notations sont celles de la section \ref{sec:IOPglob}. On suppose que l'endomorphisme standard $X$ est simple (au sens de la définition \ref{def:simple}). 
  \end{paragr}

\begin{paragr}[Mesures de Haar.] --- \label{S:Haar-G1} Soit $H$ un sous-groupe de $G$. Soit
$$H(\AAA)^1=\{h\in H(\AAA) \mid    |\chi(h)|=1 \text{  pour tout } \chi\in X^*(H)   \}.$$
C'est un sous-groupe de $H(\AAA)$. Supposons $H(\AAA)$ munit d'une mesure de Haar alors le volume de $H(F)\back H(\AAA)^1$ est fini. On appliquera cette propriété au groupe $G$ lui-même, à ses sous-groupes paraboliques ou de Levi ou au centralisateur $G_X$ de $X$. On a fixé au §\ref{S:notations-Haar} des mesures de Haar sur les points adéliques de ces groupes. Pour tout sous-groupe de Levi $L\in \lc(T_0)$, on a une suite exacte 
$$1 \to L(\AAA)^1 \to L(\AAA) \to a_L \to 0$$
induite par $H_L$. On a fixé une mesure de Haar sur $a_L$ au §\ref{S:Haar-aM}. On munit alors  $L(\AAA)^1$ de la mesure de Haar qui donne sur le quotient  $L(\AAA)^1 \back L(\AAA)\simeq a_L$ la mesure sur $a_L$. Pour le sous-groupe $M\in \lc(T_0)$ défini au §\ref{S:surZ}, on a une suite exacte analogue 
  $$1 \to G_X(\AAA)^1 \to G_X(\AAA) \to a_M \to 0$$
induite cette fois par l'un des $R_P$ pour $P\in \pc(M)$, cf. lemme \ref{lem:act-centralisateur} et éq. \eqref{eq:RP-RPv}. La surjectivité est vraie ici car $X$ est supposé simple.

Soit $L\in \lc(T_0)$. Soit $A_{L}^\infty$ la composante neutre du groupe des $\RR$-points du sous-tore $\QQ$-déployé maximal de la restriction des scalaires de $F$ à $\QQ$ de $A_{M_1}$. Soit $A_L^{G,\infty}=A_L^\infty \cap G(\AAA)^1$. Ce sont des sous-groupes de $L(\AAA)$ sur lesquels la restriction de $H_L$ induit des isomorphismes $A_L^\infty\simeq a_L$ et $A_L^{G,\infty}\simeq a_L^G$. Par transport, ces isomorphisme munit $A_L^\infty$ et $A_L^{G,\infty}$ de mesures de Haar. On a alors des décompositions en produits directs $L(\AAA)=A_L^\infty \times L(\AAA)^1$ et $L(\AAA)\cap G(\AAA)^1=A_L^{G,\infty} \times L(\AAA)^1$ compatibles aux mesures de Haar.

\end{paragr}

\begin{paragr}[Cônes et fonctions caractéristiques.] --- Soit $P\in \fc(T_0)$. Soit $\hat{\tau}_P$, resp.  $\tau_P$,  la fonction caractéristique du cône formé des $H\in a_P$ tels que $\bg \varpi,H\bd>0$ pour tout $\varpi\in \hat{\Delta}_P$, resp.  $\bg \al, H\bd >0$ pour tout $\al\in \Delta_P$. On a $\tau_P\leq \hat{\tau}_P$. Soit $a_B^+\subset a_{T_0}$ le cône défini par la condition $\bg \al, H\bd\geq 0$ pour tout $\al \in \Delta_B$. On dira qu'un point $T\in a_B^+$ est assez régulier si les rapports $\bg \al, T\bd / \|T\|$ sont assez grands pour tout $\al \in \Delta_B$.

Pour tous $T\in a_B^+$ et $Q\in \fc(B)$, Arthur a défini, dans \cite{ar1} section 6, une fonction $F^Q(\cdot ,T)$ : c'est la fonction caractéristique d'un compact de $A_{M_Q}(\AAA)N_Q(\AAA)M_Q(F)\back G(\AAA)$ qui dépend de $T$. La construction de la fonction dépend non seulement du point $T$ mais aussi  du compact $K$, du sous-groupe de Borel $B$ et de son tore maximal $T_0$ et encore d'autres choix auxiliaires, dont   un point $T_-$ dans l'opposé de la chambre $a_B^+$. En fait, du moins pour des points $T$ assez réguliers, la fonction obtenue est indépendante des choix auxiliaires (cf. lemme 2.1 \cite{ar_unipvar}). On utilisera plusieurs fois l'identité suivante (cf. \cite{ar1} lemme 6.4) valable pour $T\in a_B^+$ assez régulier et tout $g\in A_G(\AAA)G(F)\back G(\AAA)$
  \begin{equation}
    \label{eq:HN}
    1=\sum_{Q\in \fc(B)}\sum_{\delta\in Q(F)\back G(F)} \tau_Q(H_Q(\delta g)-T) F^Q(\delta g,T).
  \end{equation}
  
\end{paragr}

\begin{paragr}[Énoncé des résultats.] ---  Soit $\of$ la $G$-orbite de $X$ et  $f\in \Sc(\ggo(\AAA))$. Soit $P$ un sous-groupe parabolique semi-standard et $g\in G(\AAA)$. Soit
  \begin{equation}
    \label{eq:KP}
     K_{P,\of}(f,g)=\sum_{Y\in \pgo(F)\cap \of(F), I_P^G(Y)=\of}  f(\Ad(g^{-1})Y),
   \end{equation}
   où la somme est convergente.

 Pour tout $T\in a_{T_0}$ soit

\begin{equation}
  \label{eq:KT}
  K^T_{\of}(f,g)=\sum_{B\subset P\subset G} (-1)^{\dim(a_P^G)} \sum_{\delta\in P(F)\back G(F)} \hat{\tau}_P(H_P(\delta g)-T_P) K_{P,\of}(f,\delta g)
\end{equation}
  où la somme externe est prise sur les sous-groupes paraboliques standard $P$. La somme interne sur $\delta$ est en fait à support finie. 
\end{paragr}

On observera que les fonctions $ K_{P,\of}(f,g)$ et $K^T_{\of}(f,g)$ de la variable $g$ sont  invariantes par le centre $A_G(\AAA)$ et à droite respectivement  par $P(F)$ et $G(F)$.

Soit
\begin{equation}
  \label{eq:JoT}
  J_\of^T(f)=\int_{G(F)\back G(\AAA)^1} K^T_{\of}(f,g) \, dg.
\end{equation}

 \begin{proposition}\label{thm:nilp-globale}
   Pour toute fonction $f\in \Sc(\ggo(\AAA))$ et tout $T\in a_{T_0}$, l'intégrale $J_\of^T(f)$ est absolument convergente. On a de plus l'égalité
\begin{equation}
  \label{eq:JoT-JMX}
   J_\of^T(f)= \vol(G_X(F)\back G_X(\AAA)^1) \cdot J_{M,X}^{G,T}(f)
 \end{equation}
 où dans le membre de droite le volume est fini et l'expression $J_{M,X}^{G,T}(f)$ est définie en \eqref{eq:JMXG}. En particulier $T\mapsto J_\of^T(f)$ est un polynôme.
\end{proposition}

Pour tout $T\in a_B^+$, il est naturel depuis les travaux d'Arthur d'introduire l'intégrale orbitale tronquée

\begin{equation}
  \label{eq:IO-tronquee}
  \int_{G(F)\back G(\AAA)^1} F^G(g,T) \sum_{Y\in \ggo(F)\cap \of(F)} f(g^{-1}Yg)   \, dg.
\end{equation}
 L'intégrale est évidemment convergente puisque  $F^G(g,T)$ induit sur $G(F)\back G(\AAA)^1$ la fonction caractéristique d'un compact. Le théorème suivant décrit l'asymptotique en $T$ de cette intégrale tronquée.

 \begin{theoreme}\label{thm:asymp}
   Pour tout $\eps'>0$, il existe $\eps>0$ et pour toute fonction $f\in \Sc(\ggo(\AAA))$ une constante $C>0$ tels que pour tout $T\in a_{B}^+$ qui vérifie
$$\bg \al,T\bd \geq \eps' \|T\|$$
pour tout $\al\in \Delta_{B}$, on a

  $$\int_{G(F)\back G(\AAA)^1} |K_{\of}^T(f,g) -F^G(g,T) \sum_{Y\in  \of(F)} f(g^{-1}Yg) |  \, dg \leq C\cdot \exp(-\eps \|T\|).$$
\end{theoreme}

\begin{remarque}
  D'après la proposition \ref{thm:nilp-globale} combinée avec le théorème \ref{thm:cv-pondere}, le théorème \ref{thm:asymp} implique que l'intégrale orbitale tronquée \eqref{eq:IO-tronquee} est asymptotique comme fonction de $T$ à un polynôme en $T$, en l'occurrence $J_{\of}^T(f)$. Sous cette forme faible, le théorème \ref{thm:asymp} est  une légère généralisation du théorème 4.2 de \cite{ar_unipvar} pour l'orbite $\of$. Cependant, on a mieux : on dispose de la  formule intégrale \eqref{eq:JoT-JMX} pour ce polynôme. 
\end{remarque}

La suite de cette section est consacrée à la démonstration de la proposition \ref{thm:nilp-globale} et du théorème  \ref{thm:asymp}.
\end{paragr}

\begin{paragr}[Démonstration de la  proposition \ref{thm:nilp-globale}.] --- La première étape dans la démonstration de la  proposition \ref{thm:nilp-globale} consiste à réécrire l'expression \eqref{eq:KT} comme une somme sur la $G(F)$-orbite de $X$. Auparavant, nous aurons besoin du lemme suivant.

\begin{lemme}\label{lem:lesLS}Soit $\of$ la $G$-orbite d'un endomorphisme standard $X$ (l'hypothèse $X$ simple n'est pas requise ici). Soit $\fc_\of$ l'ensemble des couples $(Q,\of')$ tels que
  \begin{itemize}
  \item $Q\in \fc(B)$
  \item $\of'$ est une $M_Q$-orbite nilpotente dans $\mgo_Q$
  \item $I_Q^G(\of')=\of$.
  \end{itemize}

Soit
\begin{equation}
  \label{eq:LS-Fo}
      \lc\Sc(X) \to \fc_\of
    \end{equation}
    l'application définie ainsi : l'image de $P\in \lc\Sc(X)$ est le couple $(Q,\of')$ où 
\begin{itemize}
\item $Q$ est l'unique sous-groupe parabolique standard conjugué à $P$ ;
\item soit $g\in G$ tel que $gPg^{-1}=Q$ ; la $M_Q$-orbite $\of'$ contient la projection  de $gXg^{-1}\in \qgo=\mgo_Q\oplus\ngo_Q$ sur $\mgo_Q$. 
\end{itemize}
L'application \eqref{eq:LS-Fo} est une bijection.
  \end{lemme}

\begin{preuve}
L'application est bien définie : l'élément $gXg^{-1}$ est uniquement défini à $Q$-conjugaison près. Sa projection sur $\mgo_Q$ définit donc une unique $M_Q$-orbite $\of'$. La condition $I_P^G(X)=\of$ implique $I_Q^G (gXg^{-1})=\of$ d'où      $I_Q^G(\of')=\of$.

L'application est injective. En effet, si $P_1$ et $P_2$ ont même image $(Q,\of')$, ils sont conjugués par $g_1$ et $g_2$ à $Q$. De plus, $g_1Xg_1^{-1}$ et $g_2Xg_2^{-1}$ appartiennent à la même $Q$-classe de conjugaison à savoir $(\of'\oplus \ngo_Q)\cap \of$. Donc, il existe $q\in Q$ tel que  $h=g_1^{-1}q g_2\in G_X\subset P_1\cap P_2$. On a donc
$$P_2=g_2^{-1}Qg_2= (q^{-1}g_1 h)^{-1}Q q^{-1}g_1 h=h^{-1} g_1^{-1} Q g_1 h=h^{-1} P_1 h=P_1.$$

L'application est surjective. Soit $(Q,\of')\in \fc_\of$ et $Y'\in \of'$. Il existe $U\in \ngo_Q$  tel que $Y=Y'+U\in \of$. Soit $g\in G$ tel que $gXg^{-1}=Y$. Alors $g^{-1}Pg\in \lc\Sc(X)$. 

  \end{preuve}

Pour tout $g\in G(\AAA)$, on pose
\begin{equation}
  \label{eq:sigma}
  \sigma^T(g)= \sum_{Q\in \fc(M)} (-1)^{\dim(a_Q^G)}  \hat{\tau}_Q(R_Q(g)-T_Q).
\end{equation}

\begin{lemme} \label{lem:reecriture}
  Pour tout $g\in G(\AAA)$, 
$$K^T_{\of}(f,g)= \sum_{\delta\in G_X(F)\back G(F)}\sigma^T(\delta g)f(\Ad( (\delta g)^{-1})X).$$
\end{lemme}

\begin{preuve}
  On a des bijections

$$
\begin{array}{ccccc}
  \fc(M) &\to& \lc\Sc(X) &\to& \fc_{\of}\\
Q&\mapsto &\tilde{Q} &\mapsto &(P,\of')
\end{array}
$$
donnée par \eqref{eq:lappli2} (car $X$ est simple) et  \eqref{eq:LS-Fo}. Soit $P$ un sous-groupe parabolique standard. Soit $\fc(M)[P]\subset \fc(M)$ l'ensemble des éléments de $\fc(M)$ qui ont pour image dans $\fc_{\of}$ un couple de la forme $(P,\of')$. Autrement dit $\fc(M)[P]$ est l'ensemble des éléments de $\fc(M)$ qui sont conjugués à $P$. Pour tout $Q\in \fc(M)[P]$, soit $\tilde{Q}$ son image par \eqref{eq:lappli2} et $w_Q$ et $s_Q$ des éléments de $W$ tels que  $\tilde{Q}=w_Q^{-1}Qw_Q$ et $P=s_Q^{-1}\tilde{Q}s_Q$. 
On a 
$$ K_{P,\of}(f,g)= \sum_{\{\of'\mid (P,\of')\in \fc_\of \}} \sum_{Y \in (\of'(F)\oplus\ngo_P(F))\cap \of(F) } f(\Ad(g^{-1})Y).
$$
Si $Q$ a pour image $(P,\of')$, l'ensemble sur lequel on somme $Y$ est encore le conjugué sous $s_Q^{-1}$ de la $\tilde{Q}(F)$-orbite de $X$.   On a donc
$$ K_{P,\of}(f,g)= \sum_{Q\in \fc(M)[P]} \sum_{\nu\in G_X(F)\back \tilde{Q}(F)} f(\Ad(\nu s_Q g) ^{-1} X).$$

On a alors
$$  K^T_{\of}(f,g)= $$
$$ \sum_{B\subset P\subset G} (-1)^{\dim(a_P^G)} \sum_{\delta\in P(F)\back G(F)} \hat{\tau}_P(H_P(\delta g)-T_P) \sum_{Q\in \fc(M)[P]} \sum_{\nu\in G_X(F)\back \tilde{Q}(F)} f(\Ad(\nu s_Q \delta g) ^{-1} X)
$$
$$=
\sum_{B\subset P\subset G} (-1)^{\dim(a_P^G)} \sum_{Q\in \fc(M)[P]} \sum_{\nu\in G_X(F)\back G(F)} \hat{\tau}_P(H_P(s_Q^{-1}\delta g)-T_P) f(\Ad(\delta g) ^{-1} X)
$$
$$= \sum_{\nu\in G_X(F)\back G(F)}f(\Ad(\delta g) ^{-1} X) [  \sum_{B\subset P\subset G} (-1)^{\dim(a_P^G)} \sum_{Q\in \fc(M)[P]}    \hat{\tau}_P(H_P(s_Q^{-1}\delta g)-T_P)   ].
$$

Il suffit pour conclure de prouver que l'expression entre crochets est égale à $\sigma^T(\delta g)$. Pour cela, il suffit d'observer qu'on a pour $Q\in \fc(M)[P]$
$$H_P(s_Q^{-1}g)=s_Q^{-1}\cdot H_{\tilde{Q}}(g)=(w_Qs_Q)^{-1}\cdot R_Q(g)$$
et 
$$\hat{\tau}_P(H_P(s_Q^{-1}\delta g)-T_P)= \hat{\tau}_Q( R_Q( g)- T_Q).$$
 \end{preuve}

Le lemme suivant va nous donner une prise sur la fonction $\sigma^T$.

 \begin{lemme}\label{lem:cvx}(Arthur) Soit $\yc=(Y_P)_{P\in \pc(M)}$ une famille orthogonale.
    \begin{enumerate}
    \item L'application
$$H \in a_M^G \mapsto  \sum_{P\in \fc(M)}  (-1)^{\dim(a_P^G)}  \hat{\tau}_P(H -Y_P)$$
s'annule en $H$ sauf si $H$ appartient à l'image de l'enveloppe convexe des points $(Y_P)_{P\in \pc(M)}$ par la projection $a_M\to a_M^G$. Si, de plus, la famille $\yc$ est positive, alors cette application est la fonction caractéristique de la projection sur $a_M^G$ de  cette enveloppe convexe.
     \item L'intégrale 
$$\vol(\yc)= \int_{a_M^G}   \sum_{P\in \fc(M)}  (-1)^{\dim(a_P^G)}  \hat{\tau}_P(H -Y_P) \,dH       $$  est égale à la valeur $v_M$ associée à la $(G,M)$-famille définie pour $P\in \pc(M)$ par 
$$v_P(\la,g)=\exp(\bg \la, Y_P\bd).$$
    \end{enumerate}
  \end{lemme}
  
\begin{preuve} Il s'agit juste de réécrire l'application de façon à pouvoir citer Arthur. Soit $P_1\in \pc(M)$. On a une partition
$$\fc(M)=\coprod_{P\in \pc(M)} \fc(P,P_1)$$
où $$\fc(P,P_1)=\{Q\in \fc(P) \mid P\cap M_Q= P_1\cap M_Q \}.$$
 Soit $\Lambda$ un point de la chambre ouverte $a_{P_1}^+$ définie par la positivité sur $\Delta_{P_1}$. Pour tout $P\in \pc(M)$, soit $\varphi_P^\Lambda$ la fonction caractéristique des $H\in a_M$ tels que
\begin{itemize}
\item $\bg \varpi_\al,H\bd >0$ pour $\al\in \Delta_P$ tel $\bg \Lambda , \al^\vee \bd <0$ ;
\item $\bg \varpi_\al,H\bd \leq 0$ pour $\al\in \Delta_P$ tel $\bg \Lambda , \al^\vee \bd >0$.
\end{itemize}
On a alors
\begin{equation}
  \label{eq:tau-phi}
  \sum_{Q\in \fc(P,P_1)} (-1)^{\dim(a_Q^G)} \hat{\tau}_P(H -Y_Q)= (-1)^{|\Delta_P^\Lambda|}\varphi_P^\Lambda(H-Y_P),
\end{equation}
où $\Delta_P^\Lambda\subset \Delta_P$  est formé des racines $\al$ telles que $\bg \Lambda , \al^\vee \bd <0$.

Il s'ensuit qu'on a
\begin{eqnarray*}
  \sum_{P\in \fc(M)}  (-1)^{\dim(a_P^G)}  \hat{\tau}_P(H -Y_P)&=& \sum_{P\in \pc(M)} \sum_{Q\in \fc(P,P_1)}  (-1)^{\dim(a_Q^G)}\hat{\tau}_P(H -Y_Q)\\
&=&\sum_{P\in \pc(M)} (-1)^{|\Delta_P^\Lambda|}\varphi_P^\Lambda(H-Y_P).
\end{eqnarray*}
L'assertion 1 est alors un résultat d'Arthur et Langlands (cf. \cite{localtrace} p. 22 eq. (3.8) \emph{infra}). L'assertion 2 s'en déduit également (cf. \cite{dis_series}).
    \end{preuve}

On est maintenant en mesure de donner la démonstration de la proposition \ref{thm:nilp-globale}. Traitons d'abord la convergence de $J_\of^T(f)$. D'après le lemme \ref{lem:reecriture}, il s'agit de prouver la convergence  de l'intégrale suivante
$$\int_{G_X(F)\back G(\AAA)^1} |\sigma^T(g)|\cdot |f(\Ad(g^{-1})X )| \, dg$$
$$=\int_{G_X(\AAA)\back G(\AAA)}   |f(\Ad(g^{-1})X )|  \int_{G_X(F)\back (G_X(\AAA)\cap G(\AAA)^1)} |\sigma^T(hg)| \, dh  \, dg.$$
On a muni  $G_X(\AAA)\cap G(\AAA)^1$ de la mesure de Haar de sorte que l'identification naturelle $(G_X(\AAA)\cap G(\AAA)^1) \back G(\AAA)^1 \simeq G_X(\AAA)\back G(\AAA)$ soit compatible aux mesures quotients.

Traitons l'intégrale intérieure. Par le choix des mesures et le lemme \ref{lem:act-centralisateur}, on a 

$$\int_{G_X(F)\back (G_X(\AAA)\cap G(\AAA)^1)}  |\sigma^T(hg)|\, dh$$
$$= \vol(G_X(F)\back G_X(\AAA)^1) \int_{a_{M}^G}  |\sum_{Q\in \fc(M)} (-1)^{\dim(a_Q^G)}  \hat{\tau}_Q(H+R_Q(g)-T_Q)|\, dH.$$
Le volume qui apparaît ci-dessus est fini (cf. §\ref{S:Haar-G1}). D'après le lemme  \ref{lem:cvx}, l'expression
 $$|\sum_{Q\in \fc(M)} (-1)^{\dim(a_Q^G)}  \hat{\tau}_Q(H+R_Q(g)-T_Q)|$$
est nulle sauf si $H$ appartient à l'image de l'enveloppe convexe des points $(T_P-R_P(g))_{P\in \pc(M)}$ par la projection $a_M\to a_M^G$. Il est facile d'en déduire qu'il existe une constante $C>0$ telle que pour $k=\dim(a_M^G)$, on a 
$$ \int_{a_{M}^G}  |\sum_{Q\in \fc(M)} (-1)^{\dim(a_Q^G)}  \hat{\tau}_Q(H+R_Q(g)-T_Q)|\, dH \leq C\cdot \sum_{(P_1,P_2)\in \pc^{}(M)^{\text{adj}}} \|R_{P_1}(g)-T_{P_1}-R_{P_2}(g)+T_{P_2}\|^k$$
où la somme porte sur les couples de paraboliques adjacents. La convergence s'obtient alors comme dans la preuve du théorème \ref{thm:cv-pondere}. Une fois la convergence  acquise, le reste de la proposition résulte de la formule
 $$
v_{M,X}^G(g,T)=\int_{a_{M}^G}  \sum_{Q\in \fc(M)} (-1)^{\dim(a_Q^G)}  \hat{\tau}_Q(H+R_Q(g)-T_Q)\, dH
$$
(cf. lemme \ref{lem:cvx} assertion 2).
\end{paragr}

\begin{paragr}[Début de la démonstration du théorème \ref{thm:asymp}.] --- 
En utilisant  le lemme \ref{lem:reecriture}, on a l'égalité 

$$K_{\of}^T(f,g) -F^G(g,T) \sum_{Y\in  \of(F)}  f(\Ad(g^{-1})Y)= \sum_{\delta \in G_X(F)\back G(F)} [\sigma^{T}(\delta g)-F^T(\delta g)]f(\Ad(g^{-1})X) $$
On a donc la majoration
 $$  \int_{G(F)\back G(\AAA)^1} |K_{\of}^T(f,g) -F^G(g,T) \sum_{Y\in \ggo(F)\cap \of(F)}  f(\Ad(g^{-1})Y) |  \, dg$$
$$\leq \int_{G_X(F)\back G(\AAA)^1} |\sigma^T(g)-F^G(g,T)| \cdot | f(\Ad(g^{-1})X) |  \, dg.$$

Au moins pour des  $T$ assez réguliers, on a la propriété suivante :  pour tout $g\in G(\AAA)$ tel que  $F^G(g,T)=1$ et   tout sous-groupe parabolique standard $Q$ \emph{propre}, on a 
$$\hat{\tau}_Q(H_Q(g)-T_Q)=0.$$
Cette propriété vaut encore pour les sous-groupe paraboliques semi-standard. En particulier, pour un tel $g$, on a $\sigma^T(g)=1$. On a donc la majoration suivante pour tout $g\in G(\AAA)$
$$ |\sigma^T(g)-F^G(g,T)| \leq (1-F^{G}(g,T)) |\sigma^T(g)|.$$

Introduisons la fonction $I^T$ sur $A_G(\AAA)G_X(\AAA)\back G(\AAA)$ qui est la fonction caractéristique des $g\in A_G(\AAA)G_X(\AAA)\back G(\AAA)$ tels que la famille $(T_P-R_P(g))_{P\in \pc(M)}$ soit positive. Lorsqu'on a $I^T(g)=1$, l'expression $\sigma^T(g)$ vaut $0$ ou $1$ d'après le lemme \ref{lem:cvx}.

On est donc ramené à majorer les deux intégrales suivantes, ce qui sera fait dans les deux paragraphes qui suivent.

 \begin{equation}
   \label{eq:first-int}
    \int_{G_X(F)\back G(\AAA)^1} (1-F^{G}(g,T)) \cdot\sigma^T(g)\cdot I^T(g)  \cdot |f(g^{-1}Xg) |  \, dg
 \end{equation}
et 
\begin{equation}
   \label{eq:second-int}
    \int_{G_X(F)\back G(\AAA)^1} (1-F^{G}(g,T)) \cdot|\sigma^T(g)|\cdot (1-I^T(g))  \cdot |f(g^{-1}Xg) |  \, dg.
 \end{equation}
\end{paragr}

\begin{paragr}[Majoration de \eqref{eq:first-int}.] --- L'expression  \eqref{eq:first-int} se réécrit
$$\int_{G(F)\back G(\AAA)^1}  (1-F^{G}(g,T))   \sum_{\delta \in  G_X(F)\back G(F)}\sigma^T(\delta g)\cdot I^T(\delta g)  \cdot |f((\delta g)^{-1}X\delta g) |  \, dg.$$
L'égalité \eqref{eq:HN} donne un développement indexé par $Q\in\fc(B)$ du facteur  ($1-F^{G}(g,T))$. On est donc ramené à majorer pour tout  $Q\in\fc(B)$
  \begin{equation}
   \label{eq:first-int2}
 \int_{Q(F)\back G(\AAA)^1}  \tau_Q(H_Q(g)-T) F^Q(g,T)   \sum_{\delta \in  G_X(F)\back G(F)}\sigma^T(\delta g)\cdot I^T(\delta g)  \cdot |f(\Ad (\delta g)^{-1}X) |  \, dg.
 \end{equation}

L'égalité \eqref{eq:HN} a un analogue relatif à $Q$ : soit $T_1\in a_{B}^+$ assez régulier, fixé une fois pour toutes :  pour tout $A_{M_Q}(\AAA)N_Q(\AAA)M_Q(F)\back G(\AAA)$ on a 
$$1=\sum_{\{Q_1 \mid B \subset Q_1 \subset Q\}  } \sum_{\nu \in Q_1(F)\back Q(F)} F^{Q_1}(\nu g,T_1) \tau_{Q_1}^Q(H_{Q_1}(\nu g)-T_1).$$
En conséquence, l'expression \eqref{eq:first-int2} est alors la somme sur les sous-groupes paraboliques standard $Q_1$ inclus dans $Q$ de 

 \begin{equation}
   \label{eq:first-int3}
\int_{Q_1(F)\back G(\AAA)^1}  F^{Q_1,Q}(g,T_1,T)  \cdot \tau_{Q_1}^{Q_1,Q}(g, T_1,T)\sum_{\delta \in  G_X(F)\back G(F)}\sigma^T(\delta g)\cdot I^T(\delta g)  \cdot |f((\delta g)^{-1}X\delta g) |  \, dg
 \end{equation}
où l'on pose
$$F^{Q_1,Q}(g,T_1,T)= F^{Q_1}(g,T_1)  F^Q(g,T)$$
et
$$\tau_{Q_1}^{Q_1,Q}(g, T_1,T) = \tau_{Q_1}^Q(H_{Q_1}(g)-T_1)   \tau_Q(H_Q(g)-T).$$
Dans la suite on fixe un tel $Q_1$ et il nous suffit de majorer \eqref{eq:first-int3}.

Pour tout sous-groupe  parabolique $Q_2$ contenant $Q$, soit
\begin{equation}
  \label{eq:q2}
  \qgo_2'=\qgo_2 - \bigcup_{Q\subset R \subsetneq Q_2}\rgo
\end{equation}
Soit 
$$ \Upsilon_{Q_2} =\{\delta \in  G_X(F)\back G(F)\mid  \delta^{-1}X\delta  \in \qgo_2'(F) \}.$$
La collection des $\Upsilon_{Q_2}$, pour $Q\subset Q_2 \subset G$ forme une partition de  $G_X(F)\back G(F)$.

L'intégrale \eqref{eq:first-int2} se majore par la somme sur $Q\subset Q_2 \subset G$ de 
\begin{equation}
   \label{eq:first-int4}
\int_{Q_1(F)\back G(\AAA)^1}  F^{Q_1,Q}(g,T_1,T)  \cdot \tau_{Q_1}^{Q_1,Q}(g, T_1,T)  \sum_{\delta \in \Upsilon_{Q_2} }\sigma^T(\delta g)\cdot I^T(\delta g)  \cdot |f(\Ad(\delta g)^{-1}X) |  \, dg
 \end{equation}

Dans la suite, on fixe un tel $Q_2$.

 \begin{lemme}
   Soit $\delta \in  \Upsilon_{Q_2}$ tel que l'expression
   \begin{equation}
     \label{eq:unproduit}
    \tau_Q(H_Q(g)-T_Q) \cdot  \sigma^T(\delta g)\cdot I^T(\delta g)
  \end{equation}
  soit non nulle pour un élément $g\in G(\AAA)$. Alors pour tout sous-groupe parabolique maximal $R$ de $G$ contenant $Q_2$, on a  $R\notin \lc\Sc(\delta^{-1}X\delta)$.
 \end{lemme}

 \begin{preuve}
On va supposer la conclusion en défaut et on va trouver une contradiction. On suppose donc  qu'il existe $R$, un  sous-groupe parabolique maximal de $G$ contenant $Q_2$, tel que $R\in \lc\Sc(\delta^{-1}X\delta)$. Il suffit en fait de trouver un sous-groupe parabolique maximal $S\in \fc(M)$ tel que
\begin{equation}
  \label{eq:Snonnul}
  \hat{\tau}_S(R_S(\delta g)-T_S)=1.
\end{equation}
En effet, la non-annulation de \eqref{eq:unproduit} donc de  $I^T(\delta g)$ entraîne que la famille $(T_P-R_P(\delta g))_{P\in \pc(M)}$ est orthogonale positive. Donc le lemme \ref{lem:cvx} et la non-nullité de $\sigma^T(\delta g)$ impliquent que $0$ appartient à l'enveloppe convexe des projections sur $a_M^G$ des points $T_P-R_P(\delta g)$ pour ${P\in \pc(M)}$. En particulier, pour tout élément maximal $S\in \fc(M)$, on doit avoir (cf.  \cite{dis_series} lemme 3.2), pour $\varpi$ l'unique élément de $\hat{\Delta}_S$,
$$0\leq \bg \varpi, T_S -R_S(\delta g)\bd.$$ 
ce qui  contredit manifestement \eqref{eq:Snonnul}.
Soit  $S\in \fc(M)$ dont l'image dans $\lc(X)$ par l'application \eqref{eq:lappli2} est  $\tilde{S}=\delta R \delta^{-1}$. Pour un tel $S$, les groupes qui interviennent sont semi-standard et  conjugués ; ils sont donc conjugués sous $W$ et  on a 

\begin{equation}
  \label{eq:longue-eg}
  \hat{\tau}_S(R_S(\delta g)-T_S)=\hat{\tau}_{\tilde{S}}(H_{\tilde{S}}(\delta g)-T_{\tilde{S}})=\hat{\tau}_R(H_R(g)-T_R).
\end{equation}

Or $R$ contient $Q_2$ donc $Q$ et, puisque l'expression \eqref{eq:unproduit} n'est pas nulle, on a ${\tau}_Q(H_Q(g)-T_Q)=1$ d'où \emph{a fortiori} $\hat{\tau}_Q(H_Q(g)-T_Q)=1$ et $\hat{\tau}_R(H_R(g)-T_R)=1$. Mais alors l'égalité \eqref{eq:longue-eg} entraîne \eqref{eq:Snonnul}.
 \end{preuve}

En tenant compte du lemme précédent, on peut dans \eqref{eq:first-int4} tout d'abord restreindre la somme aux $\delta\in \Upsilon_{Q_2}$ tels que  pour tout sous-groupe parabolique maximal $R$ de $G$ contenant $Q_2$, on a  $R\notin \lc\Sc(\delta^{-1}X\delta)$. On majore ensuite  les fonctions caractéristiques $\sigma^T$ et $I^T$ par $1$. Ainsi on peut majorer \eqref{eq:first-int4} par 
\begin{equation}
   \label{eq:first-int5}
\int_{Z(\AAA)Q_1(F)\back G(\AAA)}  F^{Q_1,Q}(g,T_1,T)  \cdot \tau_{Q_1}^{Q_1,Q}(g, T_1,T)  \sum_{Y \in \of_{2}}  |f(g^{-1}Y g) |  \, dg
 \end{equation}
où l'on introduit l'ensemble $\of_2$ formé des éléments $Y\in \qgo'_2(F)\cap \of(F)$ tels que $ Y\notin I_R^G(Y)$ pour tout sous-groupe parabolique maximal $R$ contenant $Q_2$. On observera que $\of_2$ est stable par $Q_2(F)$-conjugaison.

 Pour alléger un peu, on pose $M_1=M_{Q_1}$,  $N_1=N_{Q_1}$ et $A_1^{G,\infty}=A_{M_1}^{G,\infty}$. En utilisant la décomposition d'Iwasawa, on voit que \eqref{eq:first-int5} est égal à 
\begin{eqnarray}
   \label{eq:first-int6}\nonumber
\int_{M_1(F)\back M_1(\AAA)^1} \int_{A_1^{G,\infty}}  \int_{N_1(F)\back N_1(\AAA)} \exp(-\bg 2\rho_{Q_1},H_{Q_1}(a)\bd ) F^{Q_1,Q}(nam,T_1,T)   \\
\tau_{Q_1}^{Q_1,Q}(a, T_1,T)  \int_{K} \sum_{Y \in \of_{2}}  |f((namk)^{-1}Y namk) |  \, dk .
 \end{eqnarray}
Sous la condition  $F^{Q_1,Q}(nam,T_1,T)\not=0$, on a en particulier $F^{Q_1}(m,T_1)\not=0$ et $m$ reste dans un compact fixé, qui dépend de $T_1$, inclus dans  $M_1(F)\back M_1(\AAA)^1$. Le quotient $N_1(F)\back N_1(\AAA)$ est compact. Écrivons alors 
$$\qgo_2=\mgo_1\oplus \bigoplus_{\al \in \Sigma(Q_2,A_1)}\ggo_{\al}$$
la décomposition de $\qgo_2$ en espaces propres sous l'action par conjugaison du centre  $A_1$ de $M_1$. Pour tout $Y\in \qgo_2(F)$, on écrit 
$$Y=Y_1+\sum_{\al \in  \Sigma(Q_2,A_1)}Y_\al$$ 
selon cette décomposition. Il existe alors des fonctions de Bruhat-Schwartz positives $\varphi_1$ sur $\mgo_1(\AAA)$ et  $\varphi_\al$ sur $\ggo_\al(\AAA)$ pour tout $\al \in  \Sigma(Q_2,T_0)$ telles qu'on ait 
$$\sum_{Y \in \of_{2}}  |f((namk)^{-1}Y namk) | \leq \sum_{Y \in \of_{2}} \varphi_1(Y_1)\prod_{\al\in  \Sigma(Q_2,A_1) }\varphi_{\al}(\al(a)^{-1}Y_\al)$$
pour tout $k\in K$, $n\in N_1(F)\back N_1(\AAA)$ et $m\in M_1(F)\back M_1(\AAA)^1$ tel que $F^{Q_1}(m,T_1)\not=0$.

Introduisons 
$$
c_-=\inf_m \inf_{\al\in \Delta_B^Q} (-\bg \al, H_{B}(m)\bd ) \ \text{    et     }\   c_+=\sup_m\sup_{\varpi \in \hat{\Delta}^{Q}_B}( -\bg \varpi, H_{B}(m)\bd)  
$$
où les bornes inférieure et supérieure  $c_+$ et $c_-$ sont  prises sur  les $m\in  M_1(F)\back M_1(\AAA)^1$ tels que $F^{Q_1}(m,T_1)\not=0$. Elles  ne dépendent que de $T_1$. 

Soit $a_{M_1}^G(Q,T)$ l'ensemble des $H\in a_{M_1}^G$ qui vérifient les conditions ci-dessous

\begin{enumerate}
\item pour tout $\al \in \Delta^{Q}_B$, $\bg \al, H \bd \geq \bg \al, T_-\bd +c_-$ ;
\item pour tout $\varpi\in \hat{\Delta}^{Q}_B$, $\bg \varpi, H \bd \leq \bg \varpi, T\bd +c_+$ ;
\item pour tout $\al \in \Delta_{Q_1}^Q$, $\bg \al, H \bd > \bg \al, T_1\bd$ ;
\item pour tout $\al \in \Delta_{Q}$, $\bg \al, H \bd > \bg \al, T\bd$.
\end{enumerate}
La condition $F^Q(nam, T)\not=0$ implique que $H_{Q_1}(a)$ vérifie les conditions 1 et 2. La condition $\tau_{Q_1}^{Q_1,Q}(a,T_1,T)\not=0$ implique que $H_{Q_1}(a)$ satisfait les conditions 3 et 4. Soit $A_1^{G,\infty}(Q,T)$  l'image inverse par $H_{M_1}$ dans $A_1^{G,\infty}$ de $a_1^G(Q,T)$. La condition $F^Q(nam, T)\not=0$ implique donc $H_{Q_1}(a)\in A_1^{G,\infty}(Q,T)$.

Avec le lemme \ref{lem:une-maj} ci-dessous, on voit, qu'à une constante près, l'intégrale \eqref{eq:first-int6} est majorée par 
\begin{equation}
  \label{eq:first-integ7}
  \int_{ a_1^{G}(Q,T)  } \prod_{\al \in \Delta_{1}-\Delta_{1}^{Q}} \exp(-c_\al \bg \al, H\bd) \, dH.
\end{equation}
On a une partition $\Delta_1^Q\cup (\Delta_1-\Delta_1^Q)$ de la base $\Delta_1^Q$ de $a_1^{G,*}$. En prenant les espaces orthogonaux on a une décomposition en somme directe notée
$$a_1^{G}=a_Q^G\oplus b.$$
Soit $H\in a_1^G$ qu'on écrit $H'+H''$ selon cette décomposition. Alors $H\in  a_1^{G}(Q,T)$ si et seulement si $H''$ appartient à l'ensemble $b(Q,T)$ des éléments de $b$ qui vérifient les conditions 1 à 3 ci-dessus et $H'+H''$ vérifie la condition 4 ci-dessus. Supposons donc $H\in  a_1^{G}(Q,T)$. On a 
$$H=H_Q+\sum_{\al \in \Delta^{Q}}  \bg\varpi_\al,H \bd \al^\vee$$
où $H_Q$ est la projection de $H$ sur $a_Q^G$ selon la décomposition $a_{1}^G=a_{1}^Q\oplus a_Q^G$. En utilisant l'inégalité $\bg \al, \beta^\vee \bd\leq 0$ pour $\al\not=\beta$ et la condition 2 pour $H''$, il vient  pour tout $\beta\in \Delta_{Q_1}-\Delta_{Q_1}^Q$
\begin{eqnarray}\label{eq:minor-beta}
  \bg \beta, H' \bd &=& \bg \beta, H \bd \\
&=& \bg \beta, H_Q\bd + \sum_{\al \in \Delta^{Q}}  \bg\varpi_\al,H \bd \bg \beta, \al^\vee\bd \nonumber \\
&=& \bg \beta, H_Q\bd + \sum_{\al \in \Delta^{Q}}  \bg\varpi_\al,H'' \bd \bg \beta, \al^\vee\bd \nonumber\\
&\geq & \bg \beta, T_Q \bd + \sum_{\al \in \Delta^{Q}}  (\bg\varpi_\al,T \bd +c_+ )  \bg\beta, \al^\vee\bd \nonumber\\
&= &   \bg \beta, T_+\bd .\nonumber
\end{eqnarray}
où l'on $T_+=T+ c_+ \sum_{\al \in \Delta^{Q}}    \al^\vee$. Soit $a_Q^G(T)$ l'ensemble des $H'\in a_Q^G$ tel que pour tout  $\beta\in \Delta_{Q_1}-\Delta_{Q_1}^Q$ on ait  $\bg \beta, H' \bd\geq  \bg \beta, T\bd$.

L'intégrale \eqref{eq:first-integ7} se majore alors par 

\begin{equation}
  \label{eq:first-integ8}
  \vol(b(Q,T))  \int_{ a_Q^{G}(T_+)  } \prod_{\al \in \Delta_{1}-\Delta_{1}^{Q}} \exp(-c_\al \bg \al, H \bd) \, dH.
\end{equation}
qui vaut à une constante près qui ne dépend que des choix des mesures
\begin{equation}
  \label{eq:first-integ9}
 = \vol(b(Q,T))  \prod_{\al \in \Delta_{1}-\Delta_{1}^{Q}}c_\al^{-1} \exp(-c_\al \bg \al, T_+\bd) .
\end{equation}
Le volume  $\vol(b(Q,T)) $ dépend polynomialement de $T$. Il existe $C>0$ et $\eps>0$ tels que pour tout $T$ tel que  $\bg \al, T \bd\geq \eps' \| T\|$, l'expression \eqref{eq:first-integ9} se majore par 
$$C \cdot \exp(-\eps \|T\|).$$
 Cela conclut la majoration de  \eqref{eq:first-int} modulo le lemme suivant.

\begin{lemme}\label{lem:une-maj}  Pour tout $c>0$, il existe $C>0$ et 
  \begin{itemize}
  \item $c_\al >c$ pour  tout $\al\in \Delta_{1}^{2}-\Delta_1^Q$ ;
  \item $c_\al >0$  pour  tout $\al \in \Delta_{1}-\Delta_{1}^{2}$,
  \end{itemize}
tels que pour $a \in    A_1^{G,\infty}(Q,T)$ l'expression 
$$ \exp(-\bg 2\rho_{Q_1},H_{Q_1}(a)\bd )\sum_{Y \in \of_{2}} \varphi_1(Y_1)\prod_{\al\in  \Sigma(Q_2,A_1) }\varphi_{\al}(\al(a)^{-1}Y_\al)$$
est majorée par 
$$ C\cdot \prod_{\al \in \Delta_{1}-\Delta_{1}^{Q}}\exp(-c_\al \bg \al, H_{Q_1}(a)\bd).$$

\end{lemme}

\begin{remarque}
  Rappelons que $Q$ est un sous-groupe parabolique propre de $G$. L'ensemble $\Delta_{1}-\Delta_{1}^{Q}$ n'est donc pas vide  et la majoration obtenue n'est jamais triviale.
\end{remarque}

\begin{preuve} Soit $\beta \in  \Delta_{1}-\Delta_{1}^{2}$. Il suffit de démontrer la variante plus faible du lemme où l'on suppose  $c_\al=0$ si  $\al \in \Delta_{1}-\Delta_{1}^{2}$ et $\al\not=\beta$. En faisant varier $\beta$, on obtient aisément  le résultat cherché. La racine $\beta$ détermine un sous-groupe parabolique maximal $R$ contenant $Q_2$. On va utiliser le lemme suivant, dont la démontration (omise) est en tout point semblable à celle de la proposition 5.3.1 de \cite{scfhn}.

\begin{lemme}
  \label{lem:zero-de-py} Soit $d$ un entier. Il existe une constante $c_R>0$ telle que pour tout $T\in a_{B_0}^+$ et $a\in A_1^{G,\infty}(Q,T)$ et toute famille de polynômes non tous nuls $g=(g_i)_{i\in I} \in F[\ngo_R]$ sur $\ngo_R$ dont le degré total de chaque $g_i$ est borné par $d$, on a l'inégalité
$$\exp(-\bg 2\rho_R,H_{Q_1}(a)\bd) \sum_{Y\in \ngo_R(F)\cap \vc(g)} \prod_{\al \in \Sigma(N_R,A_1)}\varphi_{\al}(\al(a)^{-1}Y_\al)  \leq c_R\cdot \exp(-\bg \al,H_{Q_1}(a)\bd)$$
où  $\vc(g)$ est le fermé de $\rgo$ défini par la famille $g$ de polynômes. 
\end{lemme}

L'expression 
$$ \exp(-\bg 2\rho_{Q_1},H_{Q_1}(a)\bd )\sum_{Y \in \of_{2}} \varphi_1(Y_1)\prod_{\al\in  \Sigma(Q_2,A_1) }\varphi_{\al}(\al(a)^{-1}Y_\al)$$
se majore par
\begin{equation}
  \label{eq:majI}
  \exp(-\bg 2\rho_{Q_1},H_{Q_1}(a)\bd )\sum_{(Y,U)  } \varphi_1(Y_1)\prod_{\al\in  \Sigma(Q_2\cap M_R,A_1) }\varphi_{\al}(\al(a)^{-1}Y_\al)  \prod_{\al\in  \Sigma(N_R,A_1) }\varphi_{\al}(\al(a)^{-1}U_\al)
\end{equation}
où
\begin{itemize}
\item $Y$ parcourt les éléments nilpotents de $q_2'(F)\cap \mgo_R(F)$
\item $U$ décrit les éléments de  $\ngo_R(F)$ tels que $Y+U\notin I_R^G(Y)$.
\end{itemize}
Pour tout élément nilpotent $Y\in \mgo_R(F)$, la condition  $Y+U\notin I_R^G(Y)$ définit un fermé propre de $\ngo_R(F)$ décrit par une famille de polynômes dont le degré est majoré indépendamment de $Y$ (il n'y a qu'un nombre fini d'orbites nilpotentes possibles).  Le lemme \ref{lem:zero-de-py} nous donne une  constante $c_R>0$ telle que \eqref{eq:majI} se majore par 
\begin{equation}
  \label{eq:majII}
 c_R\cdot  \exp(-\bg \al,H_{Q_1}(a)\bd)\cdot  \exp(-\bg 2\rho_{Q_1}^R,H_{Q_1}(a)\bd )\sum_{Y \in \qgo_2'(F)\cap \mgo_R(F)} \varphi_1(Y_1)\prod_{\al\in  \Sigma(Q_2\cap M_R,A_1) }\varphi_{\al}(\al(a)^{-1}Y_\al)  
\end{equation}
pour $a\in A_1^{G,\infty}(T,T_1)$.

Pour terminer, on observe que pour $a\in A_1^{G,\infty}(Q,T)$, l'expression 
$$\exp(-\bg 2\rho_{Q_1}^R,H_{Q_1}(a)\bd )\sum_{Y \in \qgo_1(F)\cap \mgo_R(F)} \varphi_1(Y_1)\prod_{\al\in  \Sigma(N_1 \cap M_R,A_1) }\varphi_{\al}(\al(a)^{-1}Y_\al) $$
est bornée et l'expression
$$\sum_{Y \in \qgo_2'(F)\cap \bar{\ngo}_{1}(F)\cap \mgo_R(F) } \prod_{\al\in  \Sigma( \bar{N}_1\cap M_R,A_1) }\varphi_{\al}(\al(a)^{-1}Y_\al),$$
où $\bar{N}_1$ est le radical unipotent du sous-groupe parabolique opposé à $Q_1$, se majore par
$$ C\cdot \prod_{\al \in \Delta_{1}^2-\Delta_{1}^{Q}}\exp(-c_\al \bg \al, H_{Q_1}(a)\bd)$$
où $C$ est une constante $>0$ et les coefficients $c_\al$ sont aussi grands que l'on veut. Pour le voir, on observe que pour tout $\al \in \Delta_{1}$, l'ensemble des $\bg \al, H_1(a)\bd$ est borné inférieurement quand $a\in A_1^{G,\infty}(Q,T)$ : cela résulte de la condition 1 si $\al\in \Delta_1^Q$ et de l'inégalité \eqref{eq:minor-beta} sinon. Puis on utilise les majorations élémentaires suivantes (la première repose sur la formule de Poisson). Soit $t_0 >0$. Pour toute fonction de Schwartz $f$ sur $\RR^d$ et pour tout $n\in \NN$ il existe $N\geq n$ et  $C>0$ tels que les majorations suivantes soient valables  pour tout $t\geq t_0$
 \begin{enumerate}
 \item $\displaystyle t^{-d} \sum_{X\in \ZZ^d} f(t^{-1}X ) \leq C.$
 \item $\displaystyle \sum_{X\in \ZZ^d-\{0\}} f(t X ) \leq C\cdot t^{-n}.$
 \item $\displaystyle \sum_{X\in \ZZ^d} f(t X ) \leq C.$
 \end{enumerate}
Cela termine la démonstration.
\end{preuve}
\end{paragr}

\begin{paragr}[Majoration de \eqref{eq:second-int}.] --- Pour tout couple de sous-groupes  paraboliques $(P_1,P_2)$ adjacents dans $\pc(M)$, tout $T\in a_{B_0}$ et tout  $g\in A_G(\AAA)G_X(\AAA)\back G(\AAA)$, il existe des réels $\gamma_{P_1,P_2}(T)$ et $\gamma_{P_1,P_2}(g)$ uniquement déterminés par la condition
$$T_{P_1}-T_{P_2}= \gamma_{P_1,P_2}(T) \al^\vee$$
et 
$$-R_{P_1}(g)+R_{P_2}(g)=  \gamma_{P_1,P_2}(g) \al^\vee$$
où $\al^\vee$ est l'unique élément de $\Delta_{P_1}^\vee\cap (-\Delta_{P_2}^\vee)$.
Le calcul de  $\gamma_{P_1,P_2}(g)$ a été donné au lemme \ref{lem:Y}. On laisse au lecteur le soin de préciser la constante $\gamma_{P_1,P_2}(T) $ : il nous suffit de savoir que c'est une certaine somme non vide de termes $\bg \beta, T \bd$ où $\beta$ est une racine de $T_0$ dans $N_{B_0}$. On supposera que $T$ vérifie l'inégalité $\bg \al,T\bd \geq \eps' \|T\|$ pour tout $\al\in \Delta_{B_0}$. En particulier, on retiendra qu'on a $\gamma_{P_1,P_2}(T) \geq  \eps' \|T\|.$ Si la famille orthogonale $(T_P-R_P(g))_{P\in \pc(M)}$ n'est pas positive, il existe un couple de sous-groupes  paraboliques $(P_1,P_2)$ adjacents dans $\pc(M)$ tels que $\gamma_{P_1,P_2}(g) < -\gamma_{P_1,P_2}(T).$   En particulier, on doit avoir
  \begin{equation}
    \label{eq:ineg-eps}
    \gamma_{P_1,P_2}(g)<  -\eps' \|T\|.
  \end{equation}
  Il s'ensuit qu'on a la majoration suivante
$$1-I^T(g) \leq \sum_{(P_1,P_2)}  I_{P_1,P_2}^T(g)$$
où $(P_1,P_2)$ parcourt les couples d'éléments adjacents de $\pc(M)$ et où $I_{P_1,P_2}(g)$ est la fonction caractéristique des éléments $g\in  A_G(\AAA)G_X(\AAA)\back G(\AAA)$ qui satisfont l'inégalité \eqref{eq:ineg-eps}. Dans la suite, on fixe un tel couple $(P_1,P_2)$ et on va majorer

\begin{equation}
   \label{eq:second-int2}
    \int_{G_X(F)\back G(\AAA)^1} (1-F^{G}(g,T)) \cdot|\sigma^T(g)|\cdot I_{P_1,P_2}^T(g)  \cdot |f(g^{-1}Xg) |  \, dg.
 \end{equation}
Tout d'abord, on majore brutalement $1-F^{G}(g,T)$ par $1$. On est donc ramené à majorer
\begin{equation}
   \label{eq:second-int3}
    \int_{G_X(\AAA)\back G(\AAA)}  v(g,T)\cdot I_{P_1,P_2}^T(g)  \cdot |f(g^{-1}Xg) |  \, dg.
 \end{equation}
où l'on a introduit 
$$v(g,T)=\int_{G_X(F)\back G_X(\AAA)^1  }|\sigma^T(hg)|\, dh.$$
Pour obtenir la majorant du théorème \ref{thm:asymp}, il suffit d'obtenir un majorant analogue pour le carré de \eqref{eq:second-int3}. Par l'inégalité de Cauchy-Schwartz, on est ramené à majorer les deux expressions suivantes 
\begin{equation}
   \label{eq:second-int4}
  \int_{G_X(\AAA)\back G(\AAA)}  v(g,T)^2  |f(g^{-1}Xg) |  \, dg.
 \end{equation}
et
\begin{equation}
   \label{eq:second-int5}
  \int_{G_X(\AAA)\back G(\AAA)}   I_{P_1,P_2}^T(g)   |f(g^{-1}Xg) |  \, dg.
 \end{equation}

Comme dans la démonstration du théorème \ref{thm:nilp-globale}, on montre qu'il existe une constante  $C>0$ tel que
  $$v(g,T)^2\leq C\cdot(\|T\|^{2k}+ \sum_{(P,P')\in \pc(M)^{\mathrm{adj}}}   \|R_P(g)-R_{P'}(g) \|^{2k} )$$
où $k=\dim(a_M^G)$. En utilisant les résultats de convergence du théorème \ref{thm:cv} et ceux qui apparaissent dans la démonstration du théorème \ref{thm:cv-pondere},  on obtient la majoration suivante pour \eqref{eq:second-int4}
\begin{equation}
   \label{eq:second-int6}
  \int_{G_X(\AAA)\back G(\AAA)}  v(g,T)^2  |f(g^{-1}Xg) |  \, dg\leq  c_1  \|T\|^{2k}+c_2
\end{equation}
où $c_1$ et $c_2$ sont indépendants de $T$.

Majorons ensuite \eqref{eq:second-int5}. On reprend les notations de la démonstration du théorème \ref{thm:cv-pondere}. En utilisant les lemmes \ref{lem:Y} et \ref{lem:alaRao}, on majore  \eqref{eq:second-int5} par  (à la constante $c$ près définie en \eqref{eq:cestc}),
\begin{equation}
  \label{eq:second-int7}
  \int_{I_{1,3}(\AAA)} \int_{\ngo_{\tilde{Q}}(\AAA)} \int_{K} f(k^{-1}(x+U)k) i^T(x)   |\det(x)|^{r_1+r_2} \, dk\, dU \, dx.
\end{equation}
où l'on introduit $i^T$  la fonction caractéristique des éléments $x\in I_{1,3}(\AAA)$ qui vérifient
$$\log|\det(x)| \leq -\eps'\|T\|.$$
Soit $0<\eta < 1$. On utilise l'inégalité  évidente 
$$ i^T_{1,3}(x)   \leq |\det(x)|^{-\eta} \exp(-\eta \eps' \|T\|)$$
pour majorer \eqref{eq:second-int7} par 

\begin{equation}
  \label{eq:second-int8}
 \exp(-\eta \eps' \|T\|) \int_{I_{1,3}(\AAA)} \int_{\ngo_{\tilde{Q}}(\AAA)} \int_{K} f(k^{-1}(x+U)k)  |\det(x)|^{r_1+r_2 -\eta} \, dk\, dU \, dx.
\end{equation}
où l'intégrale ci-dessus convergente puisque $r_1+r_2-\eta >r_1$ (rappelons qu'on a $r_2\geq 1$ lorsque $X$ est simple). Cela termine la démonstration.
\end{paragr}

\section{Une forme précisée du développement d'Arthur}\label{sec:dvpt-Arthur}

\begin{paragr} Les hypothèses et les notations sont celles des sections \ref{sec:IOPglob} et \ref{sec:asymp}. En particulier $X$ est supposé simple.
\end{paragr}

\begin{paragr}
Soit $\of$ l'orbite de $X$. On a défini en \eqref{eq:JoT} une distribution $J_{\of}^T$ polynomiale en $T$. D'après le théorème \ref{thm:asymp}, ce polynôme est asymptotique à l'intégrale orbitale tronquée \eqref{eq:IO-tronquee}. D'après Arthur, cf. \cite{ar_unipvar}, la valeur $J_\of(f)$ de ce polynôme en $T=0$ est la contribution de l'orbite $\of$ au développement géométrique de la formule des traces (ou plutôt son analogue \cite{PH1} sur les algèbres de Lie). Dans \cite{ar_unipvar}, Arthur donne un développement de $J_\of(f)$ en terme d'intégrales orbitales pondérées unipotentes locales mais avec des coefficients non spécifiés. Dans cette section, nous allons pour l'orbite $\of$ retrouver le développement d'Arthur avec en prime des formules intégrales pour les coefficients.
\end{paragr}

\begin{paragr}   Soit $\lc_X$ l'ensemble des couples $(L,\of)$ formés d'un sous-groupe $L\in \lc(T_0)$ et d'une orbite $L$-nilpotente $\of\subset \lgo$ tels que $X\in I_L^G(\of)$. Le groupe $W$ agit sur $\lc_X$ par $w\cdot (L,\of)=(wLw^{-1},w\of w^{-1})$. Soit $\lc_X/W$ le quotient. On a une application 
$$\lc(M) \to \lc_X$$
donnée par $L\mapsto (L,I_M^L(0))$. Par composition avec la projection canonique $\lc_X \to \lc_X/W$, on obtient une application
\begin{equation}
  \label{eq:bij-cas-simple}
  \lc(M) \to \lc_X/W.
\end{equation}

  \begin{lemme}\label{lem:lcX}Soit $X$ simple.
    \begin{enumerate}
    \item L'application $\lc(M) \to \lc_X/W$ est bijective. 
    \item L'orbite de $(L,\of)$ sous l'action de $W$ est de cardinal $|W|/|W^L|$.
    \end{enumerate}
  \end{lemme}
  
  \begin{preuve}
    Traitons l'assertion 1. La surjectivité a été prouvée dans la démonstration du théorème \ref{thm:egalite-IOP}. Soit $L_1$ et $L_2$ dans $\lc(M)$ et $w\in W$ tels que $L_2=wL_1w^{-1}$ et $I_{M}^{L_2}(0)=wI_{M}^{L_1}(0)w^{-1}$. Cette dernière induite n'est autre que $I_{wMw^{-1}}^{L_2}$. Il s'ensuit que $M$ et $wMw^{-1}$ sont conjugués sous $W^{L_2}$. Donc $w\in W^{L_2}\Norm_W(M)$. En particulier, si $X$ est simple on a $\Norm_W(M)=W^M$ et donc $w\in W^{L_2}$ d'où $L_1=L_2$.

    Prouvons ensuite l'assertion 2. On cherche le stabilisateur de $(L,\of)$. Par la surjectivité dans l'assertion 1, on peut supposer $L\in \lc(M)$ et $\of=I_M^L(0)$. Ce qui précède montre que le stabilisateur est $W^L$.
  \end{preuve}
\end{paragr}

\begin{paragr}[Ensemble $S$.] --- Soit $S\subset V$ un ensemble fini de places qui contient les places archimédiennes. On note $\AAA_S=\prod_{v\in S} F_v$ et $\AAA^S$ l'anneau des adèles hors $S$ de sorte qu'on a  $\AAA=\AAA_S\times \AAA^S$.   
\end{paragr}

\begin{paragr}[Intégrales orbitales pondérées semi-locales d'Arthur.] --- \label{S:semiloc} Dans ce paragraphe, on introduit des énoncés et des constructions qui généralisent ceux de la section \ref{sec:IOPloc} à $\AAA_S$. Comme il suffit soit de reprendre les démonstrations de la section  \ref{sec:IOPloc} soit de faire appel à la théorie du scindage des $(G,M)$-familles (cf. \cite{inv_loc} section 9), les détails sont laissés au lecteur. Soit $f_S\in \Sc(\ggo(\AAA_S)$. On introduit tout d'abord pour tout $L\in \lc(M)$ et tout $Q\in \fc(L)$ l'intégrale orbitale pondérée (qui est absolument convergente)
$$J_{L,X}^Q(f_S)= \int_{G_X(\AAA_S)\back G(\AAA_S)}f_S(\Ad(g^{-1})X) v_{L,X}^Q(g)\, dg.$$

Les intégrales nilpotentes d'Arthur sont introduites comme dans le cas local, cf. §\ref{S:Arthur-loc}. Soit $L\in \lc(T_0)$ et $\of\subset \lgo$ une $L$-orbite nilpotente. On pose
$$
J_L^G(\of,f)=\int_{\ngo_P(\AAA_S)} \int_{K_S} f(\Ad(k^{-1})U) \cdot w_L(U^L,U_Q)  \, dU dk.$$
C'est l'analogue pour $\AAA_S$ de la définition \eqref{eq:Jof}. Les notations sont celles de la section
  \ref{S:Arthur-loc} à la précision suivante près : on pose $K_S=\prod_{v\in S}K_v$ et le poids   $w_L(U^L,U_Q)$ est celui obtenu à partir de la $(G,L)$-famille produit
$$w_{Q'}(\la,U^L,U_Q)=\prod_{v\in S} w_{Q'}(\la,U^L_v,U_{Q,v})$$
pour tout $Q'\in \pc(L)$, l'indice $v$ désignant une composante sur $F_v$. Comme dans le théorème \ref{thm:egalite-IOP}, cette intégrale est absolument convergente et ne dépend pas des choix de $P$ et $Q$ qui interviennent dans sa définition. De plus, si $L_1\in \lc(M)$ correspond à la classe de $(L,\of)$ par la bijection \eqref{eq:bij-cas-simple}, alors on a
\begin{equation}
  \label{eq:IO-Arthur-S}
  J_{L_1,X}^G(f_S)= c_{X,S}\cdot J_{L}^G(\of,f_S),
\end{equation}
où $c_{X,S}=\prod_{v\in S}c_{X,v}$ est le produit des constantes locales qui apparaissent dans le lemme \ref{lem:Howe}.
\end{paragr}

\begin{paragr}[Intégrales orbitales pondérées hors $S$ pour la fonction unité.] --- Pour tout $L\in \lc(M)$, soit  $\mathbf{1}^S_L\in \Sc(\lgo(\AAA^S))$ la fonction caractéristique de $\lgo(\oc^S)$ où  $\oc^S=\prod_{v\notin S} \oc_v$. Si $L=G$, on note $\mathbf{1}^S=\mathbf{1}^S_G$. On introduit les intégrales orbitales suivantes pour tous $L\in \lc(M)$ et $Q\in \pc(L)$
  \begin{equation}
    \label{eq:JMXhorsS}
     J_{M,X}^{Q}(\mathbf{1}^S)=\int_{G_X(\AAA^S)\back G(\AAA^S)}\mathbf{1}^S(\Ad(g^{-1})X) v_{M,X}^Q(g)\, dg
  \end{equation}
et pour tout $P\in \pc^L(M)$
\begin{equation}
  \label{eq:JML-horsS}
  J_M^L((0),\mathbf{1}^S_L)=\int_{\ngo_{P}(\AAA_S)} \mathbf{1}^S_L(U) w^L_M(U) \,dU
\end{equation}
où le poids  $w^L_M(U)$ est le poids d'Arthur relatif à $L$ défini à l'aide de la $(L,M)$-famille $(w^L_P(\la,U))_{P\in \pc^L(M)}$, cf. \eqref{eq:wPbis}.
  
\begin{lemme}\label{lem:S-cv-inv}
Supposons $X$ simple. Les intégrales \eqref{eq:JMXhorsS} et \eqref{eq:JML-horsS} sont absolument convergentes. De plus, pour tout $Q\in \pc(L)$, on a 
\begin{equation}
  \label{eq:uneegalite}
      J_{M,X}^{Q}(\mathbf{1}^S)=c_X^S \cdot  J_M^L((0),\mathbf{1}^S_L)
    \end{equation}
    où $c_{X}^S=\prod_{v\notin S}c_{X,v}$ est le produit convergent des constantes locales qui apparaissent dans le lemme \ref{lem:Howe}.
  \end{lemme}

  \begin{preuve}
    La convergence de  \eqref{eq:JMXhorsS} se prouve comme le théorème \ref{thm:cv-pondere}. L'égalité \eqref{eq:uneegalite} ainsi que la convergence de  \eqref{eq:JML-horsS} résulte d'une variante globale hors $S$ du lemme \ref{lem:descente}. Le point clef est la convergence du produit   $c_{X}^S$ sous la condition $X$ simple. 
  \end{preuve}

Plus généralement, on définit pour tous sous-groupes de Levi semi-standard $L_1\subset L_2$ conjugués sous $W$ à $M\subset L$ une intégrale orbitale convergente par la formule
\begin{equation}
  \label{eq:JML-horsS-bis}
  J_{L_1}^{L_2}((0),\mathbf{1}^S_{L_2})=\int_{\ngo_{P}(\AAA_S)} \mathbf{1}^S_L(U) w^{L_2}_{L_1}(U) \,dU
\end{equation}
où $P\in \pc^{L_2}(L_1)$. On peut montrer qu'on a 
\begin{equation}
  \label{eq:JLsousW}
   J_{L_1}^{L_2}((0),\mathbf{1}^S_{L_2})=J_{M}^{L}((0),\mathbf{1}^S_{L}).
 \end{equation}

\end{paragr}

\begin{paragr}[Le développement d'Arthur.] --- Le théorème suivant est une généralisation aux fonctions de Schwartz du corollaire 8.4 de \cite{ar_unipvar}. Non seulement notre démonstration diffère de celle d'Arthur mais, pour les orbites considérées, on obtient aussi la formule intégrale \eqref{eq:aS} pour les coefficients $a^L(S,\of)$ (alors qu'Arthur n'affirmait que leur existence). 

  \begin{theoreme}\label{thm:dvpt} On suppose $X$ simple.
Soit $\of_X$ l'orbite de $X$  sous $G$. Pour toute fonction $f_S\in \Sc(\ggo(\AAA_S))$, on a
\begin{equation}
  \label{eq:dvpt}
  J_{\of_X}(f_S\otimes \mathbf{1}^S)=\sum_{(L,\of)\in \lc_X} \frac{|W^L|}{|W|} a^L(S,\of) J_L^G(\of, f_S),
\end{equation}
où, pour tout couple $(L,\of)\in \lc_X$, on pose
\begin{equation}
  \label{eq:aS}
   a^L(S,\of)=\vol(L_1(F)\back L_1(\AAA)^1) \cdot J_{L_1}^{L}((0),\mathbf{1}^S_{L})
 \end{equation}
où $L_1\in \lc^L(T_0)$ est un élément quelconque qui vérifie $I_{L_1}^{L}(0)=\of$. Le coefficient $ a^L(S,\of)$ est indépendant du choix de $L_1$.
  \end{theoreme}

  \begin{remarque}
    Au passage près de l'algèbre de Lie au groupe, on peut se demander si les coefficients obtenus dans \cite{ar_unipvar} sont les mêmes que ceux d'Arthur. C'est bien le cas et voici une façon de procéder. On montre  d'abord que le théorème \ref{thm:asymp} et la proposition \ref{thm:nilp-globale} valent encore si l'on remplace $G$ par un sous-groupe de Levi $L'$ et $X$ par un élément analogue $X'$ tel que $X\in I_{L'}^G(X')$. Alors $X'$ est encore simple et le théorème \ref{thm:dvpt} vaut  pour tout tel $L'$ (en lieu et place de $G$). Sous cette dernière condition, les coefficients  $a^L(S,\of)$ sont uniquement déterminés comme le montre Arthur dans \cite{ar_unipvar}. Par exemple, la dernière étape consiste à prouver que $a^G(S,\of_X)$ est uniquement déterminé par \eqref{eq:dvpt} et la donnée des   $a^L(S,\of)$ pour $(L,\of)\in \lc_X'$  où $\lc'_X=\lc_X-\{(G,\of_X)\}$. Or \eqref{eq:dvpt} se réécrit
    $$J_{\of_X}(f_S\otimes \mathbf{1}^S)-\sum_{(L,\of)\in \lc_X} \frac{|W^L|}{|W|} a^L(S,\of) J_L^G(\of, f_S)=a^G(S,\of) J_G^G(\of_X,f_S)$$
ce qui détermine $a^G(S,\of)$ vu que $J_G^G(\of_X,f_S)$ est une intégrale orbitale ordinaire donc non nulle en tant que distribution.
  \end{remarque}

  \begin{preuve}
Le point de départ est la proposition \ref{thm:nilp-globale} qui donne l'expression intégrale
\begin{equation}
  \label{eq:exp-integrale}
  J_{\of_X}(f_S\otimes \mathbf{1}^S)=\vol(G_X(F)\back G_X(\AAA)^1) \cdot  \int_{G_X(\AAA)\back G(\AAA) }  (f_S\otimes \mathbf{1}^S) (\Ad(g^{-1})X) \, v_M^G(g) \, dg.
\end{equation}

On invoque ensuite la formule suivante de scindage du poids (conséquence de formules générales d'Arthur, cf. \cite{inv_loc} corollaire 7.4)    
$$v_{M}^G(g)=\sum_{(L_1,L_2)\in \lc(M)^2} d_{M}^{G}(L_1,L_2)   v_{M,X}^{Q_1}(g_S) v_{M,X}^{Q_2}(g^S) ;$$
on a écrit $g=(g_S,g^S)\in G(\AAA)=G(\AAA_S)\times G(\AAA^S)$ ;  un certain choix que l'on ne précise pas car il ne jouera pas de rôle ici définit une application 
$$\{(L_1,L_2)\in  \lc(M)^2\mid d_{M}^{G}(L_1,L_2)\not=0 \}\to \fc(M)^2$$
notée $(L_1,L_2)\mapsto (Q_1,Q_2)$ telle que $M_{Q_i}=L_i$ pour $i=1,2$.

On en déduit que $J_{\of_X}(f_S\otimes \mathbf{1}^S)$ vaut 
$$  \vol(G_X(F)\back G_X(\AAA)^1) \times \sum_{(L_1,L_2)\in \lc(M)^2} d_{M}^{G}(L_1,L_2)  J_{M,X}^{Q_2}(\mathbf{1}^S)  \int_{G_X(\AAA_S)\back G(\AAA_S) }  f_S (\Ad(g^{-1})X) \, v_M^{Q_1}(g) \, dg$$
D'après le lemme \ref{lem:S-cv-inv}, on a    
$$J_{M,X}^{Q_2}(\mathbf{1}^S)=c_X^S\cdot J_{M}^{L_2}((0), \mathbf{1}^S_{L_2}).$$ 
En particulier, cette expression ne dépend que de $L_2$. De plus, à $L_2$ fixé, on a 
$$\sum_{L_1\in \lc(M)} d_{M}^{G}(L_1,L_2)   v_{M,X}^{Q_1}(g_S)=v_{L_2,X}^G(g_S)$$
(cf. \cite{inv_loc} section 8). En partant de \eqref{eq:exp-integrale} et en tenant successivement compte de \eqref{eq:IO-Arthur-S} et du lemme \ref{lem:calcul} ci-dessous, on obtient
\begin{eqnarray*}
  J_{\of_X}(f_S\otimes \mathbf{1}^S)&=&  c_X^S\cdot \vol(G_X(F)\back G_X(\AAA)^1) \cdot \sum_{L\in \lc(M)} J_{M}^{L}((0), \mathbf{1}^S_{L})  J_{L,X}^G(f_S)\\
&=& c_X\cdot\vol(G_X(F)\back G_X(\AAA)^1)\cdot \sum_{L\in \lc(M)} J_{M}^{L}((0), \mathbf{1}^S_{L})  J_{L}^G(I_M^L(0),f_S)\\
&=& \sum_{L\in \lc(M)} a^L(S,I_M^L(0)) J_{L}^G(I_M^L(0),f_S),
\end{eqnarray*}
où $c_X$ est la constante globale \eqref{eq:formule-c-globale}. Soit $(L,\of)\in \lc_X$. Soit $L_1\in \lc^L(T_0)$ tel que $I_{L_1}^L(0)=\of$. L'expression
$$\vol(L_1(F)\back L_1(\AAA)^1) \cdot J_{L_1}^{L}((0),\mathbf{1}^S_{L})$$
est indépendante du choix de $L_1$. En effet, tout autre choix est un conjugué sous $W^L$ de $L_1$ et on a pour tout $w\in L_1$
$$ J_{L_1^w}^{L}((0),\mathbf{1}^S_{L}) =J_{L_1}^{L}((0),\mathbf{1}^S_{L}).$$

Il existe $w\in W$ et un unique $L'\in \lc(M)$ tel que $(L,\of)=w\cdot (L',I_M^{L'}(0))$. Quitte à translater $w$ par un élément de $W^L$, on peut supposer que $L_1=wMw^{-1}$. Il résulte alors de \eqref{eq:JLsousW} qu'on a
\begin{equation}
  \label{eq:aL=aL'}
  a^L(S,\of)=a^{L'}(S,I_M^{L'}(0)).
\end{equation}
On a alors en utilisant le lemme \ref{lem:lcX}

\begin{eqnarray*}
  J_{\of_X}(f_S\otimes \mathbf{1}^S)&=& \sum_{(L,\of)\in \lc_X/W} a^L(S,\of) J_{L}^G(\of    ,f_S) \\
&=&  \sum_{(L,\of)\in \lc_X}  \frac{|W^L|}{|W|} a^L(S,\of) J_{L}^G(\of,f_S)
\end{eqnarray*}
\end{preuve}

\begin{lemme}\label{lem:calcul}
On a 
\begin{equation}
  \label{eq:verif}
  c_X\cdot \vol(G_X(F)\back G_X(\AAA)^1)=\vol(M(F)\back M(\AAA)^1).
\end{equation}
\end{lemme}

\begin{preuve}
  On a 
$$M=GL(d_1+\ldots+d_r)\times \ldots \times GL(d_{r-1}+d_r)\times GL(d_r).$$
\emph{A priori}, les mesures sur  $G_X(\AAA)^1$ et $ M(\AAA)^1$ sont celles fixées au §\ref{S:Haar-G1}. Ces mesures  dépendent toutes deux du choix de la mesure de Haar sur $a_M$. En revanche l'égalité \eqref{eq:verif} n'en dépend pas. Les déterminants des blocs $GL$ permettent d'identifier $a_M$ à $\RR^r$.  Il est plus commode pour la démonstration de supposer que la mesure de Haar sur $a_M$ correspond à la mesure euclidienne sur $\RR^r$ par cette identification.

Un calcul classique (cf. \cite{Weil-volume} ou \cite{Lan-Boulder}) montre qu'on a alors \emph{pour cette nouvelle mesure}
\begin{eqnarray*}
  \vol(M(F)\back M(\AAA)^1)&=& d_F^{\dim(M)/2}\cdot \prod_{i=1}^r   Z_{d_i+\ldots+d_r}^*(d_i+\ldots+d_r)
\end{eqnarray*}
avec les notations suivantes :
\begin{itemize}
\item $d_F$ est le discriminant du corps $F$ ;
\item  pour tout entier $n\geq 1$ on pose $Z_n^*(s)=(s-n+1)Z_n(s)$ où la fonction $Z_n$ est celle introduite en \eqref{eq:Zd-globale}.
\end{itemize}

Le centralisateur $G_X$ a une décomposition de Levi $G_X=M_XN_X$ (cf.§ \ref{S:Levi-cent}). Pour nos choix de mesures, on a 
$$\vol(G_X(F)\back G_X(\AAA)^1) =\vol(M_X(F)\back M_X(\AAA)^1)\cdot \vol(N_X(F)\back N_X(\AAA)),$$
où la mesure de Haar sur $M_X(\AAA)^1$ est celle compatible à la suite exacte (où $P\in \pc(M)$)
$$1\to M_X(\AAA)^1\to M_X(\AAA)\to_{R_P} a_M \to 1$$
et à la mesure fixée sur  $M_X(\AAA)$ et à la nouvelle mesure de Haar sur $a_M$. On a
\begin{eqnarray*}
  \vol(N_X(F)\back N_X(\AAA))&=& (d_F)^{\dim(N_X)/2}
\end{eqnarray*}
Le facteur réductif $M_X$ s'identifie à $GL(d_1)\times\ldots \times GL(d_r)$. On obtient avec nos (nouveaux) choix de mesure
\begin{eqnarray*}
  \vol(M_X(F)\back M_X(\AAA)^1)&=& d_F^{\dim(M_X)/2}\cdot \prod_{i=1}^r   Z_{d_i}^*(d_i).
\end{eqnarray*}
On a donc
$$\vol(G_X(F)\back G_X(\AAA)^1) = d_F^{\dim(G_X)/2}\cdot \prod_{i=1}^r   Z_{d_i}^*(d_i).
$$
L'égalité \eqref{eq:verif} est alors une conséquence d'une part de la formule bien connue $\dim(G_X)=\dim(M)$ et d'autre part du calcul suivant
\begin{eqnarray*}
   \prod_{i=1}^r   Z_{d_i}^*(d_i) \cdot c_X &=&   \prod_{i=1}^r   Z_{d_i}^*(d_i) \prod_{j=1}^r \prod_{i=1}^{j-1} Z_{d_j}(d_i+\ldots+d_j) \\
&=&\prod_{i=1}^r   Z_{d_i}^*(d_i) \cdot \prod_{i=1}^{r-1}  \prod_{k=d_i+1}^{d_i+\ldots+d_r} Z_1(k)\\
&=& \prod_{i=1}^r   Z_{d_i+\ldots+d_r}^*(d_i+\ldots+d_r).
\end{eqnarray*}
  \end{preuve}

\end{paragr}

\begin{paragr}[Deux applications.] --- Tout d'abord on retrouve la seule formule générale qu'Arthur donne.

  \begin{corollaire}\label{cor:calcul}
    On a 
$$a^M(S,(0))=\vol(M(F)\back M(\AAA)^1).$$
  \end{corollaire}
 
  \begin{preuve}
  La formule résulte de la formule évidente $J_M^M((0), \mathbf{1}^S_M)= \mathbf{1}^S_M(0)=1$.
  \end{preuve}

  \begin{corollaire}\label{cor:maj}
Soit $L$ un sous-groupe de Levi contenant $M$. Il existe une constante $C>0$ indépendante de $F$ telle que pour  tout ensemble $S$ fini de places contenant les places archimédiennes et tout $Q\in \pc(L)$ on a  
$$
a^L(S,I_M^L(0)) \leq C\cdot \vol(M(F)\back M(\AAA)^1)\cdot \sup_{(P_1,P_2)\in \pc(M)^{Q,\text{adj}}} (-1)^k  \frac{\frac{d^k}{ds^k}(Z_{r_1}^S)(r_1+r_2)}{Z_{r_1}^S(r_1+r_2)}
$$
où
\begin{itemize}
\item $\pc(M)^{Q,\text{adj}}$ est l'ensemble des couples $(P_1,P_2)\in \pc^Q(M)^2$ tels que $P_1$ et $P_2$ sont adjacents ; à un tel couple, on associe les entiers $r_1>0$ et $r_2>0$ (cf. \eqref{eq:r1} et  \eqref{eq:r2}) ;
\item $k=\dim(a_M^L)$ ;
\item $Z_r^S$ est la fonction zêta $Z_r$ dont on ne garde que les facteurs hors $S$.
\end{itemize}
\end{corollaire}

\begin{remarques}
Vu \eqref{eq:aL=aL'}, on en déduit une majoration de $a^L(S,\of)$ pour tout $(L,\of)\in \lc_X$. Il serait intéressant de comparer cette majoration à celle obtenue par Matz dans \cite{Matz-bounds}.  
\end{remarques}

  \begin{preuve}
 Il s'agit de démontrer qu'il existe une constante $C>0$ indépendante de $F$ telle que pour  tout $Q\in \pc(L)$ on ait 

$$
J_M^L(0,\mathbf{1}_L^S)\leq C\cdot \sup_{(P_1,P_2)\in \pc(M)^{Q,\text{adj}}} (-1)^k   \frac{\frac{d^k}{ds^k}(Z_{r_1}^S)(r_1+r_2)}{Z_{r_1}^S(r_1+r_2)}.
$$
Soit $Q\in \pc(L)$.
On peut encore écrire à l'aide de \eqref{eq:uneegalite}
$$J_M^L(0,\mathbf{1}_L^S)= J_{M,X}^{Q}(\mathbf{1}^S)\cdot J_{M,X}^{P}(\mathbf{1}^S)^{-1}$$
pour tout $P\in \pc(M)$. On observe (cf. démonstration de la proposition \ref{prop:cv-locale}) qu'il existe $C_1>0$ indépendant de $F$ tel que pour tout $g\in G(\AAA^S)$ on ait (cf. \eqref{eq:une-maj-adj})
$$|v_{M,X}^{Q}(g)|\leq C_1\cdot \sup_{(P_1,P_2)\in \pc(M)^{Q,\text{adj}}} \|R_{P_1}(g)-R_{P_2}(g)\|^k.$$

En procédant comme dans la démonstration de la proposition \ref{prop:cv-locale} dont on reprend sans plus de commentaire les notations, on montre qu'il existe $C_2>0$ indépendant de $F$ telle que 

$$ J_{M,X}^{Q}(\mathbf{1}^S)\cdot J_{M,X}^{P}(\mathbf{1}^S)^{-1}\leq C_2\cdot  \sup_{(P_1,P_2)\in \pc(M)^{Q,\text{adj}}}  J_{P_1,P_2}^Q(f)(k,\mathbf{1}^S)  (J_{P_1,P_2}^Q(f)(0,\mathbf{1}^S)^{-1}.$$
où l'on pose pour tout entier $k$
\begin{eqnarray*}
   J_{P_1,P_2}^Q(f)(k,\mathbf{1}^S)&=& \int_{I_{1,3}(\AAA^S)} \int_{\ngo_{\tilde{Q}}(\AAA^S)} \int_{K^S} \mathbf{1}^S( k^{-1}(x+U)k)|\log|\det(x)| |^k |\det(x)|^{r_1+r_2} \, dk\, dU \, dx\\
&=&  (-1)^k \int_{I_{1,3}(\AAA^S)}  \mathbf{1}^S( x)(\log|\det(x)| )^k |\det(x)|^{r_1+r_2}  \, dx.
\end{eqnarray*}
Cela conclut puisqu'on a 
$$ J_{P_1,P_2}^Q(f)(0,\mathbf{1}^S)=Z_{r_1}^S(r_1+r_2)$$
et
$$ J_{P_1,P_2}^Q(f)(k,\mathbf{1}^S)= (-1)^k\frac{d^k}{ds^k}(Z_{r_1}^S)(r_1+r_2).$$

  \end{preuve}
\end{paragr}

\bibliography{bibiog}
\bibliographystyle{plain}

\begin{flushleft}
Pierre-Henri Chaudouard \\
Université Paris Diderot (Paris 7) et Institut Universitaire de France\\
 Institut de Mathématiques de Jussieu-Paris Rive Gauche \\
 UMR 7586 \\
 Bâtiment Sophie Germain \\
 Case 7012 \\
 F-75205 PARIS Cedex 13 \\
 France
\medskip

Adresse électronique :\\
Pierre-Henri.Chaudouard@imj-prg.fr \\
\end{flushleft}

\end{document}